\documentclass[preprint,3p,12pt]{elsarticle}
\usepackage{mathrsfs}
\usepackage{amsmath}
\usepackage{stmaryrd}
\usepackage{bbding}
\usepackage{dcolumn}
\usepackage{graphicx}
\usepackage{amsfonts}
\usepackage{amssymb}
\usepackage{psfrag}
\usepackage{wrapfig}
\usepackage{subfigure}
\usepackage{makeidx}
\usepackage{bm}
\usepackage{epsf}
\usepackage{epsfig}
\usepackage{setspace}
\usepackage{graphicx}
\usepackage{epstopdf}
\usepackage{psfrag}
\usepackage{subfigure}
\usepackage{color}

\begin{document}
\title{A well-balanced gas-kinetic scheme with adaptive mesh refinement for shallow water equations}

\author[HKUST1]{Gaocheng Liu}
\ead{gliuau@connect.ust.hk}

\author[HKUST1]{Fengxiang Zhao}
\ead{fzhaoac@connect.ust.hk}

\author[HKUST1,HKUST2]{Jianping Gan}
\ead{magan@ust.hk}

\author[HKUST1,HKUST2,HKUST3]{Kun Xu\corref{cor}}
\ead{makxu@ust.hk}

\address[HKUST1]{Department of Mathematics, Hong Kong University of Science and Technology, Clear Water Bay, Kowloon, Hong Kong}
\address[HKUST2]{Center for Ocean Research in Hong Kong and Macau (CORE), Hong Kong University of Science and Technology, Clear Water Bay, Kowloon, Hong Kong}
\address[HKUST3]{Shenzhen Research Institute, Hong Kong University of Science and Technology, Shenzhen, China}
\cortext[cor]{Corresponding author}

\begin{abstract}

This paper presents the development of a well-balanced gas-kinetic scheme (GKS) with space-time adaptive mesh refinement (STAMR) for the shallow water equations (SWE). While well-balanced GKS have been established on Cartesian and triangular meshes, the proposed STAMR framework utilizes arbitrary quadrilateral meshes with hanging nodes, introducing additional challenges for maintaining well-balanced properties. In addition to spatial adaptivity, temporal adaptivity is incorporated by assigning adaptive time steps to cells at different refinement levels, further enhancing computational efficiency.
Furthermore, the numerical flux in the GKS adaptively transitions between equilibrium fluxes for smooth flows and non-equilibrium fluxes for discontinuities, providing the proposed GKS-based STAMR method with strong robustness, high accuracy, and high resolution. Standard benchmark tests and real-world case studies validate the effectiveness of the GKS-based STAMR and demonstrate its potential for interface capturing and the simulation of complex flows.

\end{abstract}

\maketitle

\section{Introduction}

The SWE are fundamental for modeling a wide range of geophysical phenomena, with applications extending from large-scale ocean circulation to regional coastal and channel flows, including the simulation of tsunamis, tidal waves, storm surges, and dam-break events. As a result, numerous numerical schemes have been developed to effectively solve the SWE \cite{leveque-1998,zhou2001-surface,xu2002-swe,CHO2007}.
A primary challenge in these numerical methods is to maintain a precise balance between the discretization of convective fluxes and source terms, thereby preserving steady-state solutions such as the lake-at-rest condition-a property known as well-balancedness. Additional complexities include the accurate capture of discontinuities, the resolution of high-wavenumber flow features, and the robust handling of moving wet-dry interfaces.
To further enhance the fidelity of SWE models for engineering applications, it is often necessary to incorporate representations of real-world physical processes. These may include dynamic riverbed topography, fluid stratification, and external forces such as wind stress, bottom friction, and the Coriolis effect \cite{zhao2021-swe,zhao2024,delg2024}.

In practical simulations involving the SWE, achieving an optimal balance between accuracy and computational efficiency is crucial. High-order schemes are particularly advantageous in this context. Due to their enhanced discretization accuracy and superior ability to resolve a broad spectrum of waves and complex flow features, these schemes can either deliver higher accuracy on the same mesh or attain comparable accuracy using coarser meshes, thereby increasing computational efficiency. As a result, numerous high-order numerical schemes have been developed and successfully applied to the SWE \cite{xing2006-DG,dambreak-2007,highorder-efficiency,zhao2021-swe,zhao2024}.

Another widely employed strategy for enhancing computational accuracy and efficiency is adaptive mesh refinement (AMR) \cite{Berger1984}. AMR-based numerical schemes are particularly effective for problems featuring sharp interfaces or discontinuities, where conventional fixed-mesh methods-even those of high order-can suffer from excessive numerical dissipation. The core principle of AMR is to dynamically adjust the mesh resolution in response to local solution features: refining the mesh in regions with large gradients, such as discontinuities or high-wavenumber waves, and coarsening it in smoother regions where high resolution is not required. This adaptability, especially when combined with unstructured meshes, makes AMR ideally suited for simulating complex flows and irregular geometries, establishing it as a leading method for high-accuracy adaptive simulations. The process of refining or coarsening mesh cells is formally known as h-adaptivity \cite{BENK2007,XU2024}. In this study, AMR will be leveraged to develop advanced numerical schemes for shallow water flows.

In addition to h-adaptivity, several other adaptive strategies are available. One notable alternative is the moving mesh method, or r-adaptivity, which maintains the total number of nodes and their connectivity while dynamically relocating node coordinates to concentrate mesh cells in regions of significant flow variation \cite{ZHANG2024,ZHANG2025}. However, a primary limitation of r-adaptivity is its tendency toward severe mesh distortion or tangling, which can undermine simulation stability.
Adaptivity can also be achieved by locally adjusting the order of the numerical scheme, a technique known as p-adaptivity. Unlike h- and r-adaptivity, p-adaptivity operates on a fixed mesh, dynamically assigning higher-order discretizations in areas with complex flow features and lower-order discretizations in smoother regions. While conceptually powerful, p-adaptivity poses significant implementation challenges, as it requires robust and accurate error estimators to guide the selection of the appropriate polynomial degree.

In a conventional AMR framework with global time stepping, the time step for the entire computational domain is dictated by the stability constraint of the most highly refined cells. This approach is inefficient, as most coarser cells are advanced using time steps much smaller than their local stability limits, resulting in unnecessary computational overhead.
To overcome this limitation, AMR can be effectively combined with local time-stepping (LTS) techniques \cite{Osher1983, Dawson2001}. LTS enables true space-time adaptivity by allowing cells at different refinement levels to advance with locally appropriate time steps. This synergy significantly enhances computational efficiency.
In this study, we develop new AMR-based numerical schemes for the shallow water equations and employ an efficient LTS strategy to advance unsteady flows in time.

This study introduces a novel well-balanced method that integrates the GKS with AMR. The kinetic-theory foundation of GKS offers unique advantages for developing numerical schemes capable of accurately simulating complex flows. The second-order GKS was initially established for compressible flow simulations \cite{xu2} and has since evolved into a unified framework for computing both inviscid and viscous fluxes. Notably, GKS features an intrinsic multi-scale mechanism for handling both smooth flows and discontinuities, and achieves second-order accuracy within a single time-stepping procedure.
A key principle of GKS is the time-accurate evolution of the gas distribution function at cell interfaces, which has enabled the recent development of high-order compact schemes and demonstrated the method's versatility and potential for high-fidelity simulations \cite{CGKSAIA, zhao2023direct}. GKS has been successfully applied to a broad range of flow problems, including multiphase flows \cite{XU1997, LI2005}, microfluidics \cite{XU_LI_2004}, and turbulence \cite{TAN2018}, among others. In addition, GKS has been extended to solve the SWE \cite{xu2002-swe, zhao2021-swe}. In the BGK equation associated with the kinetic model for shallow water flows, an acceleration term appears to represent microscopic particle motion, corresponding to the source terms in the SWE.
Designing well-balanced GKS is particularly challenging, as the time-dependent GKS flux inherently incorporates source effects. The introduction of AMR further complicates the maintenance of the well-balanced property due to the dynamic nature of mesh adaptation. Addressing these challenges is a central contribution of this study.

The objective of this paper is to develop and implement STAMR with a well-balanced second-order GKS for the SWE, demonstrating its potential for practical applications.
The paper is organized as follows: Section 2 introduces the GKS for the SWE. Section 3 details the implementation of STAMR framework with the GKS for the SWE,
while Section 4 examines the well-balanced properties of the proposed scheme. Section 5 validates the method through a series of benchmark test cases involving shallow water flows. Finally, Section 6 summarizes the main conclusions of this study.

\section{Gas-kinetic scheme for SWE}

This section provides a brief presentation of the GKS for the SWE. The SWE and the corresponding kinetic model equation are first introduced, and their mathematical connection is established. The solution of the kinetic equation is then presented, from which the GKS for the numerical discretization of the SWE is constructed.

\subsection{SWE and kinetic model equation}
The 2-D SWE are
\begin{equation}\label{SWE-macro}
\frac{\partial \textbf{W}}{\partial t}+ \frac{\partial \textbf{F}(\textbf{W})}{\partial x}+
\frac{\partial \textbf{G}(\textbf{W})}{\partial y}=\textbf{S}(\textbf{W}),
\end{equation}
where
\begin{equation*}
{\textbf{W}} =
\left(
\begin{array}{c}
h\\
hU\\
hV\\
\end{array}
\right), \\
{\textbf{F}} =
\left(
\begin{array}{c}
hU\\
hU^2+\frac{1}{2}Gh^2\\
hUV\\
\end{array}
\right),\\
{\textbf{G}} =
\left(
\begin{array}{c}
hV\\
hUV\\
hV^2+\frac{1}{2}Gh^2\\
\end{array}
\right),
{\textbf{S}} =
\left(
\begin{array}{c}
0\\
-GhB_x\\
-GhB_y\\
\end{array}
\right).
\end{equation*}
$\textbf{W}$ represent the flow variables, and $\textbf{F}$ and $\textbf{G}$ denote the corresponding fluxes in the $x$ and $y$ directions, respectively. The bottom topography is represented by $B$, and $G$ is the gravitational acceleration constant. The source terms $\textbf{S}$ arise from the bottom profile, given by $\nabla B(x,y) = (B_x, B_y)$.

The associated kinetic model equation of the SWE, the BGK model with the inclusion of the gravitational force, is given as
\begin{equation}\label{SWE-BGK}
f_t +\textbf{u}\cdot \nabla_{\textbf{x}} f +\nabla _{\textbf{x}}\Phi\cdot\nabla_{\textbf{u}} f=\frac{g-f}{\tau},
\end{equation}
where $f$ is the distribution function, $\textbf{u} = (u, v)$ represents the particle velocity, and $g$ denotes the equilibrium state.
The relaxation time is represented by $\tau$.
The term $\nabla \Phi$ describes the particle acceleration due to external forces and is related to the source terms through $\nabla \Phi = -G \nabla B$.
The equilibrium state $g$ follows a Maxwellian distribution function \cite{xu2002-swe} given by
\begin{equation*}
\begin{split}
g=h\big(\frac{\lambda}{\pi}\big)e^{-\lambda(\mathbf{u}-\mathbf{U})^2},
\end{split}
\end{equation*}
where $\lambda$ is defined by $\lambda=1/Gh$. Due to the conservation in relaxation process from $f$ to $g$, $f$ and $g$ satisfy the compatibility condition,
\begin{equation*}
\int \frac{g-f}{\tau}\pmb{\psi} \mathrm{d}\Xi=\textbf{0},
\end{equation*}
where $\pmb{\psi}=(\psi_1,\psi_2,\psi_3)^T=(1,u,v)^T$ and $\text{d}\Xi=\text{d}u\text{d}v$.

Based on the moments of the distribution function, the flow variables and their fluxes can be obtained, which are given as
\begin{equation}\label{BGK-g-to-W}
{\textbf{W}}
=\int f \pmb{\psi} \mathrm{d}\Xi,
\end{equation}
and
\begin{equation}\label{BGK-g-to-F}
{\big(\textbf{F},\textbf{G}\big)^T}
=\int f \pmb{\psi} \mathbf{u} \mathrm{d}\Xi.
\end{equation}
The source terms $\textbf{S}$ are given by
\begin{equation}\label{BGK-g-to-S}
{\textbf{S}}
=-\int \nabla\Phi\cdot\nabla_{\textbf{u}} f \pmb{\psi} \mathrm{d}\Xi.
\end{equation}

\subsection{Time-accurate evolution solution and GKS}
The integral solution of the BGK model in Eq. (\ref{SWE-BGK}) can be obtained as
\begin{equation*}
f(\textbf{x},t,\textbf{u})=\frac{1}{\tau}\int_0^t g(\textbf{x}',t',\textbf{u}')e^{-(t-t')/\tau}\mathrm{d}t'
+e^{-t/\tau}f_0(\textbf{x}_0,\textbf{u}_0),
\end{equation*}
where $\textbf{x}$ is the location at the cell interface, set $\textbf{x}$ as $\textbf{0}$ for simplicity.
The trajectory of microscopic particle is given by $\textbf{x}=\textbf{x}'+\textbf{u}'(t-t')+\frac{1}{2}\nabla\Phi(t-t')^2$, and the velocity of the particle is $\textbf{u}=\textbf{u}'+\nabla\Phi(t-t')$.
There are two unknowns in the above integral solution. One is the initial gas distribution function $f_0$ at time $t = 0$, and the other is the equilibrium state $g$.
By approximating and modeling these unknowns in the integral solution, the time-accurate distribution function $f$ at a cell interface can be explicitly expressed as
\begin{equation}\label{SWE-2nd-order}
\begin{split}
f(\textbf{x},t,\textbf{u})&=\overline{g}\big[ C_1+ C_2  \overline{\mathbf{a}}\cdot\mathbf{u} +C_3\overline{A} \big]\\
                            &+C_2\overline{g} \big[-2 \alpha \overline{\lambda} \big( \nabla\Phi^l H(u)+\nabla\Phi^r(1-H(u)) \big) \cdot (\mathbf{u}-\overline{\mathbf{U}}) \big]\\
                            &+C_4\big[g^lH(u)+g^r(1-H(u))\big] \\
                            &+C_5g^l\big[\mathbf{a}^l\cdot \mathbf{u} -2 \alpha \lambda^l \nabla\Phi^l \cdot(\mathbf{u}-\mathbf{U}^l) \big]H(u) \\
                            &+C_5g^r\big[\mathbf{a}^r\cdot \mathbf{u} -2 \alpha \lambda^r \nabla\Phi^r \cdot(\mathbf{u}-\mathbf{U}^r) \big](1-H(u)),
\end{split}
\end{equation}
where the variables $\overline{\mathbf{a}}$, $\mathbf{a}^{l,r}$, and $\overline{A}$ are related the derivatives of equilibrium state $g$.
The parameter $\alpha$ is introduced as a correction parameter to ensure strictly well-balanced evolution; for details, see \cite{xu2002-swe,zhao2021-swe}.
The coefficients $C_i~(i=1,2,\cdots,5)$ are given as
\begin{equation*}
\begin{split}
C_1&=1-e^{-t/\tau}, ~ C_2=-\tau(1-e^{-t/\tau})+te^{-t/\tau}, ~ C_3=-\tau(1-e^{-t/\tau})+t, \\
C_4&=e^{-t/\tau}, ~ C_5=-t e^{-t/\tau}.
\end{split}
\end{equation*}
The time-dependent evolution solution of the distribution function $f$ in Eq. (\ref{SWE-2nd-order}) yields a well-balanced evolution at the mesh cell interfaces \cite{zhao2021-swe}, namely,
\begin{equation*}
\begin{split}
\textbf{W}(\mathbf{x},t)=\int f(\textbf{x},t,\textbf{u}) \pmb{\psi} \mathrm{d}\Xi=\textbf{W}(\mathbf{x},t=0),\\
\widehat{\textbf{F}}(\mathbf{x},t)=\int f(\textbf{x},t,\textbf{u}) \pmb{\psi}u \mathrm{d}\Xi=\widehat{\textbf{F}}(\mathbf{x},t=0).\\
\end{split}
\end{equation*}
Here $\widehat{\textbf{F}}$ denotes the numerical flux along the interface normal, given by $\widehat{\textbf{F}}=(\textbf{F},\textbf{G}\big)\cdot \mathbf{n}$. $\mathbf{n}$ is the unit normal vector of a cell interface.

The fluxes at cell interfaces are obtained by taking moments of the distribution function solution presented above. Leveraging the time accuracy of the solution, the total transport of mass and momentum over a time step is computed by integrating the interface fluxes in time; see \cite{xu2002-swe} for further details of the formulation.
Moreover, the source terms in the momentum equations are discretized as
\begin{equation}\label{S}
\begin{split}
\frac{1}{|\Omega_j|}\int_0^{\Delta t}\int_{\Omega_j}-Gh \nabla B \text{d}\Omega\text{d}t=-\Delta tG\frac{\overline{h}_j^n+\overline{h}_j^{n+1}}{2}\nabla \overline{B}_j.
\end{split}
\end{equation}
Consequently, the GKS for the SWE can be given as
\begin{equation}\label{SWE-GKS}
\begin{split}
\textbf{W}^{n+1}_{j}= \textbf{W}^{n}_{j} -\frac{1}{\big|\Omega_j\big|}\int_{0}^{\Delta t} \int_{\partial \Omega_j} \textbf{F}\cdot \textbf{n} \mathrm{d} l -\Delta tG\frac{\overline{h}_j^n+\overline{h}_j^{n+1}}{2}\nabla \overline{B}_j.
\end{split}
\end{equation}

\section{Space-time adaptive mesh refinement}

In this section, the implementation details of space AMR is first presented. The cell-based AMR framework p4est \cite{P4est2011} is adopted.
To enable efficient computation, we then introduce an adaptive strategy in the temporal dimension tailored to space AMR, referred to as space-time AMR.
The STAMR framework developed here serves as the foundation for constructing the well-balanced GKS.

\subsection{Implementation of space AMR}

\begin{figure}[!htb]
\centering
\includegraphics[width=0.7\textwidth]{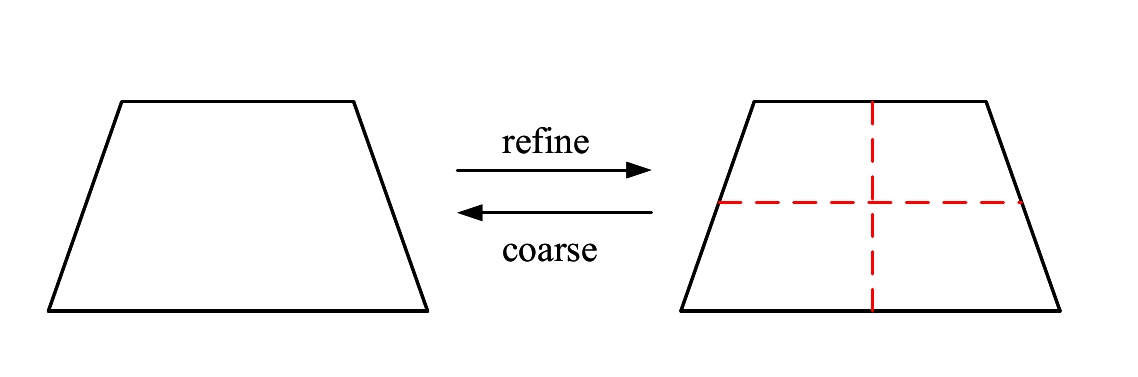}
\caption{\label{refine} Schematic of the mesh refining and coarsening process.}
\end{figure}

For two-dimensional flows, p4est uses a quadtree data structure and orders cells along a space-filling Z- (Morton) curve, enabling efficient dynamic load balancing in parallel computations. Starting from an initial coarse mesh, the AMR procedure refines any quadrilateral cell into four child cells or coarsens four child cells back into their parent cell according to a prescribed refinement criterion, as depicted in Fig. \ref{refine}.

\begin{figure}[!htb]
\centering
\includegraphics[width=1.0\textwidth]{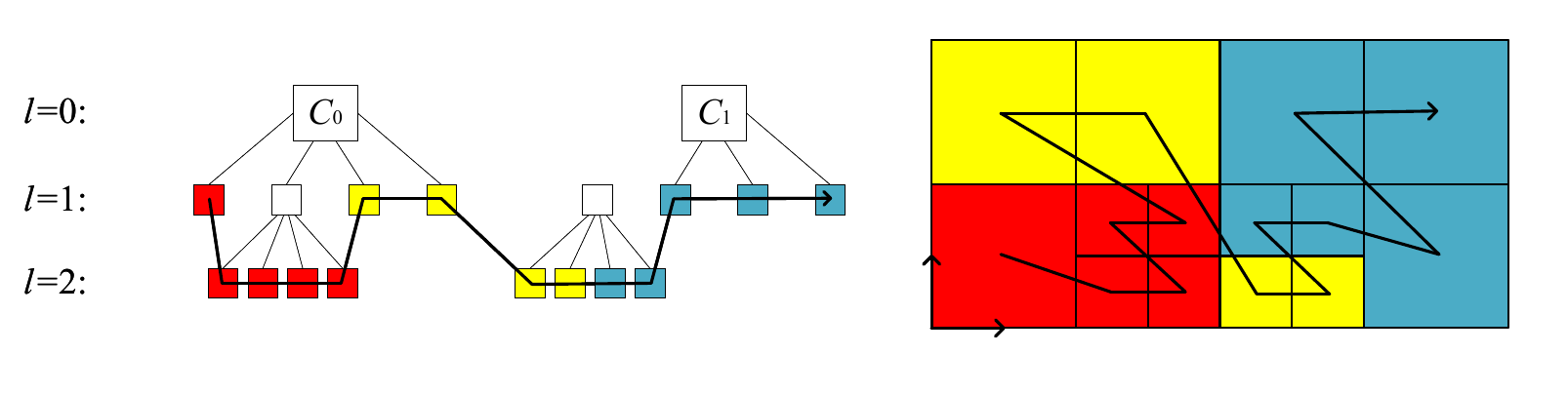}
\caption{\label{Z-order} Mesh numbering of the quadtree structure (left) and parallel partitioning (right).}
\end{figure}

Fig. \ref{Z-order} illustrates the construction of the tree data structure in the domain partitioning procedure starting from an initial mesh. The initial mesh, termed the base mesh, may consist of arbitrary quadrilateral cells, all of which are treated as parent cells.
Cells $C_0$ and $C_1$ denote two base cells that remain unrefined. The refinement level is denoted by $l$, with $l=0$ indicating the coarsest (global) level. Cells shown in white are replaced by their child cells and are not explicitly stored.
Following the ordering of the base mesh, the refined mesh is stored in one-dimensional arrays. The corresponding geometric layout is also shown in Fig. \ref{Z-order}.
A space-filling Z- (Morton) curve visits every cell in the domain, and cells along this curve are distributed across three processes (indicated by different colors) using either uniform or weighted partitioning.
Furthermore, to maintain mesh quality and simplify inter-cell communication, a strict $2:1$ balance constraint is enforced, ensuring that the refinement levels of any two face-adjacent cells differ by at most $1$.

\subsection{Efficient temporal AMR}

Space AMR effectively reduces the number of computational mesh cells, particularly in regions where flow structures exhibit significant temporal variations.
To further reduce the computational cost in the temporal dimension, time steps can be locally subdivided based on the mesh cell size.
In contrast to using a global minimum time step for time advancement, which necessitates a significantly larger number of steps and increased computational resources, temporal AMR offers a more efficient approach.
To implement temporal AMR for unsteady flows, a local preliminary time step is first calculated for each mesh cell.
The minimum value of these unweighted time steps, denoted as $ \Delta t_{min} $, is then used to determine the global time step $ \Delta t $ by applying a weighting factor based on the maximum refinement level. This relationship is determined as
\begin{equation}\label{timestep}
\Delta t = 2^{l_{max}} \Delta t_{min}.
\end{equation}
For the time interval from $ t_n $ to $ t_{n+1} $, mesh cells at refinement level $ l_i $ will perform $ 2^{l_{max} - l_i} $ time step, each using a time step size of $ 2^{l_i} \Delta t_{min} $.

\begin{figure}[!htb]
\centering
\includegraphics[width=1.0\textwidth]{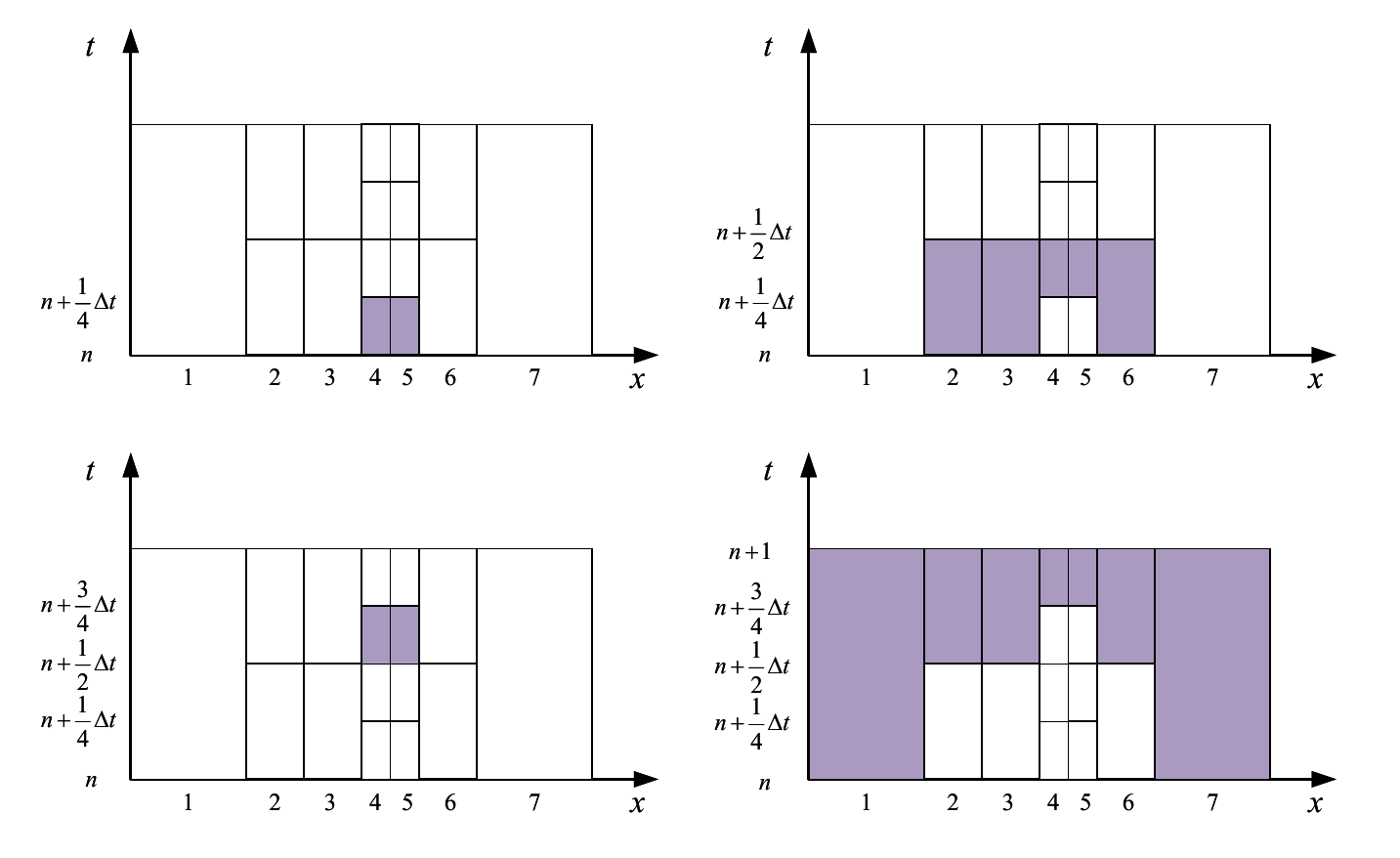}
\caption{\label{TimeAMR} The one-dimensional schematic for temporal adaptive refinement.}
\end{figure}

To illustrate the implementation of temporal adaptation, consider the one-dimensional example in Fig. \ref{TimeAMR}, where the abscissa denotes the mesh layout and the ordinate denotes time. The interval from $n$ to $n+1$ is divided into four stages. The refinement levels are: cells 1 and 7 at level 0; cells 2, 3, and 6 at level 1; and cells 4 and 5 at level 2.
\begin{enumerate}
    \item Step 1 ($n$ $\rightarrow$ $n + \frac{1}{4} \Delta t$): only cells 4 and 5 advance with $\Delta t_{min}$. Fluxes at interfaces adjacent to cells 3 and 6 are computed and stored for subsequent updates.
    \item Step 2 ($n + \frac{1}{4} \Delta t$ $\rightarrow$ $n + \frac{1}{2} \Delta t$): cells 2, 3, and 6 advance with $2\Delta t_{min}$, while cells 4 and 5 continue with $\Delta t_{min}$. When cells 2 to 6 reach $n + \frac{1}{2} \Delta t$, interface updates adjacent to cells 2 and 6, and to cells 1 and 7, are handled analogously to Step 1.
    \item Step 3 ($n + \frac{1}{2} \Delta t$ $\rightarrow$ $n + \frac{3}{4} \Delta t$): only cells 4 and 5 advance with $\Delta t_{min}$.
    \item Step 4 ($n + \frac{3}{4} \Delta t$ $\rightarrow$ $n + \Delta t$): cells 4 and 5 advance with $\Delta t_{min}$, cells 2, 3, and 6 with $2\Delta t_{min}$, and cells 1 and 7 with $4\Delta t_{min}$.
\end{enumerate}
This procedure synchronizes all cells at time $n+1$ and completes the temporal adaptation.

\section{Well-balanced GKS with STAMR}

In this section, a well-balanced GKS will be developed within an AMR framework. First, we propose a simple and effective method for discretizing source terms on quadrilateral meshes, which has not appeared in previous studies. Then, we construct a well-balanced GKS on adaptive meshes, where hanging nodes may occur along cell interfaces.

\subsection{Space discretization on quadrilateral meshes}

In this study, a piecewise linear bottom function is considered. Quadrilateral cells serve as the base mesh for adaptation within the AMR framework.
Because four arbitrary nodes are generally not coplanar, each quadrilateral is partitioned into two triangles for the sole purpose of representing the bottom topography. While all solution variables and data structures remain defined on the quadrilateral mesh.

\begin{figure}[!htb]
\centering
\includegraphics[width=0.7\textwidth]{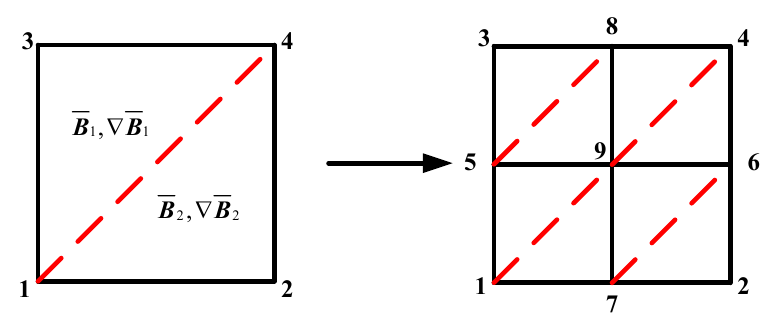}
\caption{\label{AMR_topo} Quadrilateral cell splitting on the base and adaptive meshes for source-term discretization.}
\end{figure}

First, the discretization of the source terms are presented.
As illustrated in the left plot of Fig. \ref{AMR_topo}, a quadrilateral is split into two triangles, $\Delta_{134}$ and $\Delta_{124}$. Given the bottom topography $B(x,y)$, a linear approximation over each triangle is constructed from the nodal values at its three vertices. The bottom gradient used for flux evaluations is then obtained by a linear combination of the corresponding nodal values, and cell-averaged bottom values over sub-cells are also computed. Thereby a consistent mapping between each quadrilateral cell and its triangular subcells is established.
As a result, the water height allocated to each sub-cell is given by
\begin{equation}\label{twotopo}
\begin{split}
\overline{h}_1&=\overline{h}+(1-\alpha) \overline{B}_2-\alpha \overline{B}_1,\\
\overline{h}_2&=\overline{h}-(1-\alpha) \overline{B}_2+\alpha \overline{B}_1,\\
\end{split}
\end{equation}
where $\alpha$ is the ratio of the area of triangle $\Delta_{134}$ to the quadrilateral, and $\overline{h}$ is the averaged water height over the quadrilateral.

In the SWE, the source terms in the momentum equations are linked to the bottom topography and are therefore discretized on each subcell basis.
Using the water height at times $t^n$ and $t^{n+1}$, the source terms are discretized as
\begin{equation}\label{S_2}
\begin{split}
\frac{1}{|\Omega_j|}\int_0^{\Delta t}\int_{\Omega_j}-Gh\nabla B\text{d}\Omega\text{d}t =-\Delta t G(&\alpha\frac{\overline{h}_{j,1}^n+\overline{h}_{j,1}^{n+1}}{2}\nabla \overline{B}_{j,1} \\
&+(1-\alpha)\frac{\overline{h}_{j,2}^n+\overline{h}_{j,2}^{n+1}}{2}\nabla \overline{B}_{j,2}).
\end{split}
\end{equation}

In the calculation of fluxes for the SWE, interface values of flow variables are required and obtained via spatial reconstruction.
To enhance robustness near wet-dry fronts and to handle discontinuities (e.g., shocks and pollutant interfaces) without spurious oscillations, a nonlinear reconstruction is employed.
In this study, a second-order reconstruction based on a least-square method using only immediate neighboring cells is adopted. The reconstruction is formulated for quadrilateral cells, and the cell-averaged values used to determine the reconstruction conditions are also taken over quadrilateral cells.
The reconstructed values are further limited by the limiter of Venkatakrishnan \cite{venkatakrishnan1995}, which is activated primarily in regions with discontinuities or large gradients.

In addition, the surface gradient method is used in the reconstruction, where the water surface elevation $\widetilde{h}$ ($\widetilde{h}=h + B$) is reconstructed instead of $h$ alone \cite{zhou2001-surface}. This strategy is essential for developing a well-balanced scheme, as will be demonstrated in the next subsection.
As a result, the cell interface value and gradient are given by
\begin{equation}\label{SGM}
\begin{split}
h(\mathbf{x})&=\widetilde{h}(\mathbf{x})-B(\mathbf{x}),\\
\nabla h(\mathbf{x})&=\nabla \widetilde{h}(\mathbf{x})-\nabla B(\mathbf{x}).
\end{split}
\end{equation}
Here, $\mathbf{x}$ denotes a point on the cell interface where the numerical flux is evaluated; for a second-order scheme, it is the midpoint of the interface.

\subsection{Well-balanced discretization on adaptive meshes}

In this section, the well-balanced space discretization on adaptive meshes is presented.
Under adaptive meshes, as shown in the right plot of Fig. \ref{AMR_topo}, the parent cell is subdivided into four subcells, and each subcell is further partitioned into two triangles to represent the bottom topography using piecewise linear polynomials.
Since the bottom distribution at mesh nodes is recorded, the values at the nodes of a coarsened cell are known and fixed, specifically nodes 1, 2, 3, and 4 as shown in the right plot of Fig. \ref{AMR_topo}.
For refined cells, however, the bottom elevations at the newly introduced nodes are constructed by:
\begin{equation}\label{AMR-vertice}
\begin{split}
B_5&=\frac{1}{2}(B_1+B_3),B_6=\frac{1}{2}(B_2+B_4),B_7=\frac{1}{2}(B_1+B_2),B_8=\frac{1}{2}(B_3+B_4)\\
B_9&=\frac{1}{4}(B_1+B_2+B_3+B_4)\\
\end{split}
\end{equation}
The bottom topography distribution remains consistent with that prescribed on the initial base mesh. It is recomputed only during mesh refinement and coarsening, as well as in the discretization of the source terms, on the new mesh and on the auxiliary subcells.

\begin{figure}[!htb]
\centering
\includegraphics[width=0.65\textwidth]{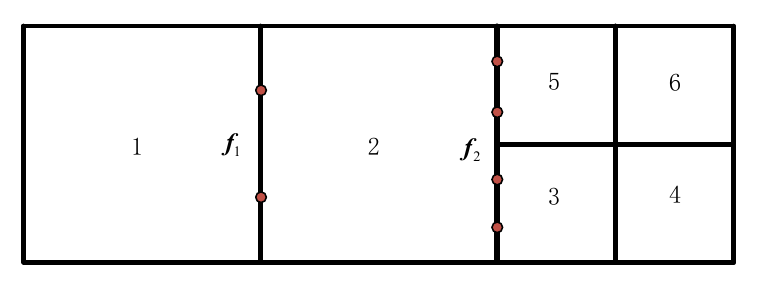}
\caption{\label{hanging-edge} Adaptive meshes, hanging nodes and gauss points for space discretization.}
\end{figure}

Hanging edges may arise during isotropic refinement or coarsening of a quadrilateral mesh, as shown in Fig. \ref{hanging-edge}.
Typically, the flux on a hanging edge, denoted by $ f_2 $, is evaluated using the adjacent fine cells 3 and 5. While the flux $ f_1 $ is evaluated based on the coarse adjacent cell 1. Consequently, within the STAMR framework, the flux $ f_1 $ generally differs from that on the hanging edge $ f_2 $ for the coarse cell 2. This discrepancy prevents the pressure contribution from the flux and the bottom force from the source terms from maintaining exact balance, thereby destroying the scheme's well-balanced property.
For a nominal second-order scheme, a single Gauss point suffices to represent the average flux across a cell interface. In the present scheme, however, two Gauss points (marked by red circles in Fig. \ref{hanging-edge}) are employed along the interface to ensure the well-balanced property by accurately computing the flux contribution.  After obtaining the exact flux, we can analyze the well-balanced property. Eq. (\ref{SWE-GKS})  can be written as

\begin{equation}\label{well-balance-flux}
\begin{split}
\textbf{W}^{n+1}_{j}= \textbf{W}^{n}_{j}+\Delta t(L_j^F+L_j^S)^n.
\end{split}
\end{equation}
where $F$ and $S$ represent flux and source terms respectively. For quadrilateral, $L_j^F$ can be written as
\begin{equation}\label{well-balance-S1}
\begin{split}
L_j^F= -\frac{1}{\big|\Omega_j\big|} \sum_{l=1}^{4}\big( \big|\Gamma_{l} \big| \sum _{k=1}^2 \frac{1}{4}Gh_k^2\cdot \textbf{n}_l \big)
\end{split}
\end{equation}

Since the quadrilateral is divided into two triangles, the flux along the red edge for each triangle in Fig. \ref{AMR_topo} will cancel each other out. Thus, Eq. (\ref{well-balance-S1}) can be further expressed as follows:
\begin{equation}\label{well-balance-S2}
\begin{split}
L_j^F=& -\frac{1}{\big|\Omega_j\big|} \big(\sum_{l=1}^{3}\big( \big|\Gamma_{l,1} \big| \sum _{k=1}^2 \frac{1}{4}Gh_k^2\cdot \textbf{n}_{l,1} \big) +\sum_{l=1}^{3}\big( \big|\Gamma_{l,2} \big| \sum _{k=1}^2 \frac{1}{4}Gh_k^2\cdot \textbf{n}_{l,2} \big)\big)\\
=&-\frac{\big|\Omega_{j,1}\big|}{\big|\Omega_j\big|} \big(\frac{1}{\big|\Omega_{j,1}\big|}\sum_{l=1}^{3}\big( \big|\Gamma_{l,1} \big| \sum _{k=1}^2 \frac{1}{4}Gh_k^2\cdot \textbf{n}_{l,1} \big) \big) \\
&-\frac{\big|\Omega_{j,2}\big|}{\big|\Omega_j\big|} \big(\frac{1}{\big|\Omega_{j,2}\big|}\sum_{l=1}^{3}\big( \big|\Gamma_{l,2} \big| \sum _{k=1}^2 \frac{1}{4}Gh_k^2\cdot \textbf{n}_{l,2} \big) \big)\\
=&-\alpha G\nabla\overline{h}_{j,1}\overline{h}_{j,1}-(1-\alpha)G\nabla\overline{h}_{j,2}\overline{h}_{j,2}.
\end{split}
\end{equation}

This condition is precisely balanced with the source terms if $( h + B )$ remains constant at each quadrature point of the cell interface after reconstruction.
This requirement is met by using the surface-gradient method for spatial reconstruction.

\section{Numerical validations}
In this section, we will test the well-balanced GKS with STAMR framework (GKS-STAMR) for SWE. The time step used in the computations is determined by the CFL condition with CFL$=0.4$. The gravitational acceleration is taken as $G=9.812$ unless otherwise specified. The collision time $ \tau $ for inviscid flow at a cell interface is defined by:
\begin{equation}\label{tau}
\begin{split}
\tau=C_1\Delta t+C_2|\frac{h_l^2-h_r^2}{h_l^2+h_r^2}|\Delta t
\end{split}
\end{equation}
where $C_1=0.05$, $C_2=5$ are used in all cases. In these simulations, the threshold values of refinement and coarsening criteria are set to the same value, and mesh refinement is performed at each time step for unsteady flows.

The scheme of SWE-coupled pollutant transport was developed and detailed in \cite{zhao2021-swe}. The present study implements the same scheme on STAMR framework for simulating pollutant transport, and related details are therefore omitted.

\subsection{Well-balanced property}
The well-balanced property of GKS-STAMR for SWE is validated by a numerical test case. The initial condition is a two-dimensional steady state with non-flat bottom topography. The bottom topography is
\begin{align*}
& B(x,y)=0.5e^{-50[(x-1)^2+(y-1)^2]}.
\end{align*}
The steady state is
\begin{align*}
\begin{split}
 h=1-B(x,y)
\end{split}
\end{align*}
and all velocities are set to $ 0 $. The computational domain is $ [0, 2] \times [0, 2] $, and the original base mesh has a cell size of $ h_{\text{mesh}} = 0.05 $. The maximum refinement level is $ l_{\text{max}} = 2 $. As time progresses, AMR is gradually applied from the lower left corner to the upper right corner, allowing us to investigate the well-balanced properties of the scheme during the dynamic changes in the mesh. A wall boundary condition is imposed on all boundaries.

\begin{figure}[!htb]
\centering
\includegraphics[width=0.45\textwidth]{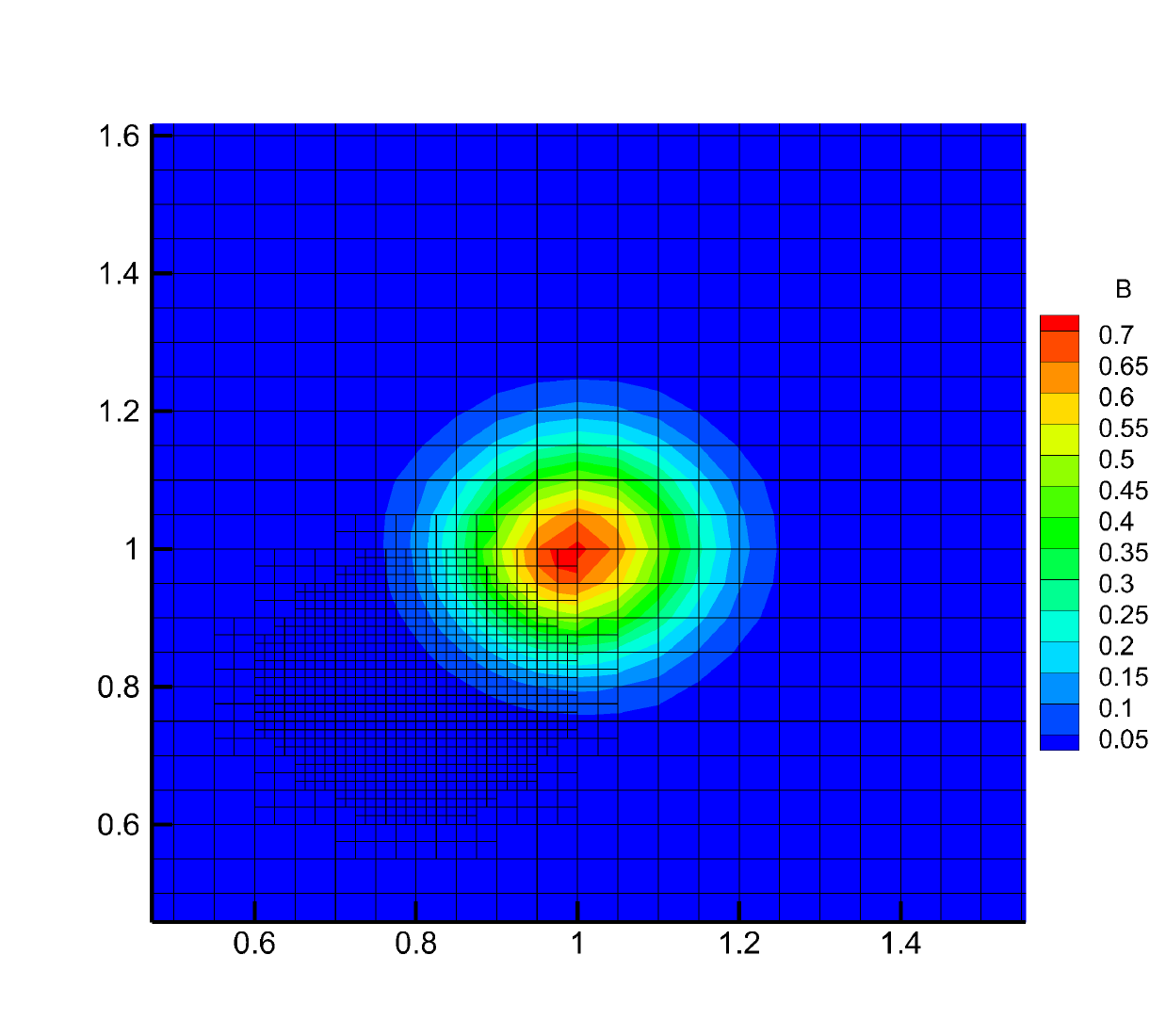}
\includegraphics[width=0.45\textwidth]{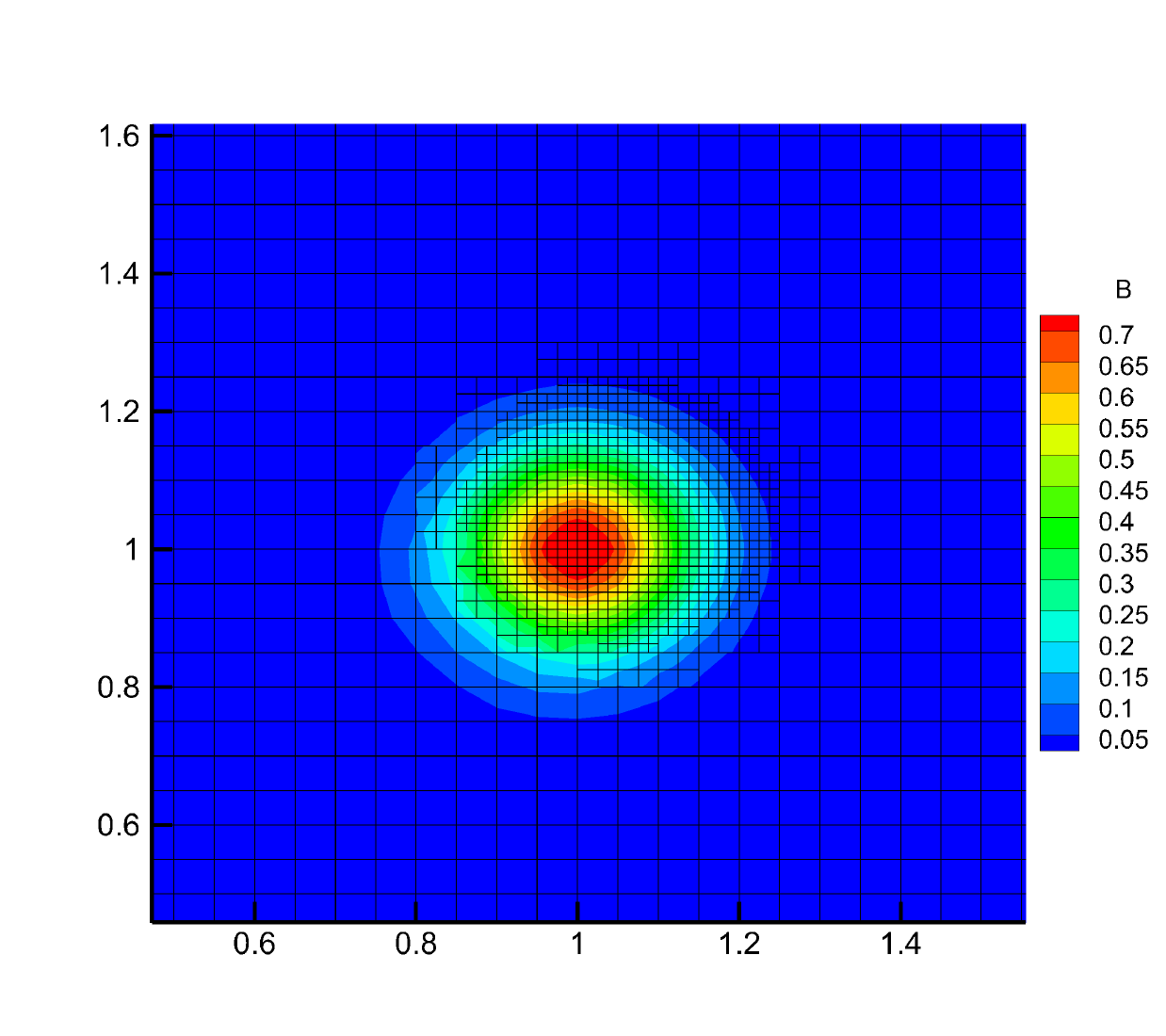}
\caption{\label{Well-balance-mesh} Well-balanced test: the bottom profile and mesh distribution, left: $t=0$, right: $t=5$.}
\end{figure}

As shown in Fig. \ref{Well-balance-mesh}, after a prolonged calculation, the areas with significant changes in the bottom topography have undergone both coarsening and refinement of the meshes, which involves the interpolation of both the bottom and the flow variables. The errors of the flow variables at different times are presented in Table \ref{Well-balance-error}. The error remains at a consistent level across different computational output times, indicating that the well-balanced properties are preserved.

\begin{table}[!h]
\renewcommand\arraystretch{1.5}
	\begin{center}
		\def\temptablewidth{0.6\textwidth}
		{\rule{\temptablewidth}{0.60pt}}
        \footnotesize
		\begin{tabular*}{\temptablewidth}{@{\extracolsep{\fill}}cccc}
         time & Error $L_1(h)$  &Error $L_1(hU)$ &Error $L_1(hV)$    \\
			\hline
            0.1    & 3.058e-16    & 6.568e-16     & 6.420e-16       \\
            1.0    & 1.710e-15    & 6.077e-15     & 6.132e-15       \\
            10.0  & 8.767e-16    & 2.385e-14     & 2.392e-14       \\
		\end{tabular*}
		{\rule{\temptablewidth}{0.60pt}}
	\end{center}
	\vspace{-6mm} \caption{\label{Well-balance-error} Well-balanced property test: the errors of flow variables of the solutions obtained by GKS-STAMR. }
\end{table}

\begin{figure}[!htb]
\centering
\includegraphics[width=0.45\textwidth]{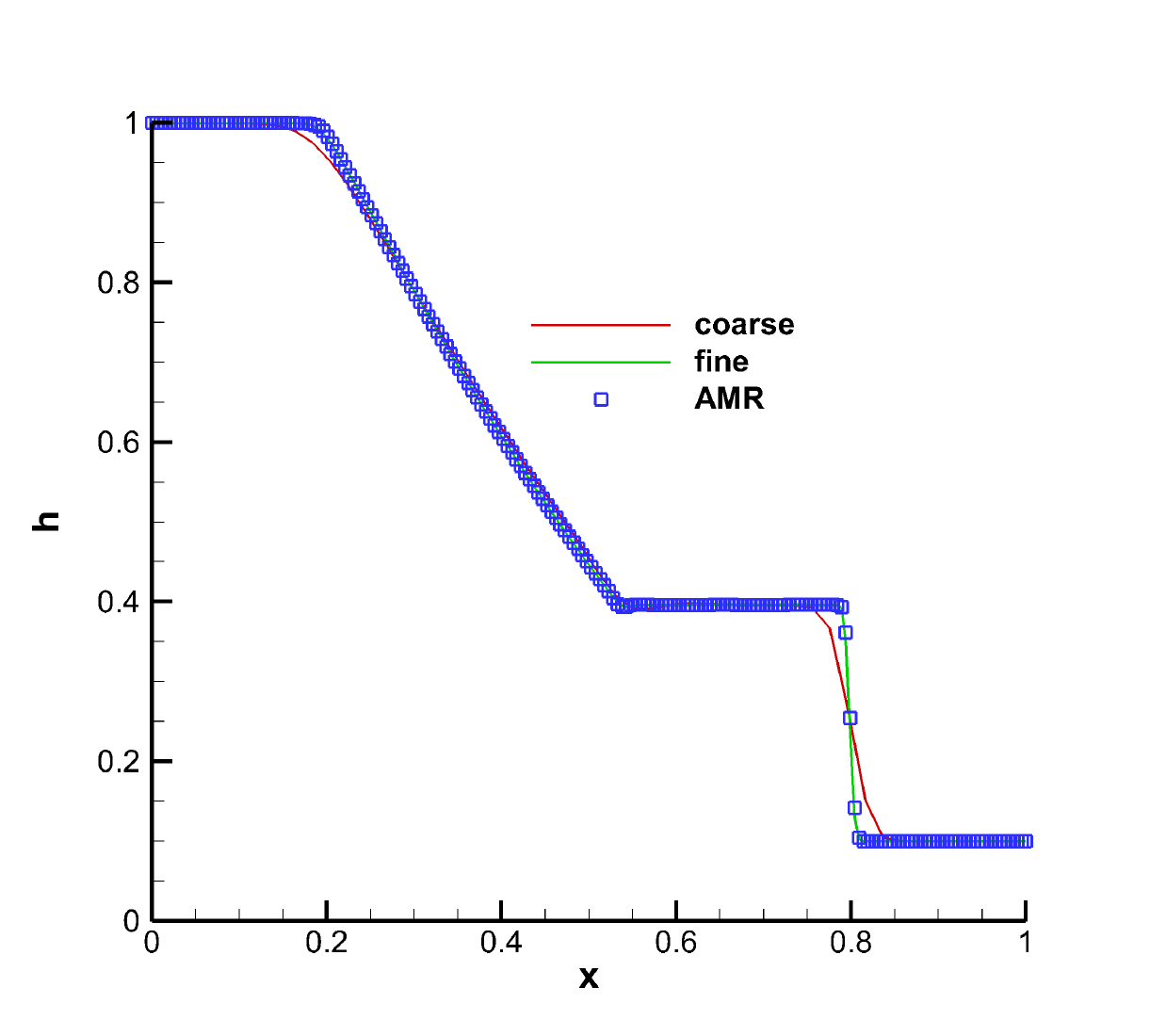}
\includegraphics[width=0.45\textwidth]{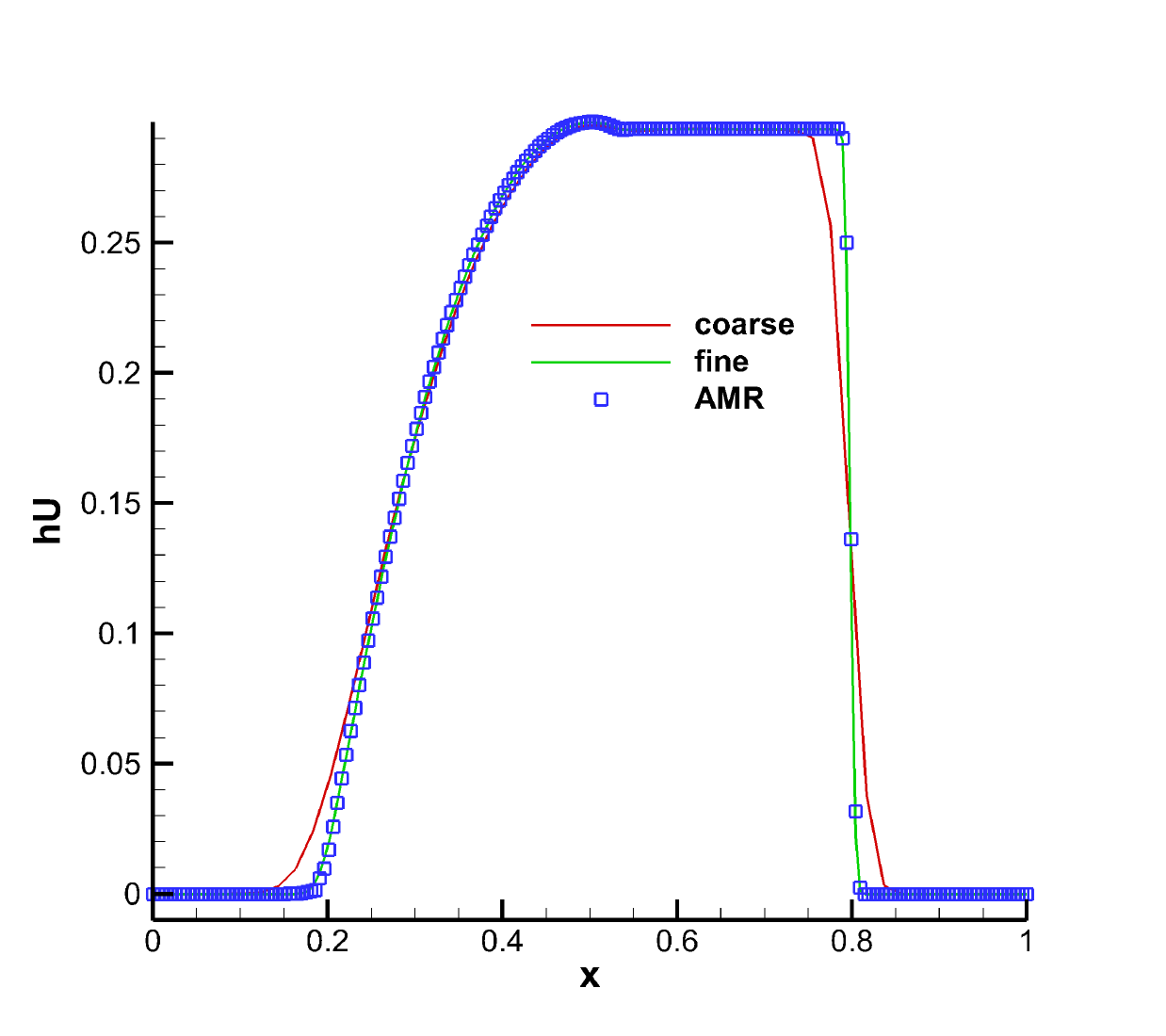}
\includegraphics[width=0.45\textwidth]{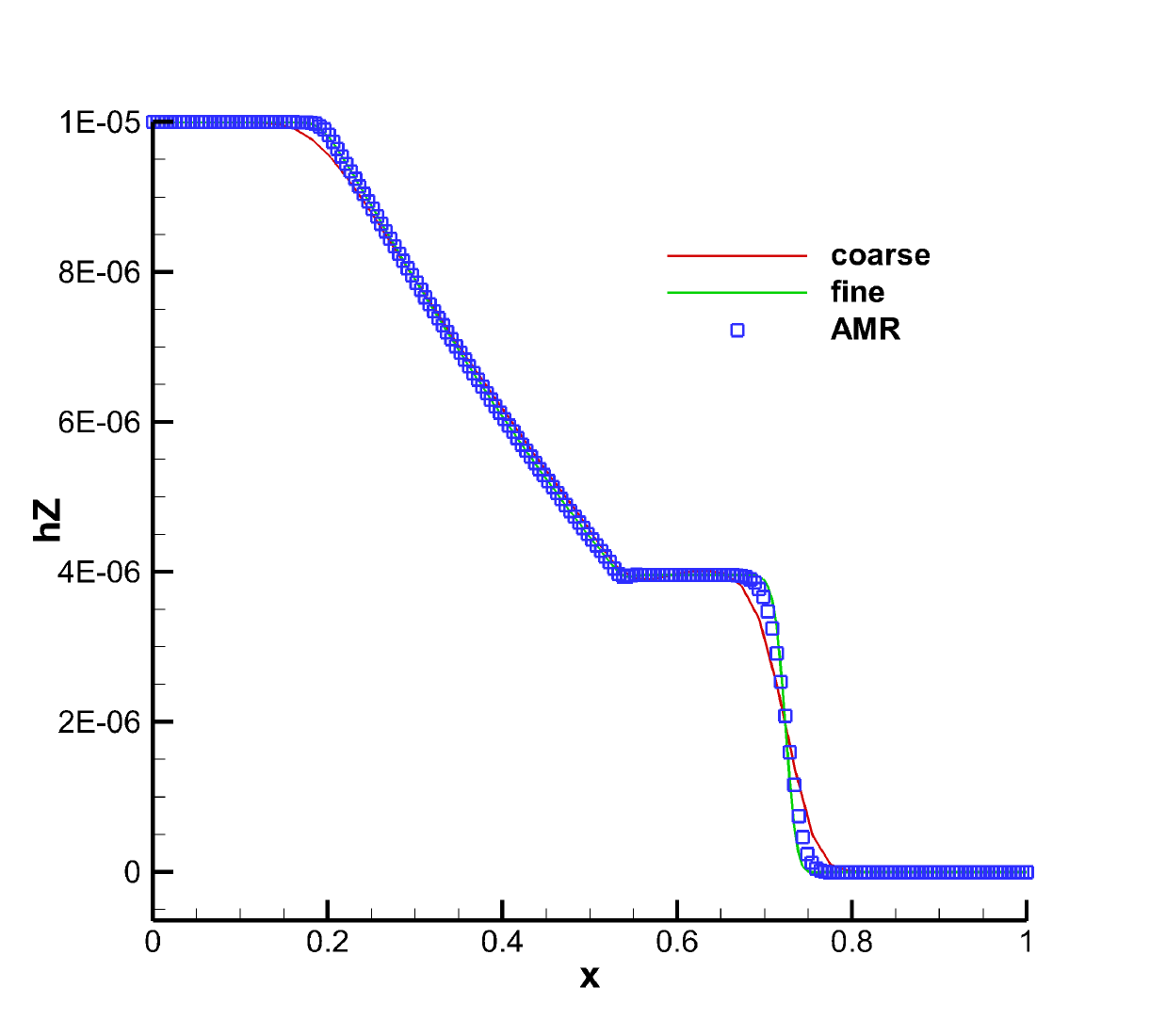}
\includegraphics[width=0.45\textwidth]{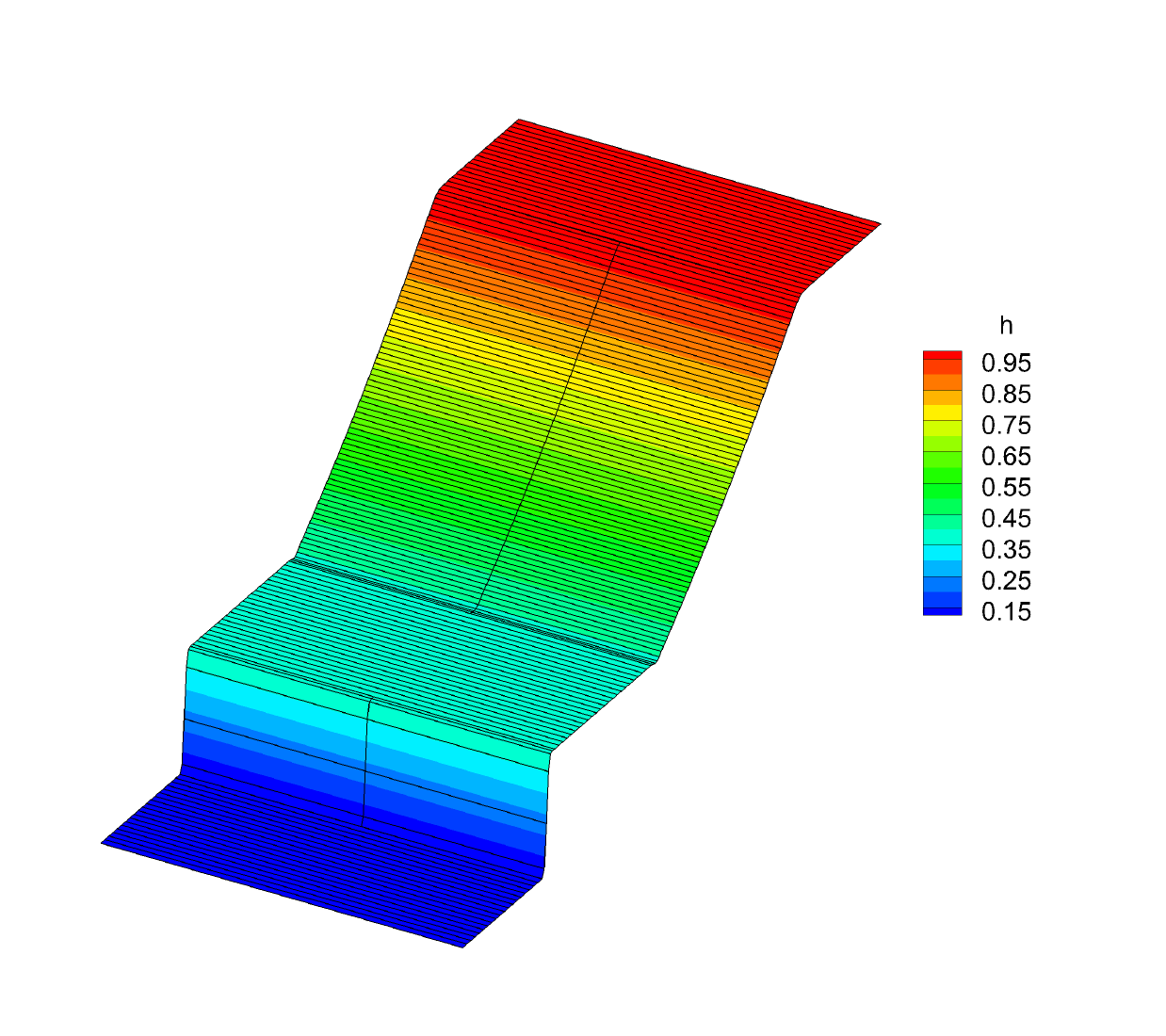}
\caption{\label{shocktube} 1-D dam-break problem with a flat bottom topography. the water surface, discharge, scalar distribution and mesh distribution at t=0.2 along the horizontal centerline of the computational domain.}
\end{figure}

\subsection{1-D dam-break problem}
The 1-D dam-break problem is simulated over flat bottom topography. The initial condition is
\begin{equation*}
(h,U,V) = \begin{cases}
(1,0,0),~~~~ &0\leq x<0.5,\\
(0.1,0,0), &0.5\leq x\leq1,
\end{cases}
\end{equation*}
The intensity of the shock and expansion waves in the shock tube problem depends on the ratio of the initial water depths of the upstream and downstream regions. The gravitational acceleration is set to $ G = 1.0 $. The computational domain is $ [0, 1]\times [0,0.04] $. The reflecting boundary condition is imposed in the $y$ direction. The free boundary condition is used on the left and right boundaries.

Three types of meshes are tested: a coarse uniform mesh with a cell size of $ h_{\text{mesh}} = 0.02 $, a fine uniform mesh with a cell size of $ h_{\text{mesh}} = 0.01 $, and AMR that begins with the coarse mesh at $ l_{\text{max}} = 1 $, the refinement criterion is chosen as $\beta=|\nabla h|/|\nabla h|_{max}$, and the threshold is 0.03. A three-layer mesh is established in the $ y $ direction, with the mesh sizes consistent with those in the $ x $ direction. Additionally, to verify the scheme for the coupled approach with scalar transport, a scalar transport simulation is included in the weak case. The initial condition for the scalar is set to:
\begin{equation*}
Z= \begin{cases}
1\times10^{-5},~~~~ &0\leq x<0.5,\\
0, &0.5\leq x\leq1,
\end{cases}
\end{equation*}
The results at $t=0.2$ are shown in Fig. \ref{shocktube}. The resolution is obviously insufficient under the coarse mesh. AMR can refine the mesh near the rarefaction waves and shock waves, and the resolution is comparable to that under the fine mesh.

\begin{figure}[!htb]
\centering
\includegraphics[width=0.45\textwidth]{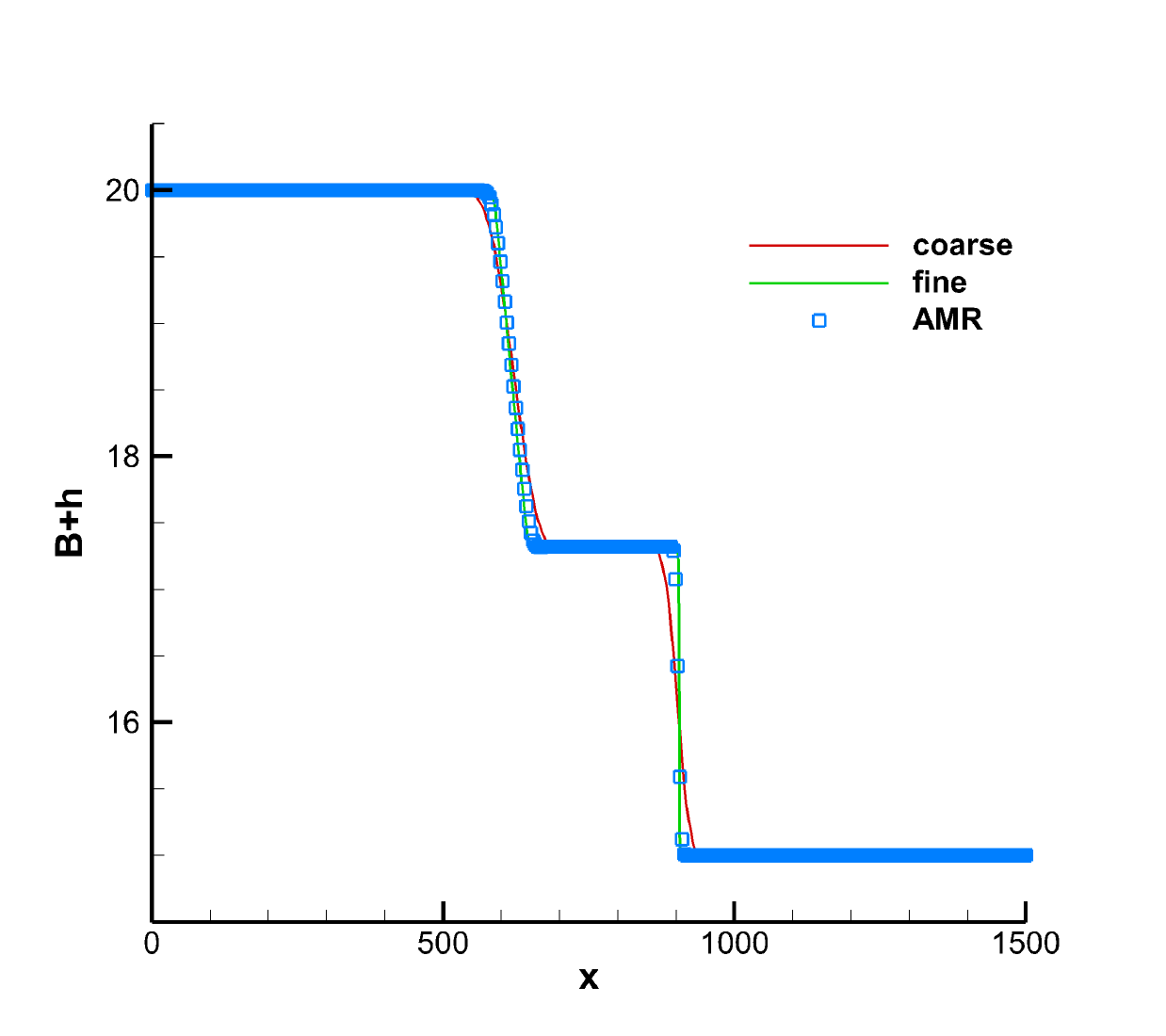}
\includegraphics[width=0.45\textwidth]{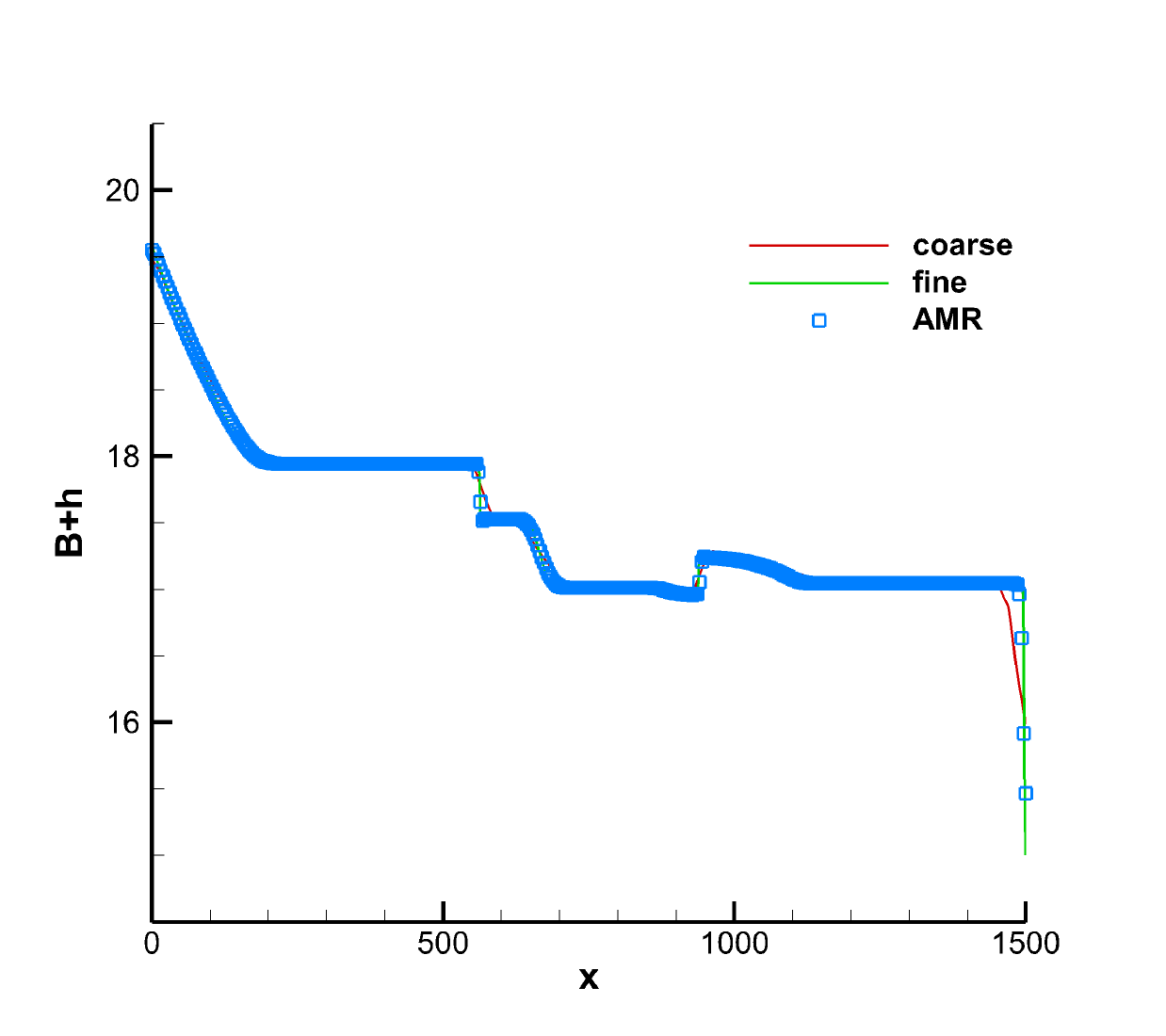}
\caption{\label{1d-dambreak} 1-D dam-break problem with a discontinuous bottom topography: the distributions
of water surface ($B+h$) along the horizontal centerline at time $t=15$ (left) and $t=60$ (right).}
\end{figure}

\subsection{1-D dam-break problem with non-flat bottom topography}
The third one is the dam-break problem with a non-flat bottom topography, where a rectangular bump is considered. The bottom topography is taken as
\begin{equation*}
B(x,y)= \begin{cases}
8, &|x-750|\leq1500/8,\\
0, &\mathrm{otherwise}.
\end{cases}
\end{equation*}
The computational domain is $[0, 1500]\times[0,7.5]$. The same boundary conditions as in the dam-break problem are adopted.

Three types of meshes are tested: a coarse uniform mesh with $N_{\text{mesh}}=400$, a fine uniform mesh with $N_{\text{mesh}}=1600$, and AMR that begins with the coarse mesh, where the maximum refinement level is $l_{\text{max}}=2$, the refinement criterion is chosen as $\beta=|\nabla h|/|\nabla h|_{max}$, and the threshold is $0.03$.
The initially discontinuous water depth is
\begin{equation*}
h= \begin{cases}
20-B(x,y), &x\leq750,\\
15-B(x,y), &x>750.
\end{cases}
\end{equation*}
The initial velocities are $U = V = 0$. The results of different initial meshes are shown in Fig.\ref{1d-dambreak}. The results obtained on the AMR mesh are in good agreement with those on the fine mesh, demonstrating a significant improvement over the coarse mesh.

\begin{figure}[!htb]
\centering
\includegraphics[width=0.45\textwidth,trim=30 10 10 10,clip]{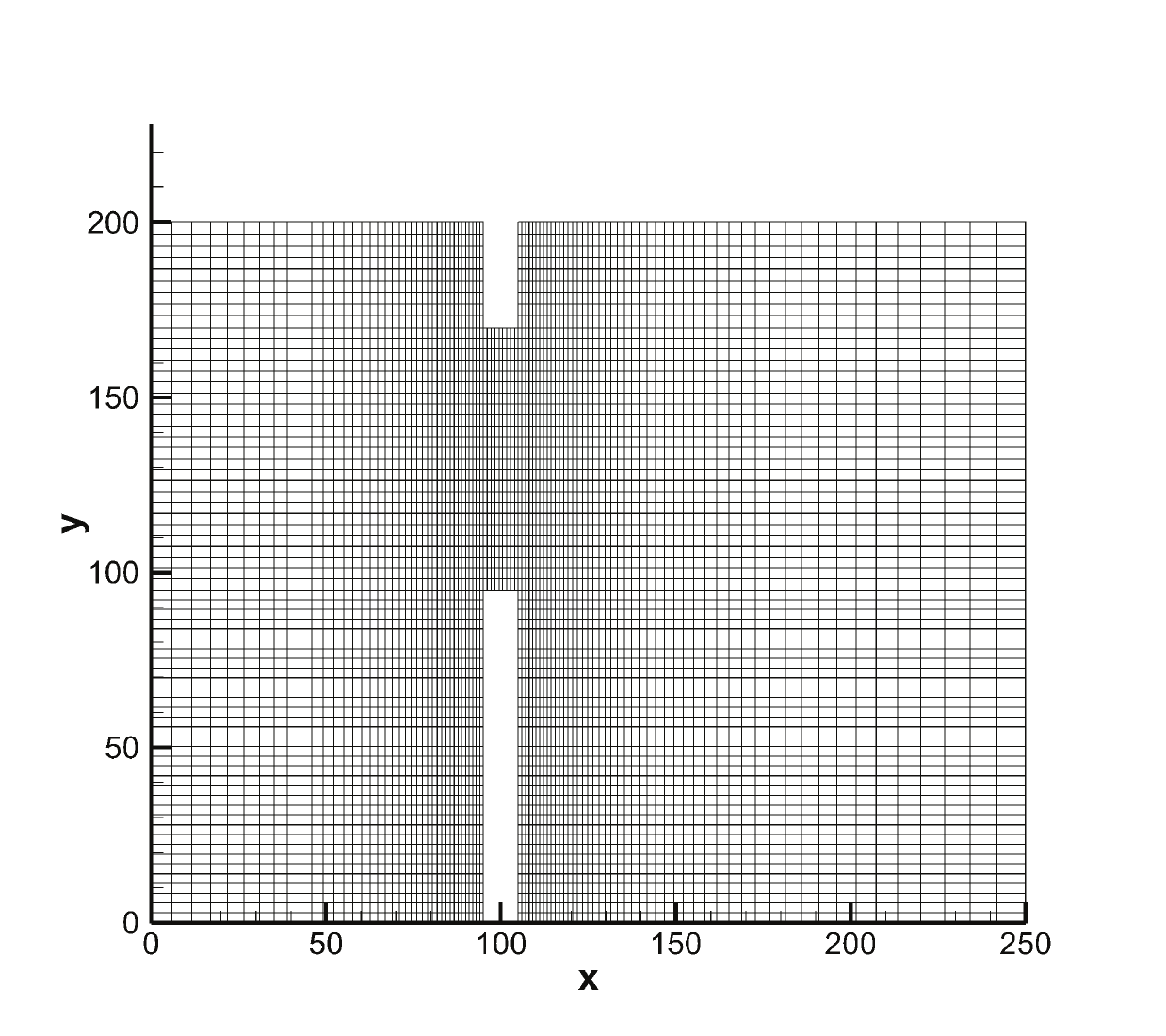}
\includegraphics[width=0.45\textwidth,trim=30 10 10 10,clip]{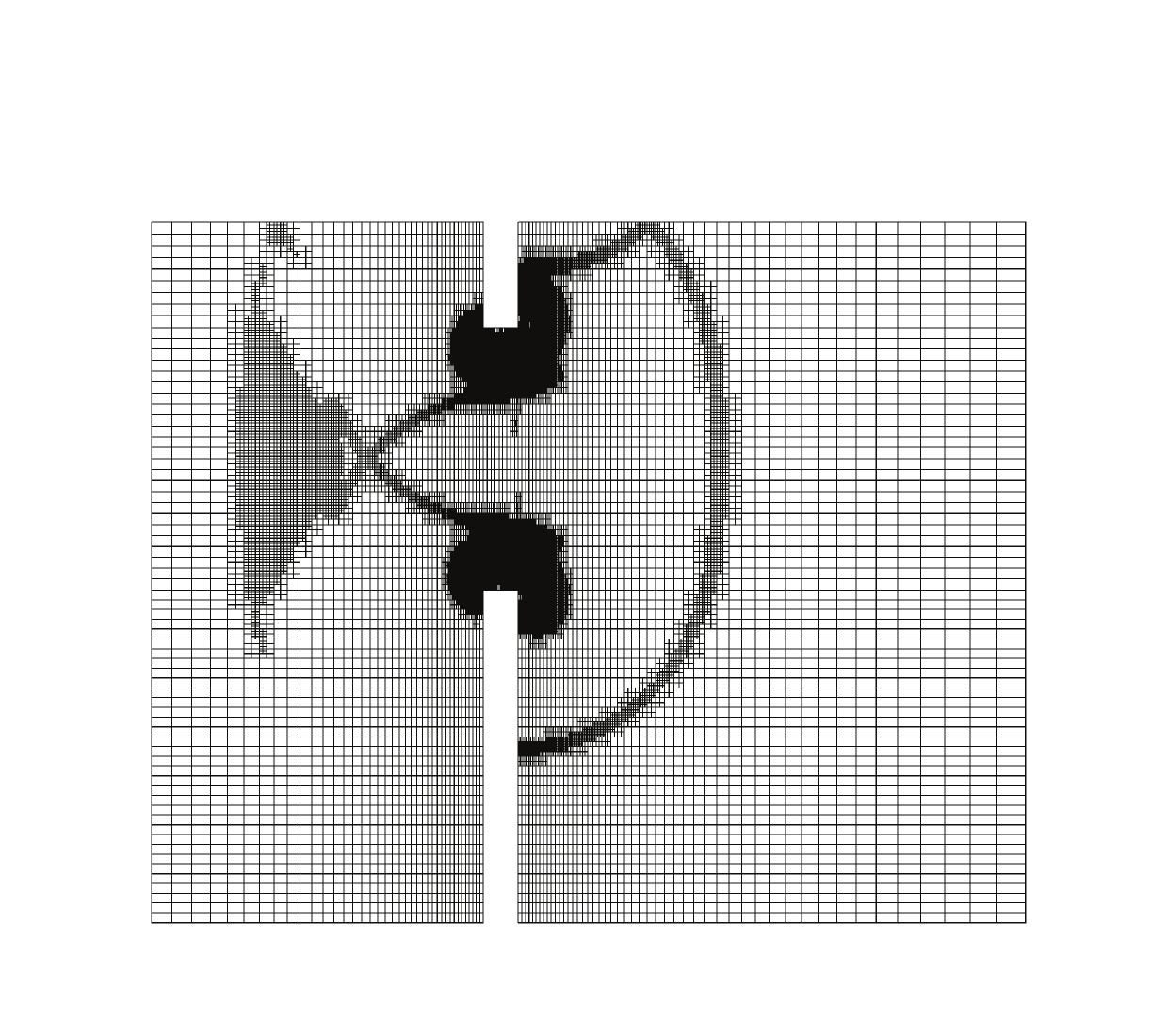}
\includegraphics[width=0.45\textwidth,trim=30 10 10 10,clip]{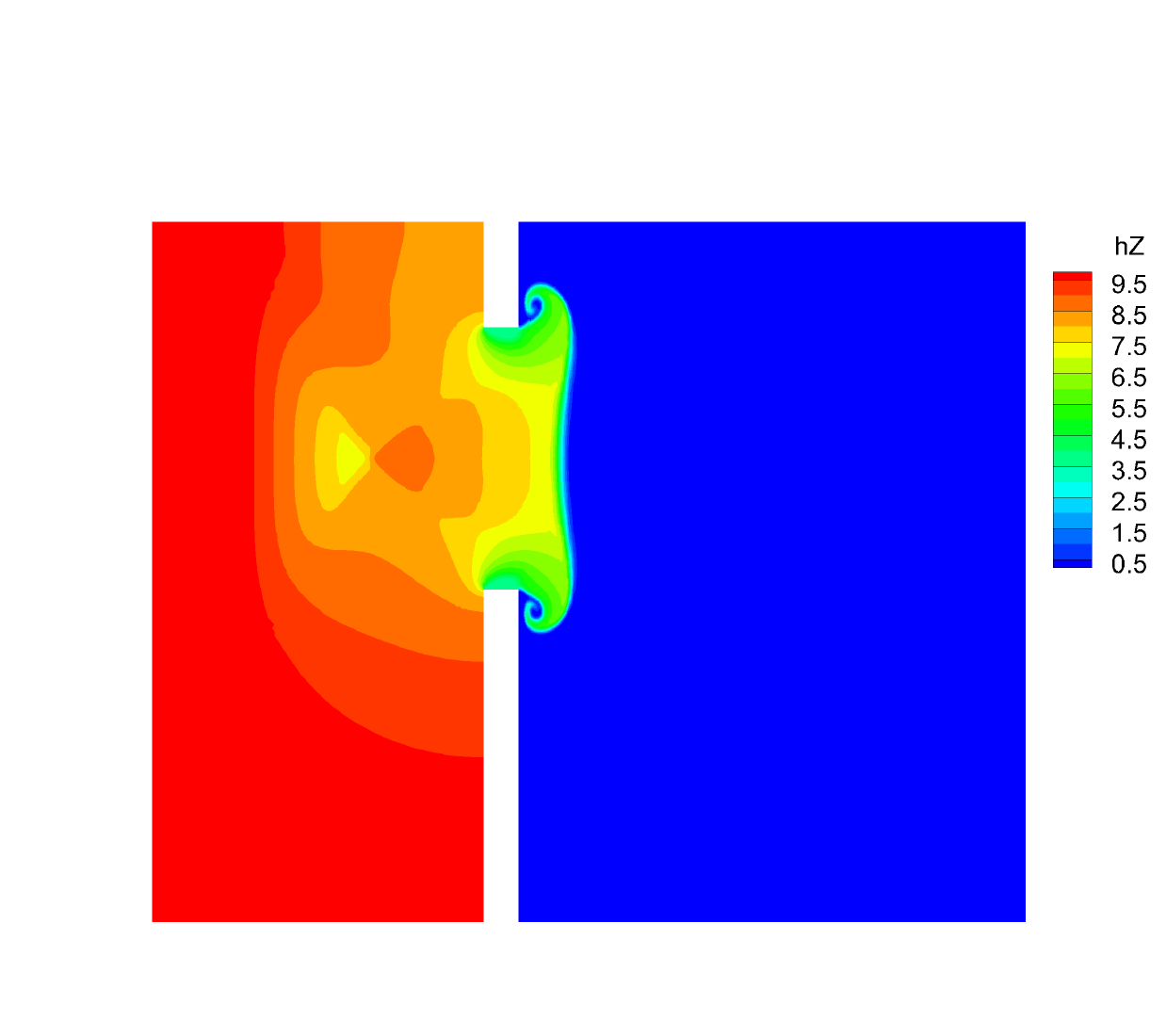}
\includegraphics[width=0.45\textwidth,trim=30 10 10 10,clip]{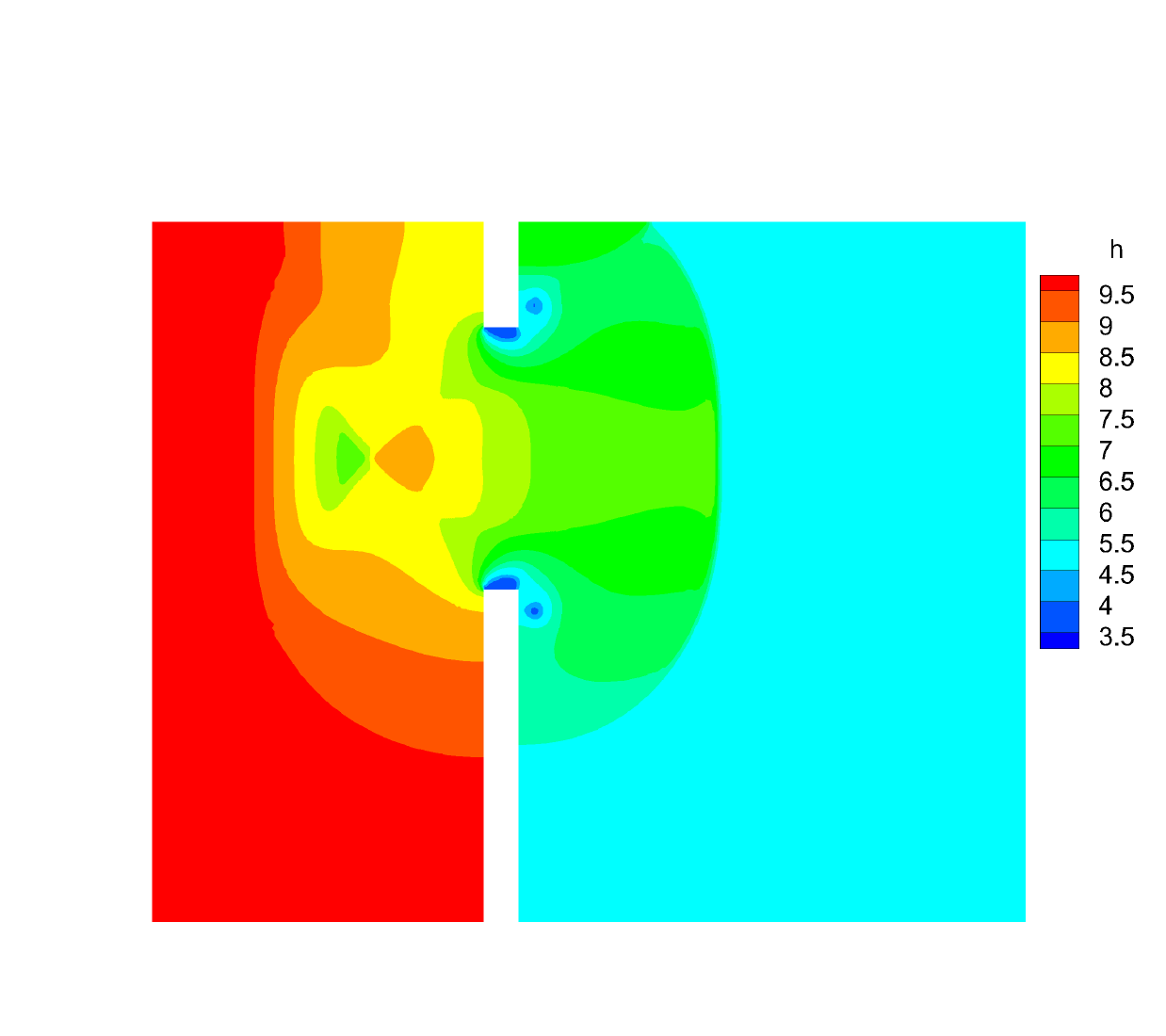}
\caption{\label{2d-dambreak} 2-D dam-break problem: initial and final mesh distribution, contours of  water height and scalar at $t = 7.2$.}
\end{figure}

\begin{figure}[!htb]
\centering
\includegraphics[width=0.45\textwidth,trim=00 10 10 10,clip]{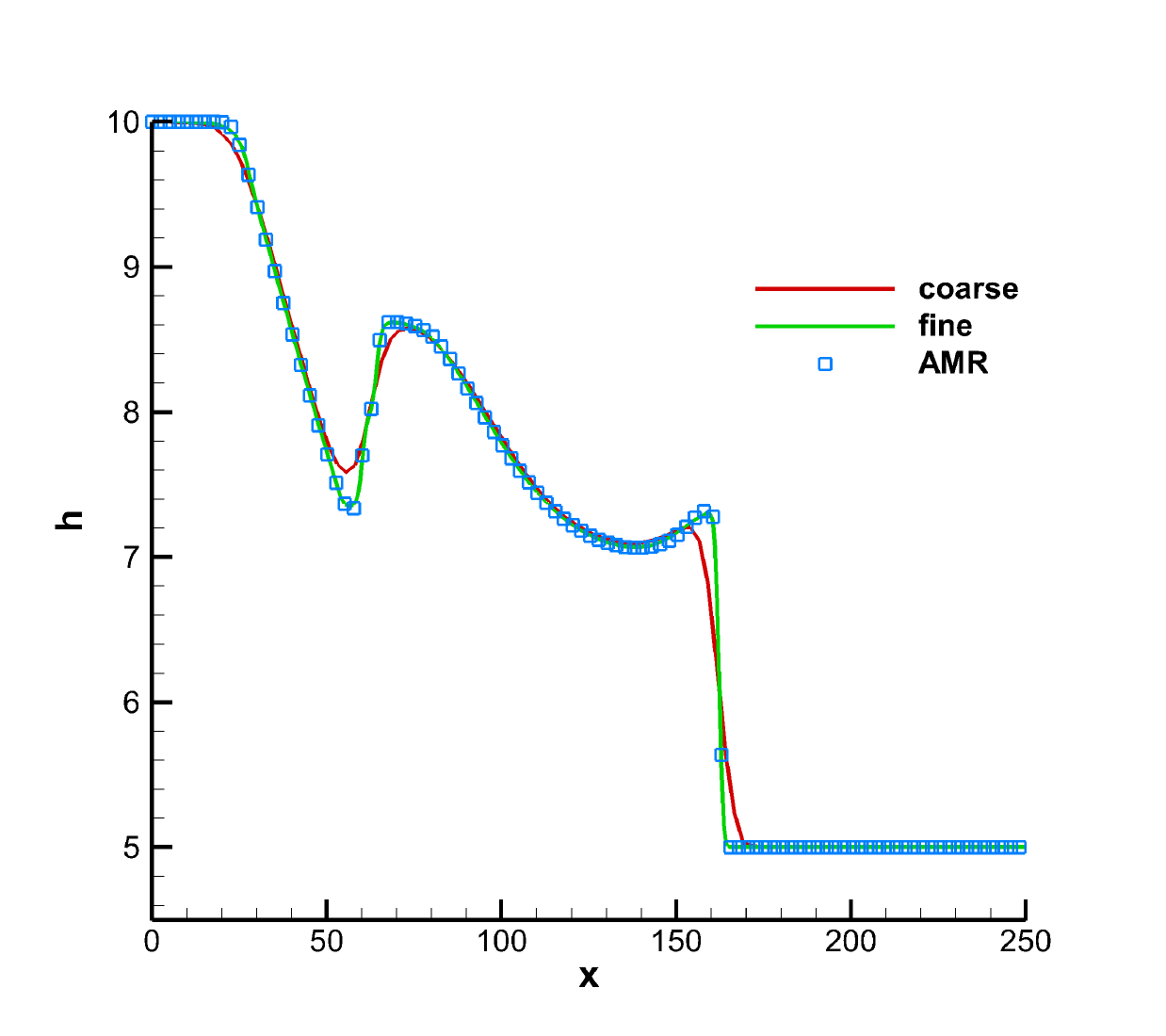}
\includegraphics[width=0.45\textwidth,trim=00 10 10 10,clip]{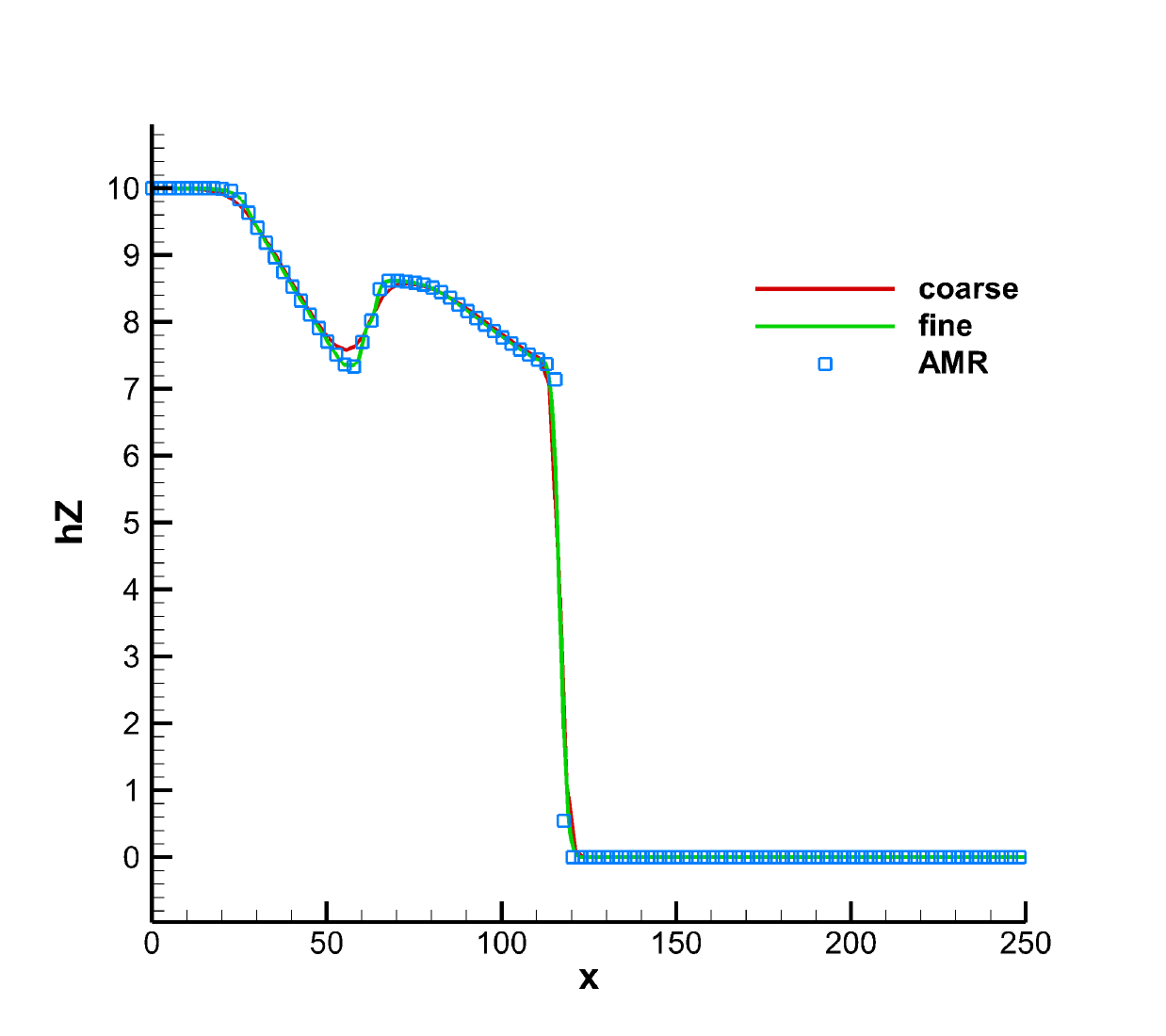}
\caption{\label{2d-dambreak-line} 2-D dam-break problem: distribution of water height and scalar along the centerline at $t = 7.2$.}
\end{figure}

\subsection{2-D dam-break problem with flat bottom topography}
The 2-D dam-break problem is tested on a complex domain. The length of the dam breach is $75$, starting at $y = 95$. The dam itself has a width of $10$, with its left side also located at $y = 95$. At $t = 0$, the water surface exhibits a discontinuity, with $h_{left}= 10$ and $h_{right} = 5$ on both sides of the breach.
The initial condition of the scalar $Z$ is set as
\begin{equation*}
Z(x,y)= \begin{cases}
1, &x\leq 95,\\
0, &\mathrm{otherwise}.
\end{cases}
\end{equation*}
The boundary condition on the far right is a free boundary, while the other boundary conditions are wall boundaries.

The numerical solutions in Fig. \ref{2d-dambreak} are sufficiently smooth in the smooth regions, and the discontinuities are captured without spurious oscillations. The contours of the water surface and the ``pollutant'' height $h_Z$ are also plotted. The mesh in the complex flow structure is adaptively refined, while the other smooth areas retain the original coarse mesh. Additionally, the flow variables along the centerline are shown in Fig. \ref{2d-dambreak-line}. Three types of meshes are compared, and areas with significant changes require higher resolution, which highlights the importance of AMR.

The number of meshes at $t = 7.2$ and the CPU time costs are compared in Table \ref{2d-dambreak-efficiency}. AMR refers to the use of only space adaptation, while STAMR indicates that both space and time adaptations are employed. All begin with the same initial coarse mesh, and the minimum cell size during the calculations remains consistent. It can be observed that, with the help of AMR, the computational number of meshes is significantly reduced, and the CPU time is considerably smaller. Additionally, after adopting time adaptation, the efficiency can be further improved.

\begin{table}[!h]
\renewcommand\arraystretch{1.5}
	\begin{center}
		\def\temptablewidth{0.6\textwidth}
		{\rule{\temptablewidth}{0.60pt}}
        \footnotesize
		\begin{tabular*}{\temptablewidth}{@{\extracolsep{\fill}}cccc}
         Scheme & Number of mesh cells  &CPU time(s) &Ratio    \\
			\hline
            Uniform    & 97792   & 240   & 1.0     \\
            AMR                   & 18682   & 100   & 2.4      \\
            STAMR               & 18682   & 60     &  4.0      \\
		\end{tabular*}
		{\rule{\temptablewidth}{0.60pt}}
	\end{center}
	\vspace{-6mm} \caption{\label{2d-dambreak-efficiency} 2-D dam-break problem: the number of mesh cells at $t=7.2$ and CPU times.}
\end{table}

\begin{figure}[!htb]
\centering
\includegraphics[width=0.475\textwidth,trim=30 30 60 60,clip]{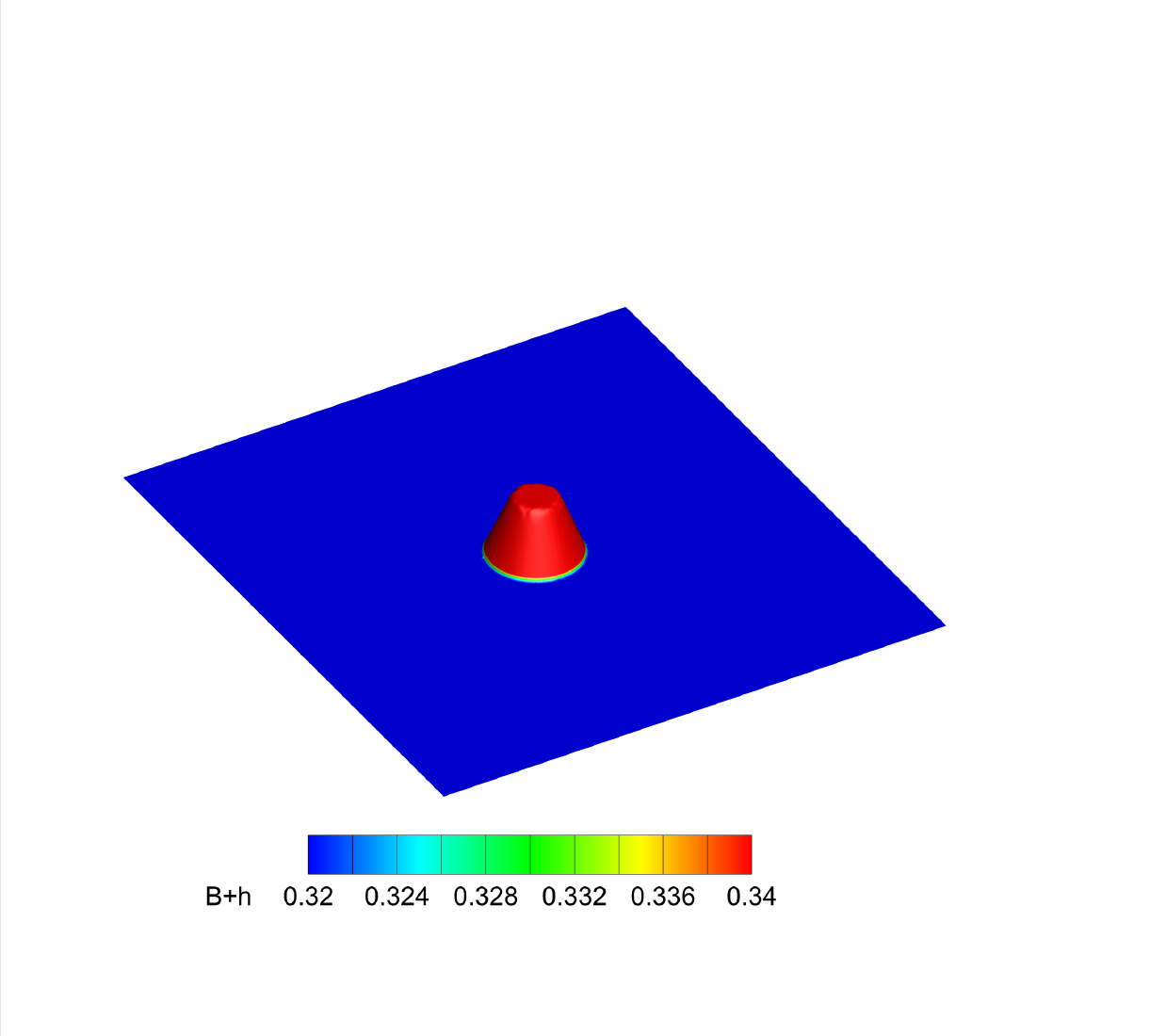}
\includegraphics[width=0.475\textwidth,trim=30 30 60 60,clip]{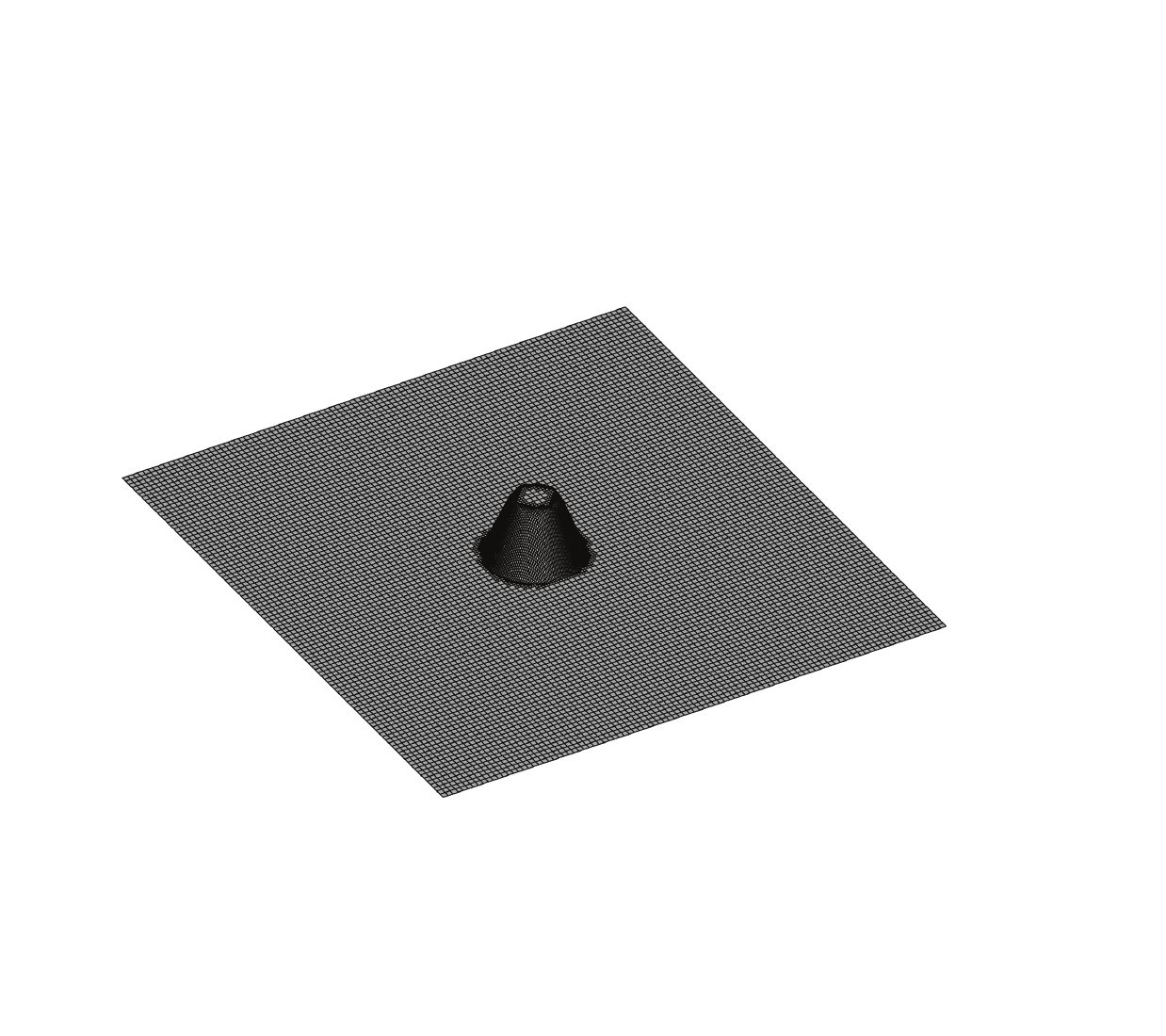}
\includegraphics[width=0.475\textwidth,trim=30 30 60 60,clip]{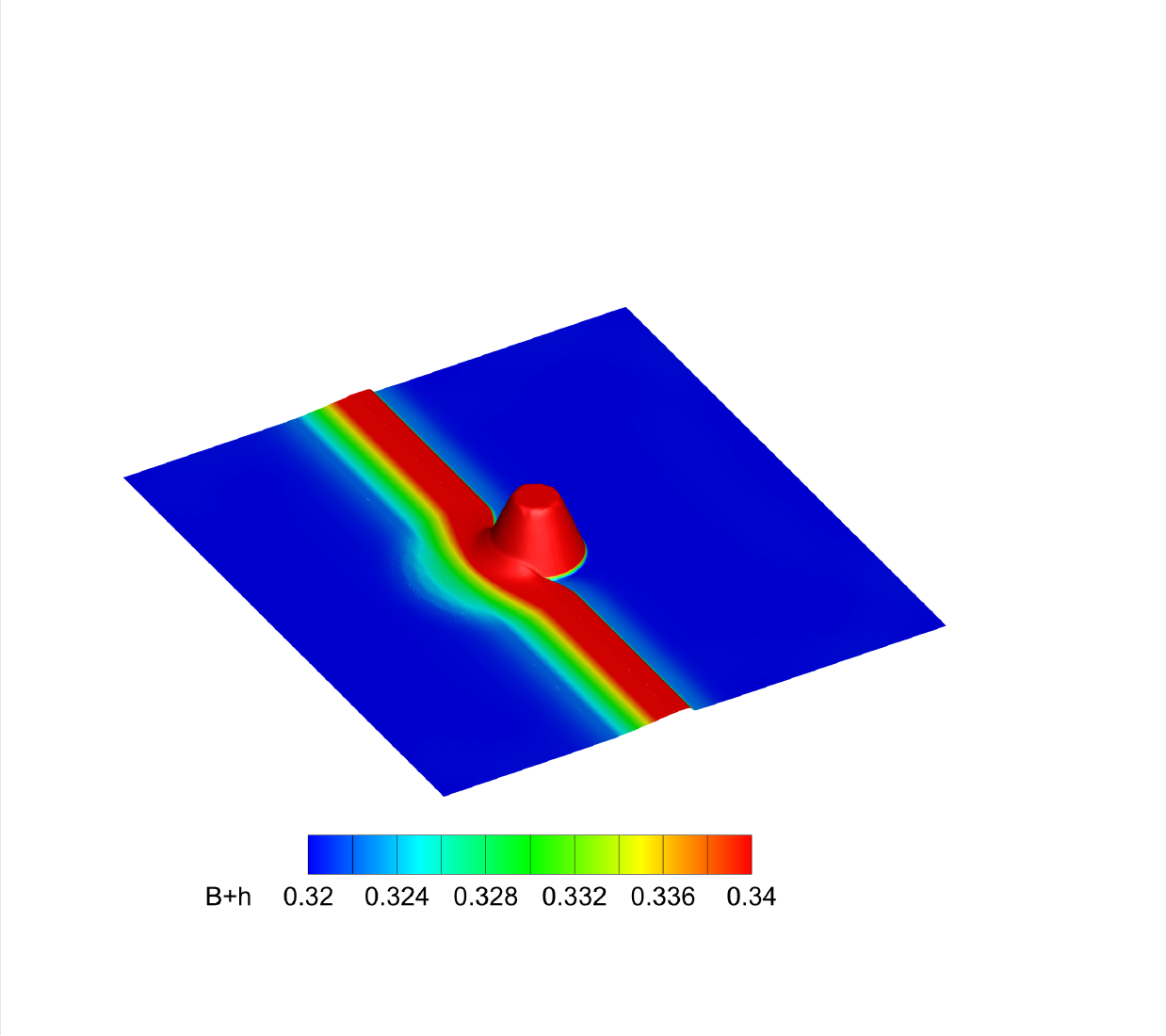}
\includegraphics[width=0.475\textwidth,trim=30 30 60 60,clip]{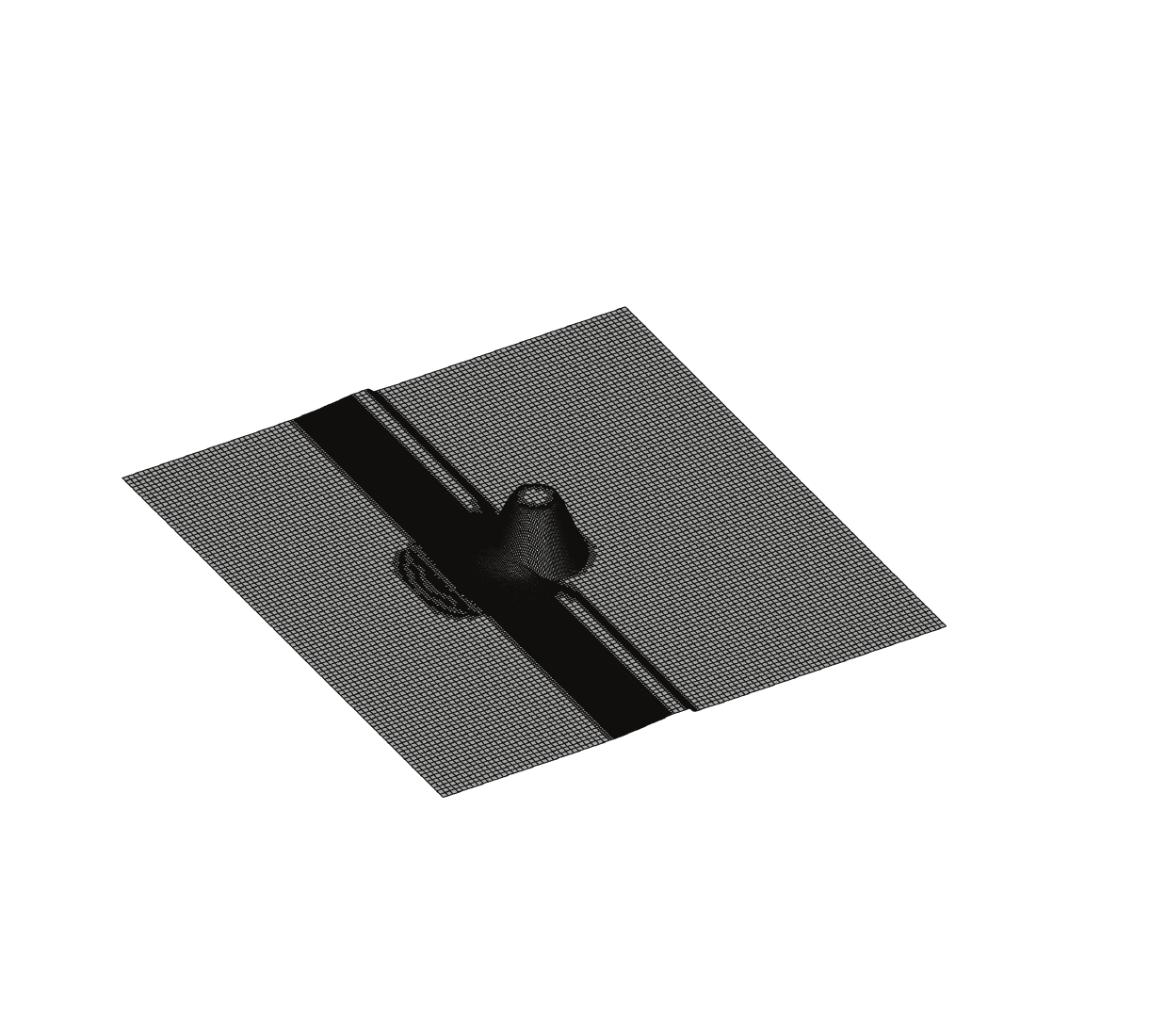}
%\end{figure}
%
%\begin{figure}[!htb]
%\centering
%\includegraphics[width=0.45\textwidth,trim=30 30 60 60,clip]{tsunami-10-height}
%\includegraphics[width=0.45\textwidth,trim=30 30 60 60,clip]{tsunami-10-mesh}
\includegraphics[width=0.475\textwidth,trim=30 30 60 60,clip]{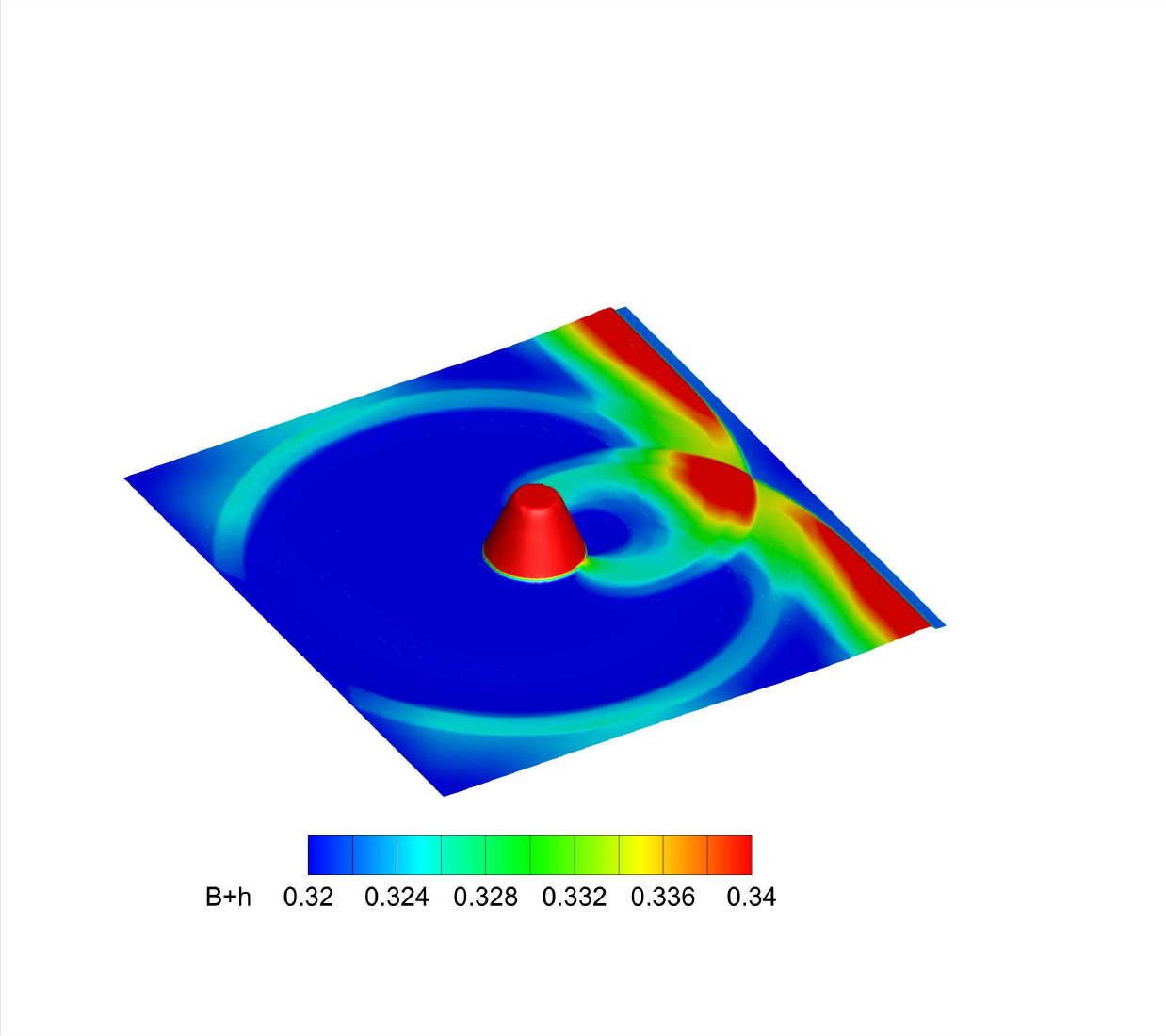}
\includegraphics[width=0.475\textwidth,trim=30 30 60 60,clip]{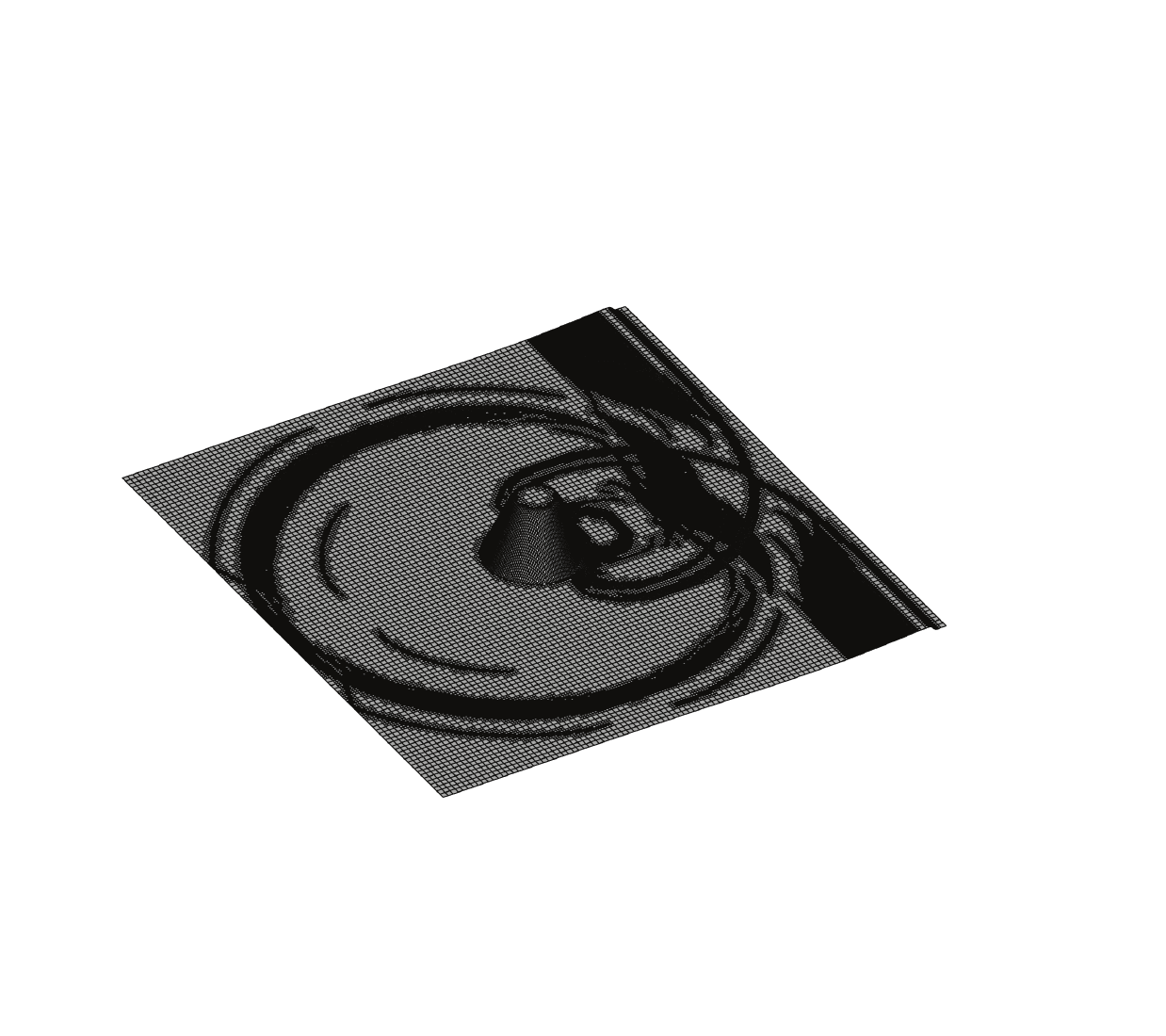}
\end{figure}

\begin{figure}[!htb]
\centering
\includegraphics[width=0.475\textwidth,trim=30 30 60 60,clip]{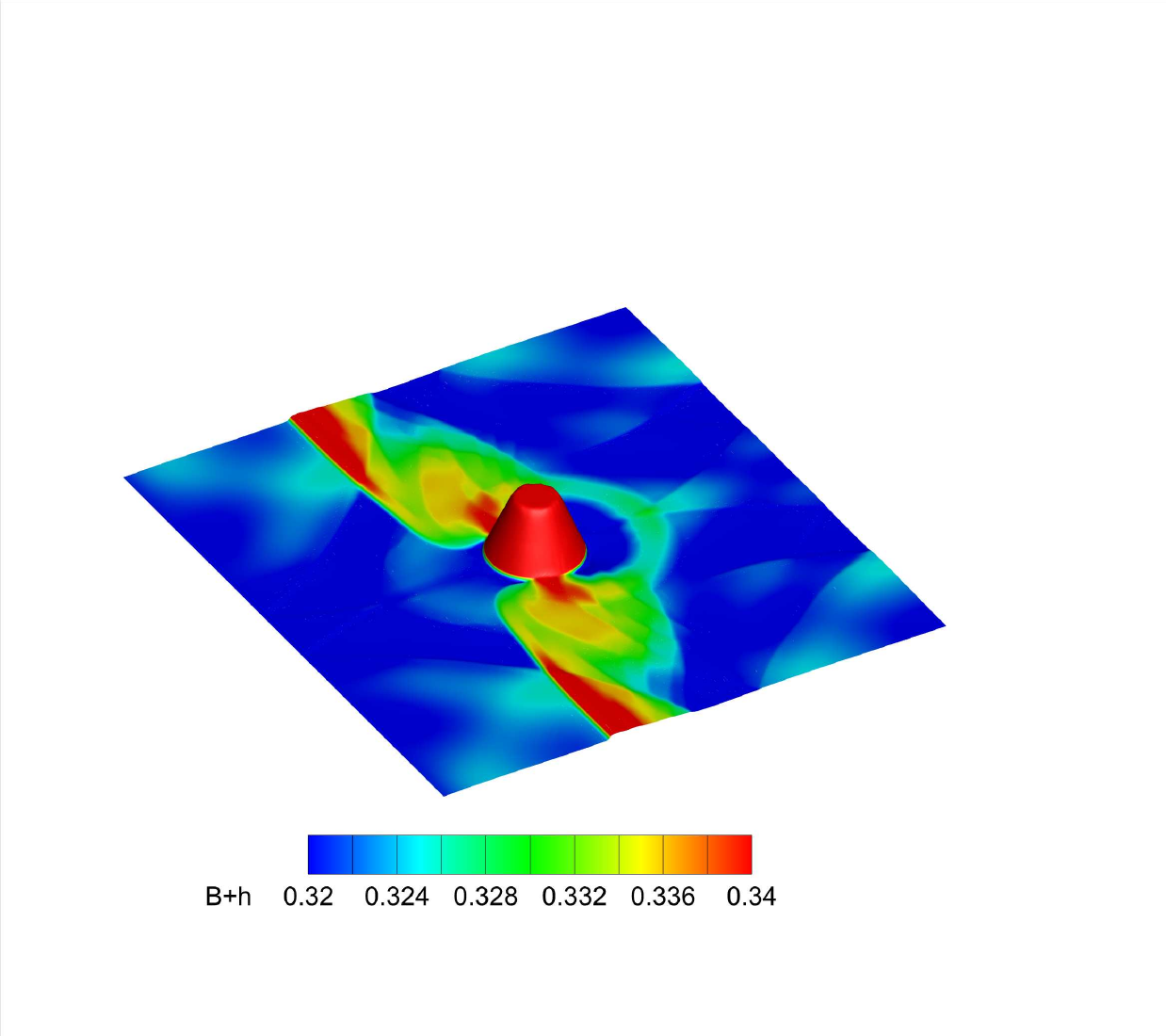}
\includegraphics[width=0.475\textwidth,trim=30 30 60 60,clip]{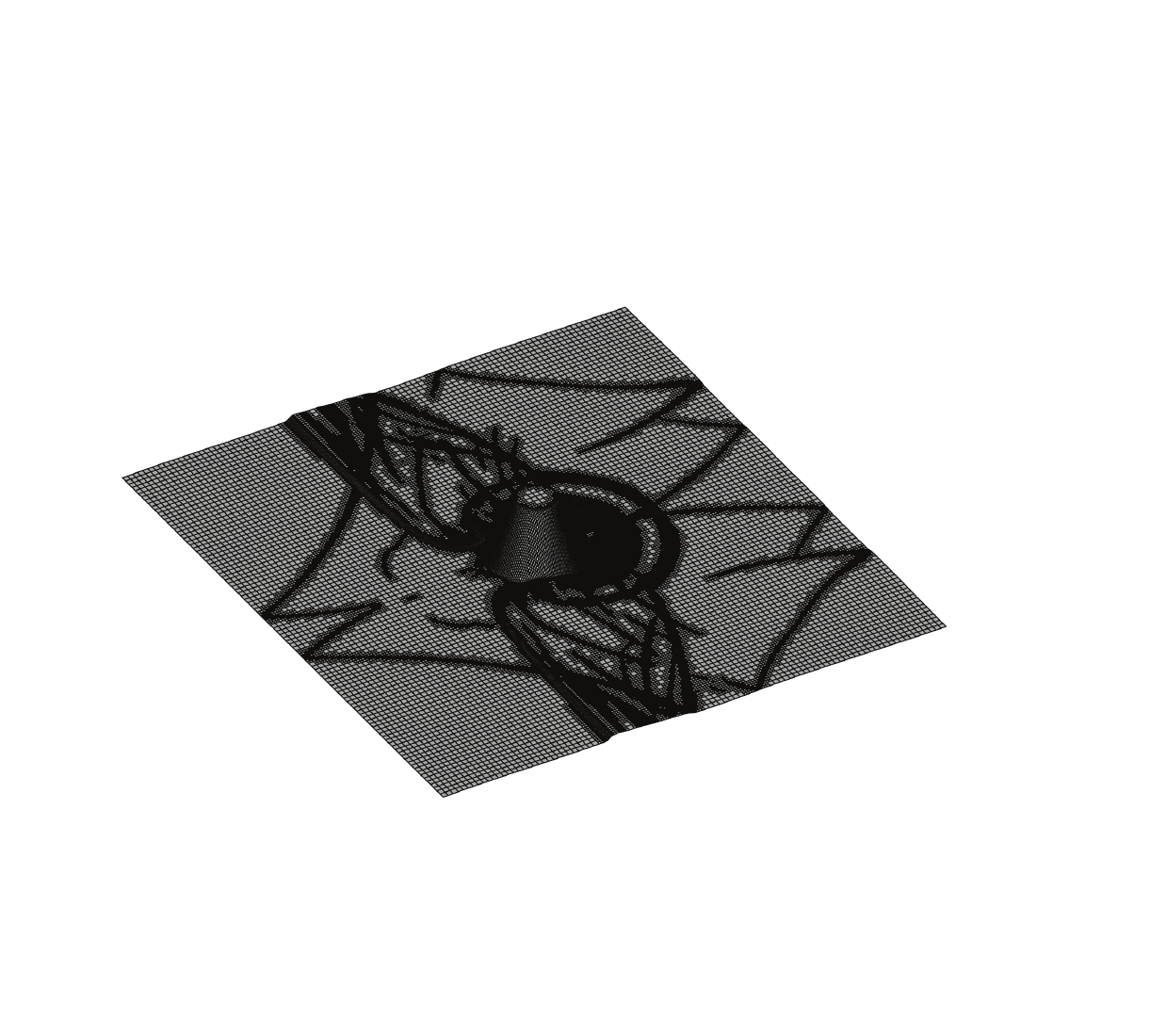}
\includegraphics[width=0.475\textwidth,trim=30 30 60 60,clip]{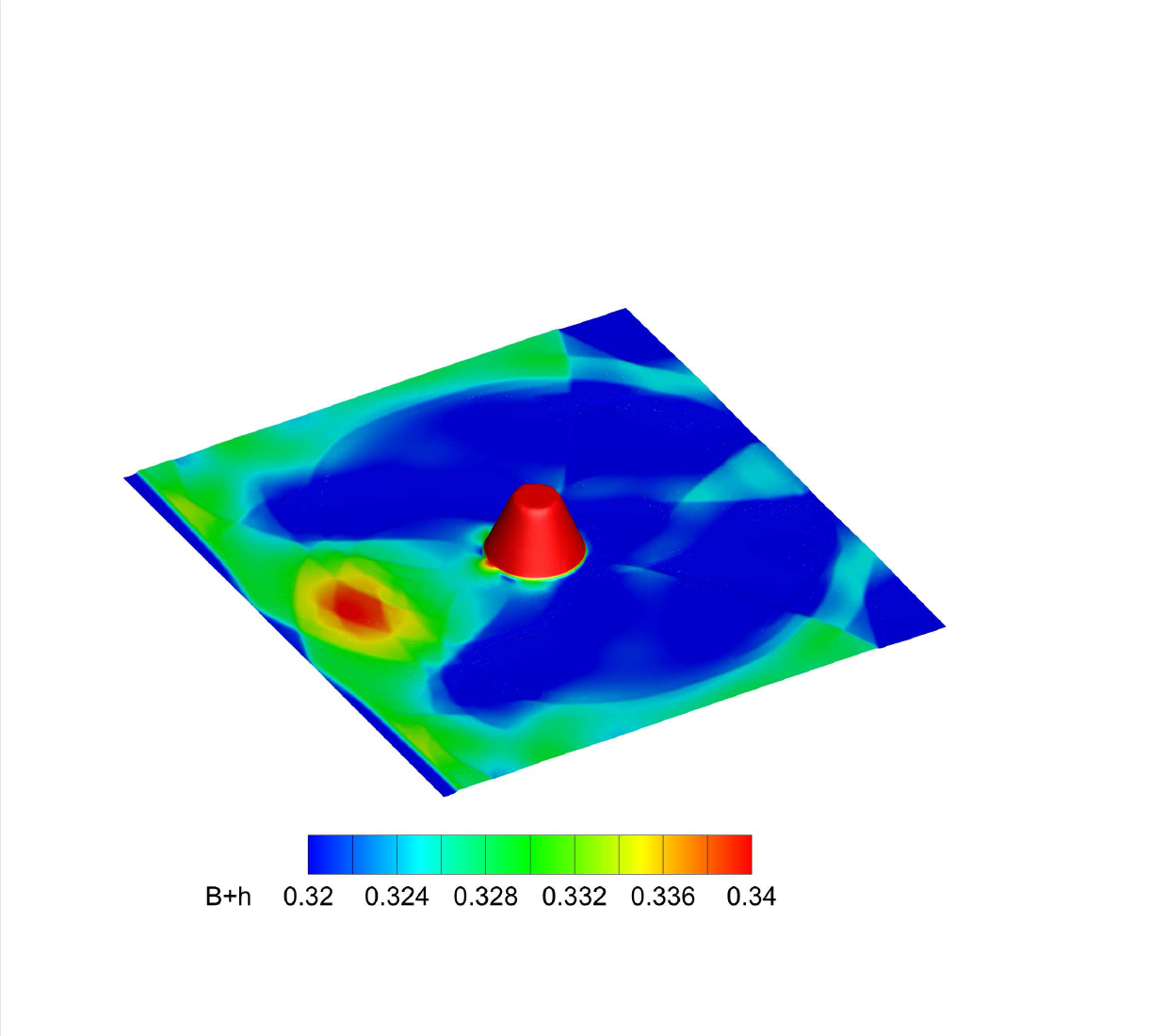}
\includegraphics[width=0.475\textwidth,trim=30 30 60 60,clip]{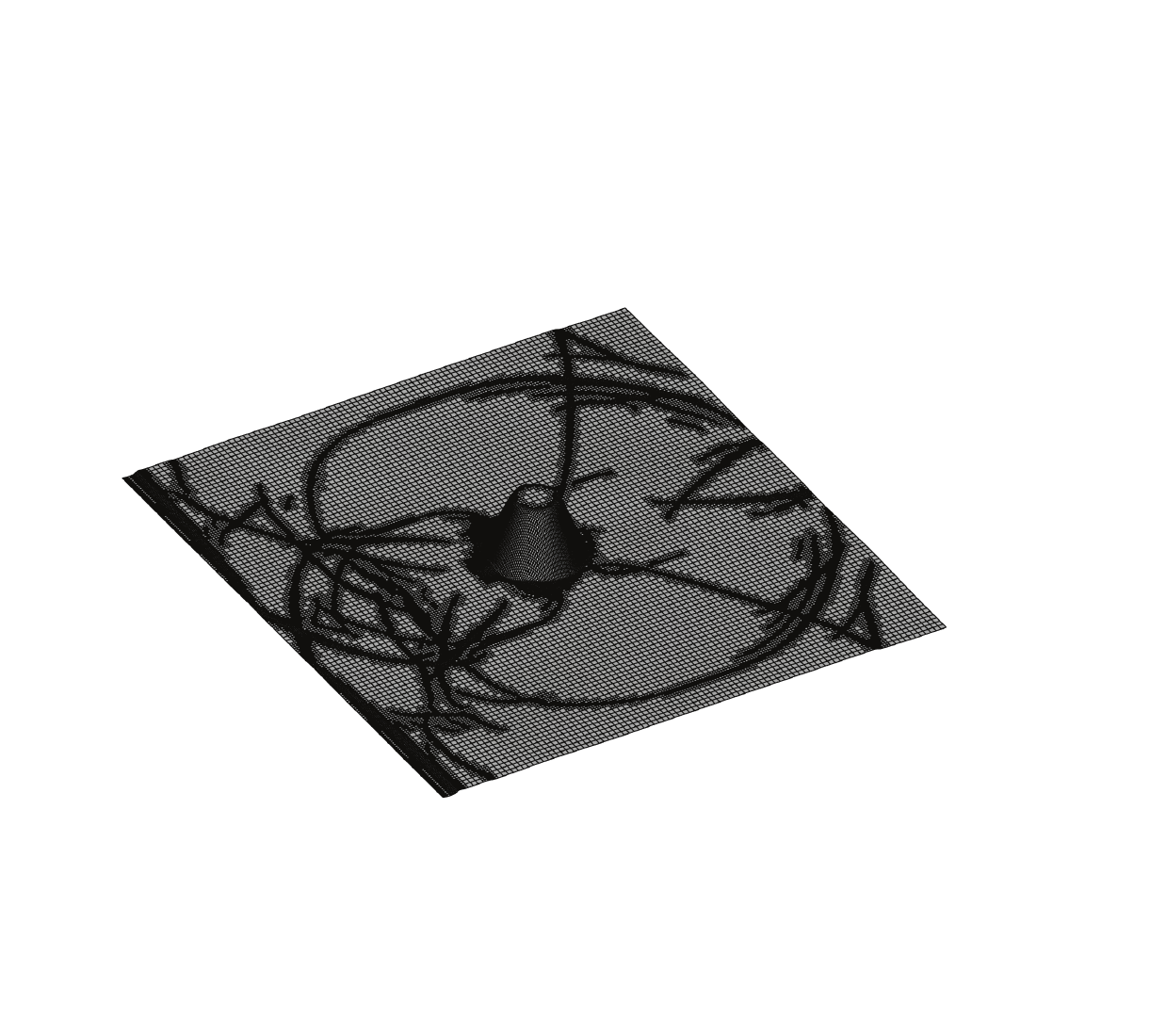}
\caption{\label{tsunami-contour}Laboratory tsunami: simulated water surface and mesh distribution at t = 0, 9 s, 15.5 s, 25 s and 30 s.}
\end{figure}

\begin{figure}[!htb]
\centering
\includegraphics[width=0.325\textwidth]{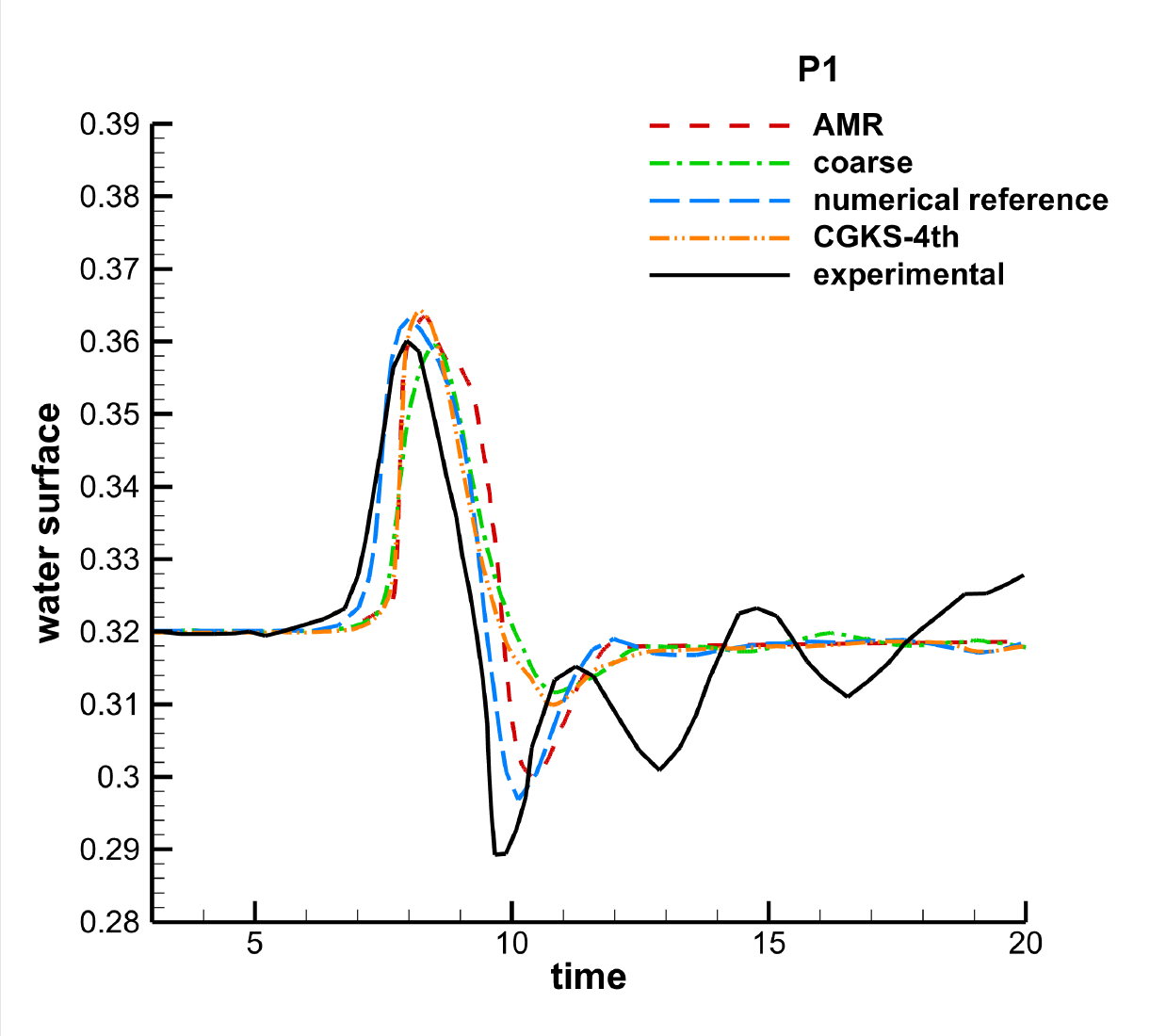}
\includegraphics[width=0.325\textwidth]{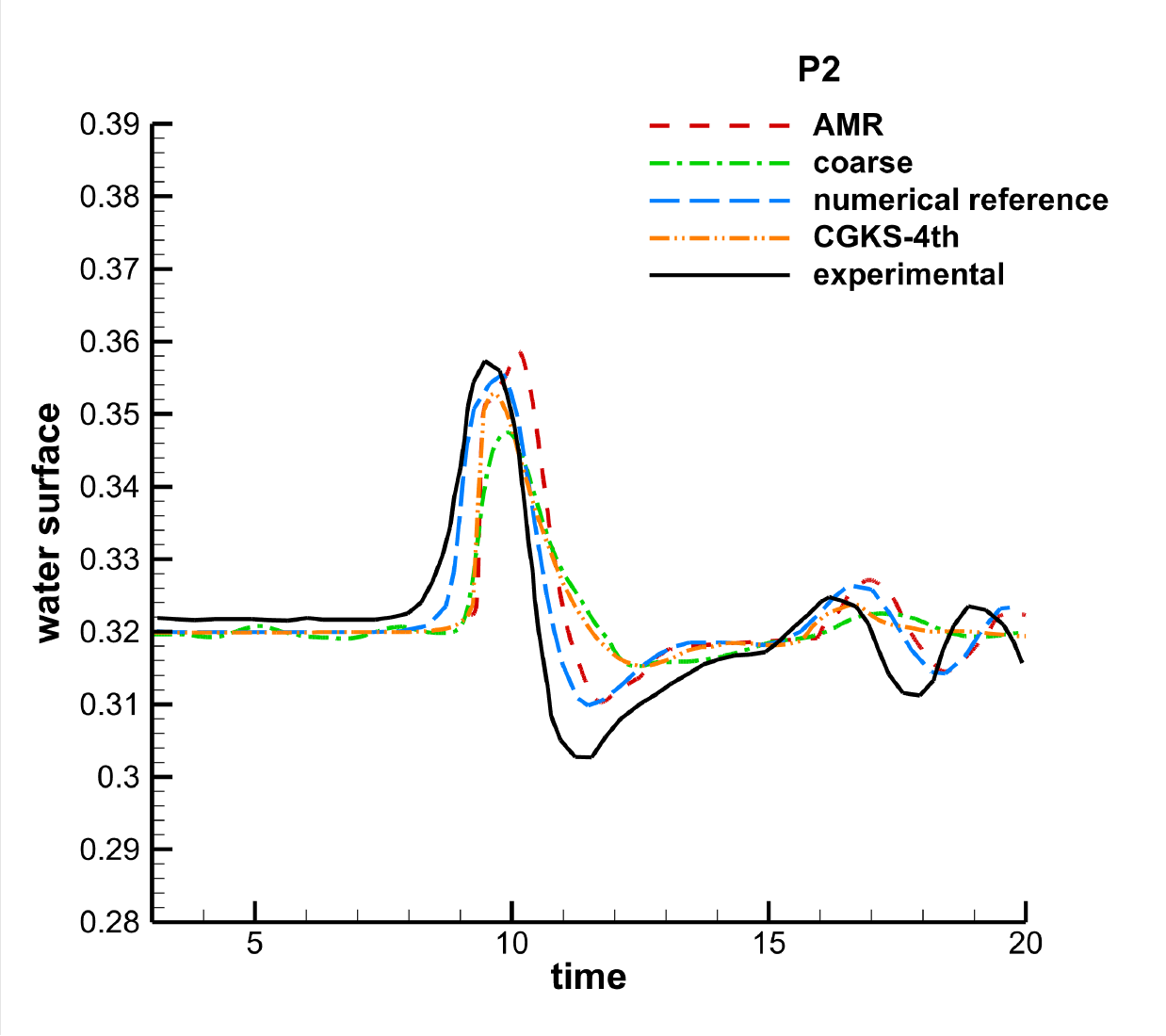}
\includegraphics[width=0.325\textwidth]{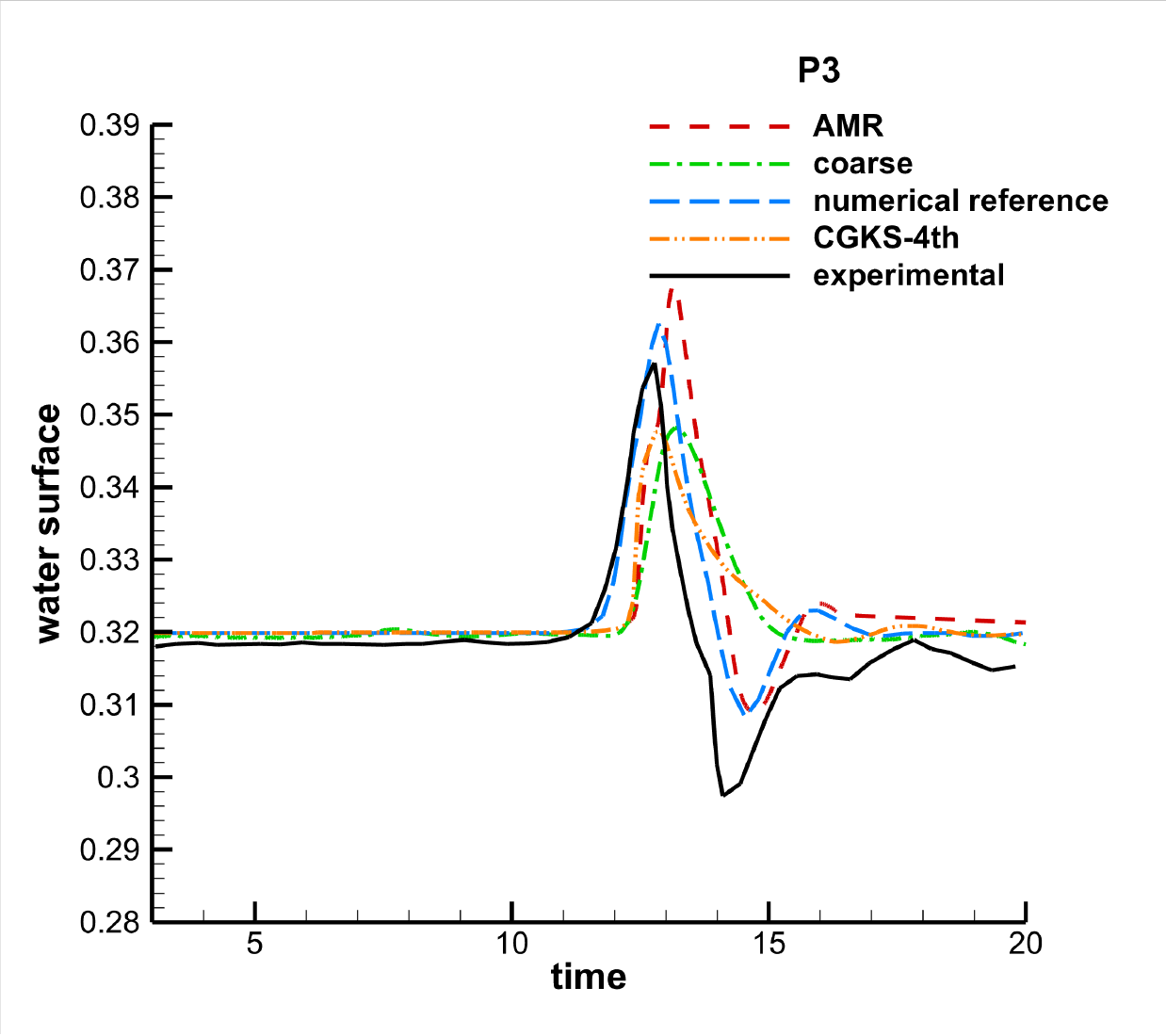}
\caption{\label{tsunami-line}Laboratory tsunami: the water surface level at three measurement gauges.}
\end{figure}

%\begin{figure}[!htb]
%\centering
%\includegraphics[width=0.3\textwidth,trim=30 10 30 10,clip]{tsunami-p1}
%\includegraphics[width=0.3\textwidth,trim=30 10 30 10,clip]{tsunami-p2}
%\includegraphics[width=0.3\textwidth,trim=30 10 30 10,clip]{tsunami-p3}
%\caption{\label{tsunami-line}Laboratory tsunami: the water surface level at three measurement gauges.}
%\end{figure}

\subsection{Laboratory tsunami wave around a conical island}

Furthermore, we use GKS-STAMR for SWE to simulate a more challenging experiment to verify the robustness and capture complex structures. The experiment involves a tsunami wave bypassing a conical island \cite{Briggs1995}. In the simulation, the island is represented by a frustum. The diameter of the upper circle is $2.2 \, \text{m}$, the diameter of the lower circle is $7.2 \, \text{m}$, and the height is $0.625 \, \text{m}$. The calculation domain is $[0, 26 \, \text{m}] \times [0, 27.6 \, \text{m}]$. A still water level of $H_0 = 0.32 \, \text{m}$ is assumed as the initial condition, and the inlet wave is set as a boundary condition for water level $H$ as follows:
\begin{equation*}
H= H_0+Asech ^2[\frac{C_1(t-T)}{C_2}],
\end{equation*}
where
\begin{equation*}
C_1=\sqrt{gH_0}(1+\frac{A}{2H_0}), C_2=H_0\sqrt{\frac{4H_0C_1}{3A\sqrt{gH_0}}}.
\end{equation*}
Let $A$ be the solitary wave amplitude and $T$ the time at which the wave peak enters the domain, with $A = 0.032 \, \text{m}$ and $T = 2.84 \, \text{s}$ for the case considered in this work. The other boundaries are solid walls. Two types of meshes are tested: a coarse uniform mesh with $100 \times 100$ cells and an AMR mesh that starts with the coarse mesh. The maximum refinement level is set to $l_{\text{max}} = 2$, with the refinement criterion defined as $\beta = |\nabla (B+h)|/|\nabla (B+h)|_{\text{max}}$, and the threshold is set at 0.01.

Fig. \ref{tsunami-contour} illustrates the water level at selected times. It is evident that when the tsunami reaches the island, the dry-wet interface changes dramatically. The water surface at the front end of the island decreases, while the water levels on both sides increase significantly. Following the tsunami flow, a low water surface area forms behind the island, as the water accumulated by the bypass flow is transmitted further. After the tsunami impacts the back wall, it reflects and strikes the island again, resulting in a more complex structure. Our scheme demonstrates robustness and effectively captures this dynamic complexity.

Fig. \ref{tsunami-line} presents a comparison between numerical and experimental data for the three water level gauges. The numerical results include those by 4th-order CGKS \cite{zhao2021-swe}, as well as numerical simulation results from the literature \cite{Garc2019}. It is evident that the coarse mesh fails to accurately capture the main water level peaks and reflections, while the AMR shows good agreement with both the experimental and numerical results in areas of local refinement.

\begin{figure}[!htb]
\centering
\includegraphics[width=0.45\textwidth,trim=30 10 30 10,clip]{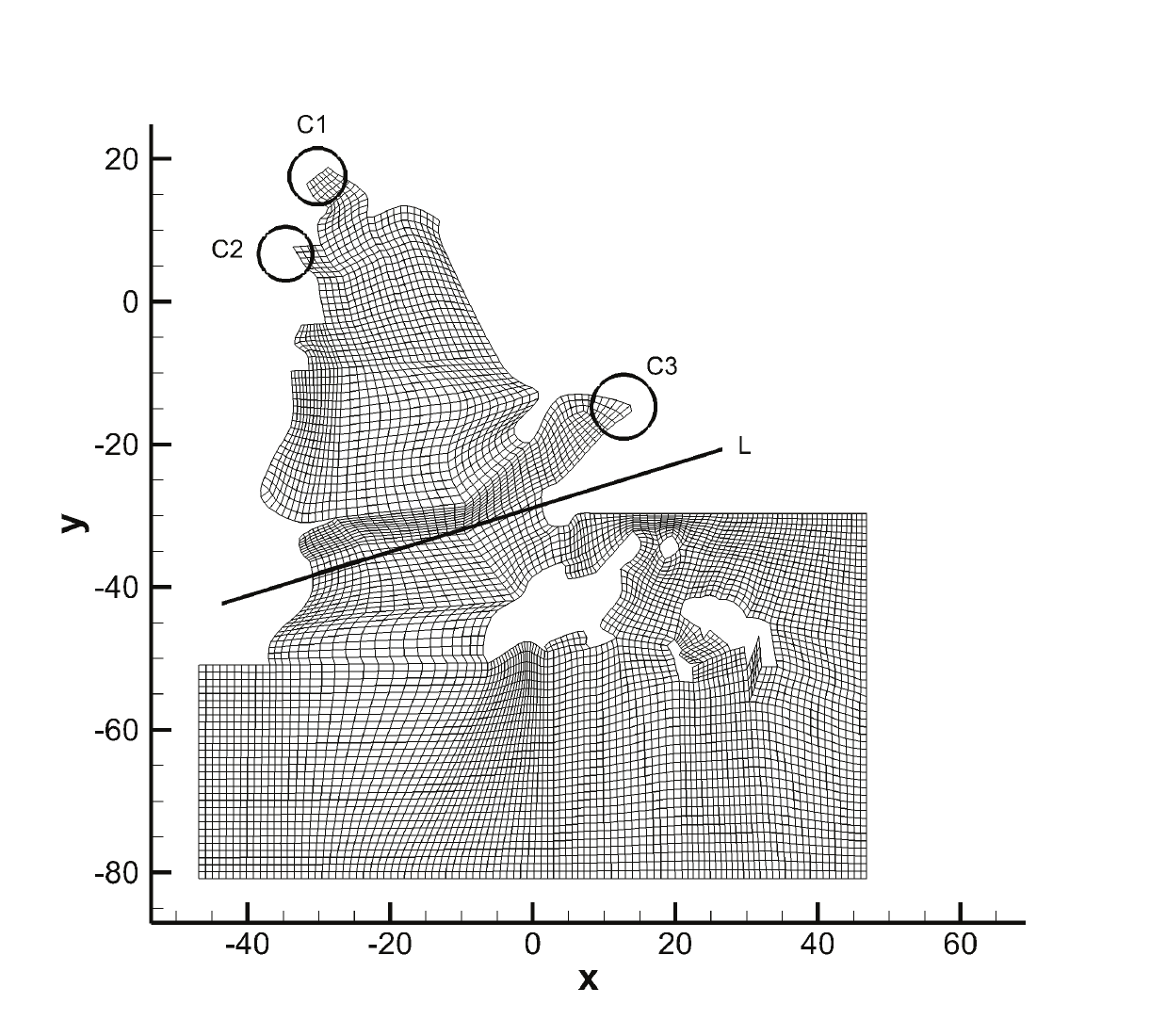}
\includegraphics[width=0.45\textwidth,trim=30 10 30 10,clip]{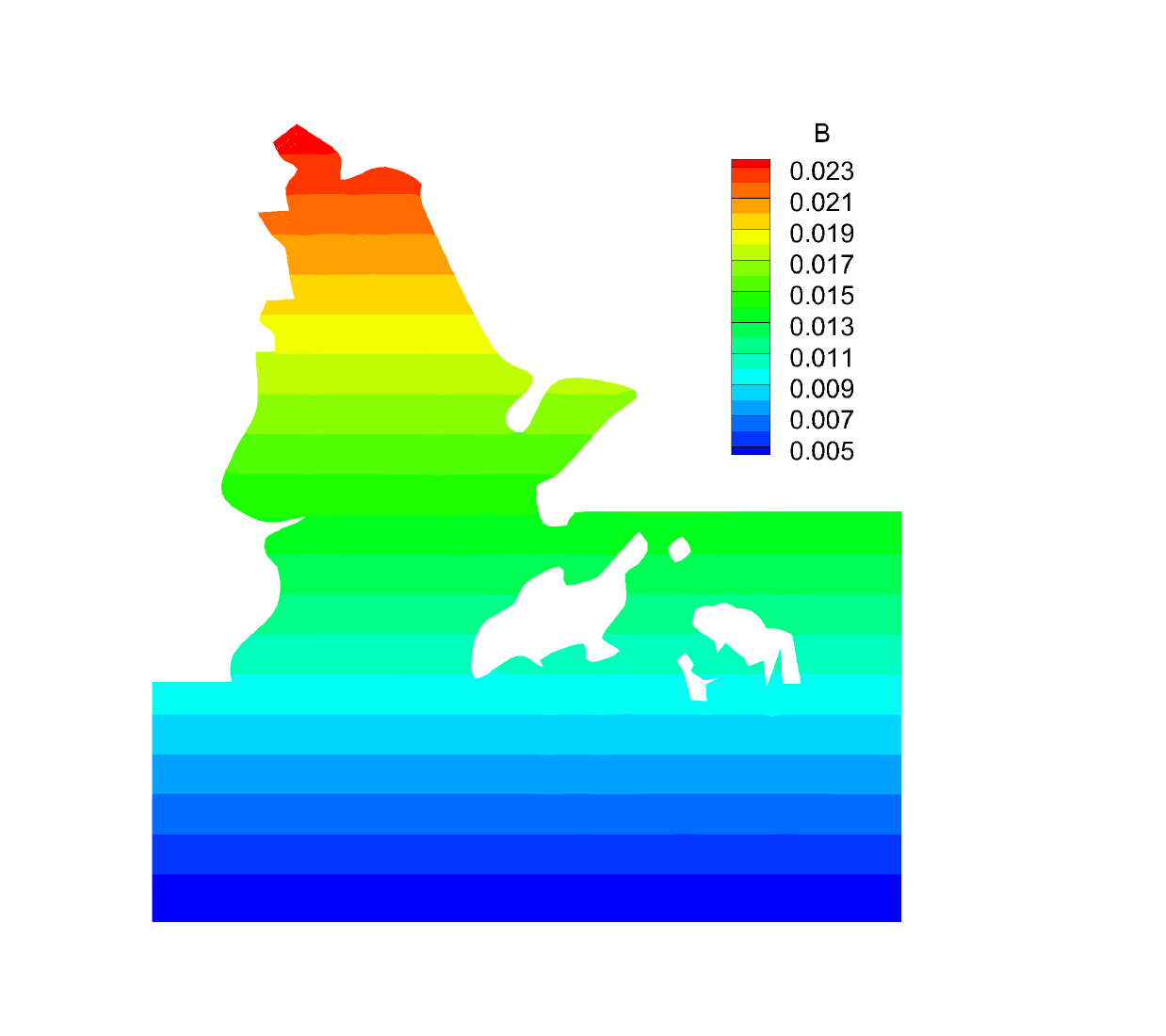}
\caption{\label{pearl-river-mesh} Flow in an estuary: the computational mesh and the bottom topography.}
\end{figure}

\begin{figure}[!htb]
\centering
\includegraphics[width=0.45\textwidth,trim=30 10 30 10,clip]{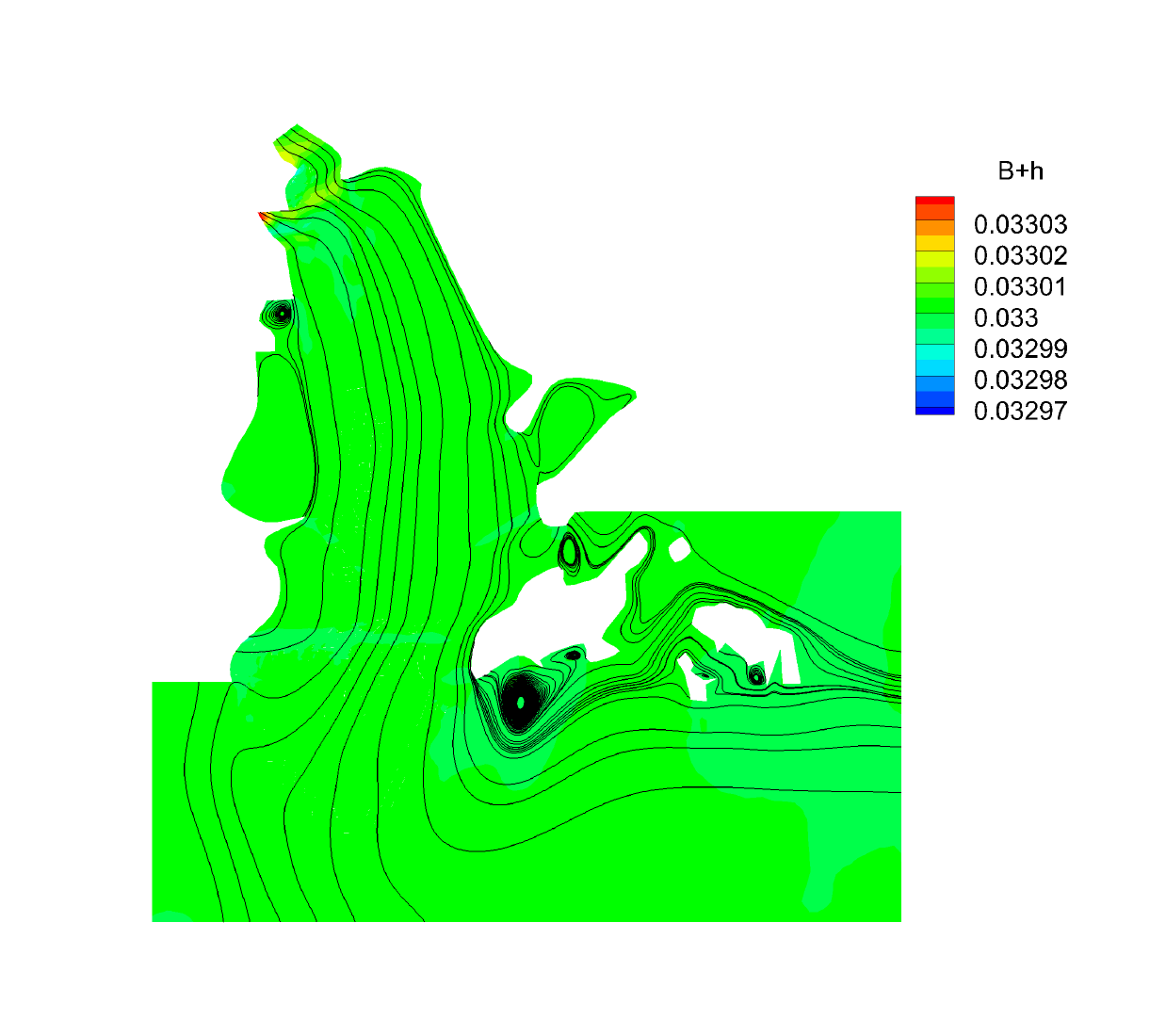}
\includegraphics[width=0.45\textwidth,trim=30 10 30 10,clip]{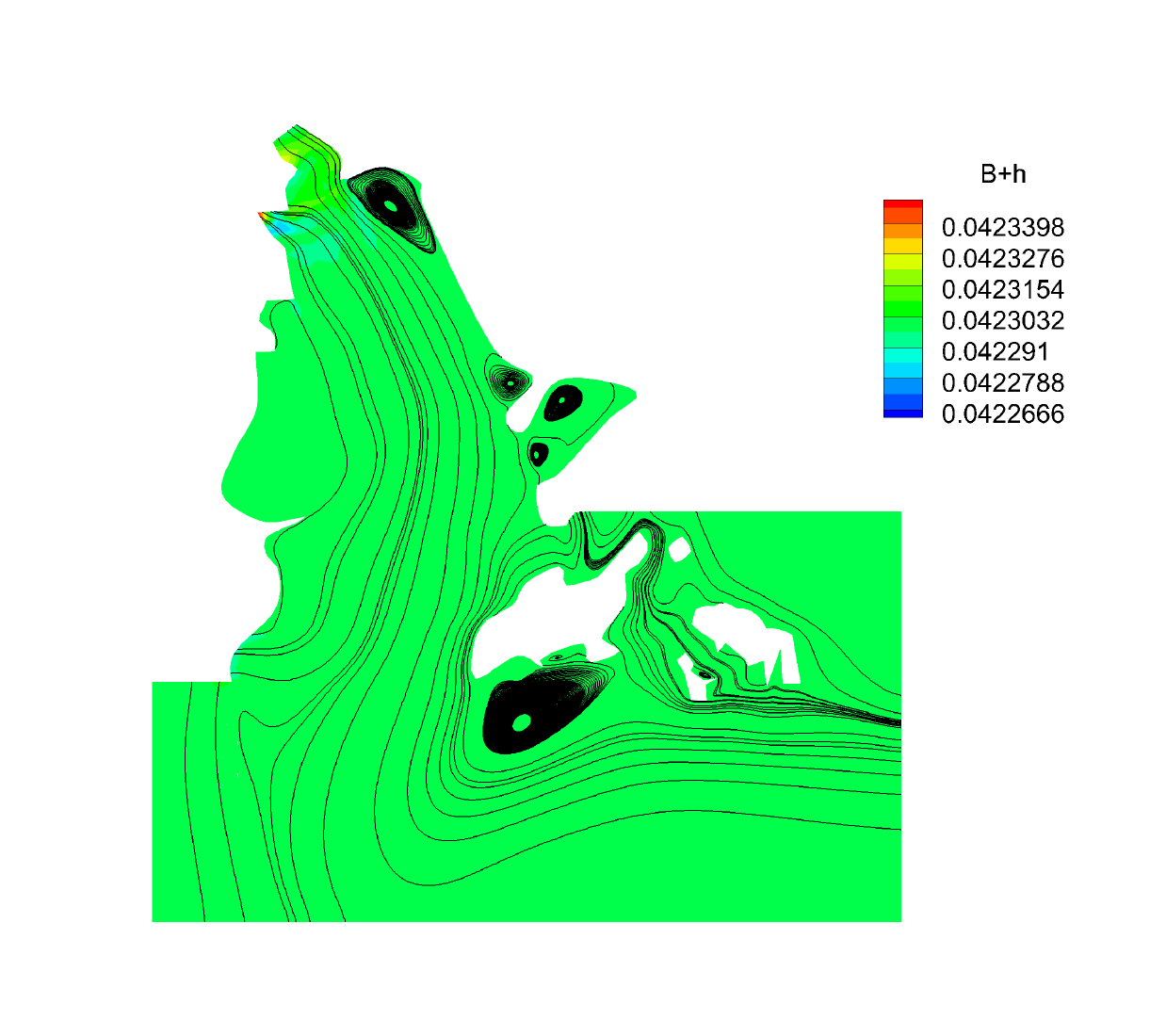}
\caption{\label{pearl-river-compare} Flow in an estuary: streamlines and contours of water surface ($B+h$) over flat and non-flat bottom topography.}
\end{figure}

\begin{figure}[!htb]
\centering
\includegraphics[width=0.45\textwidth,trim=30 20 30 20,clip]{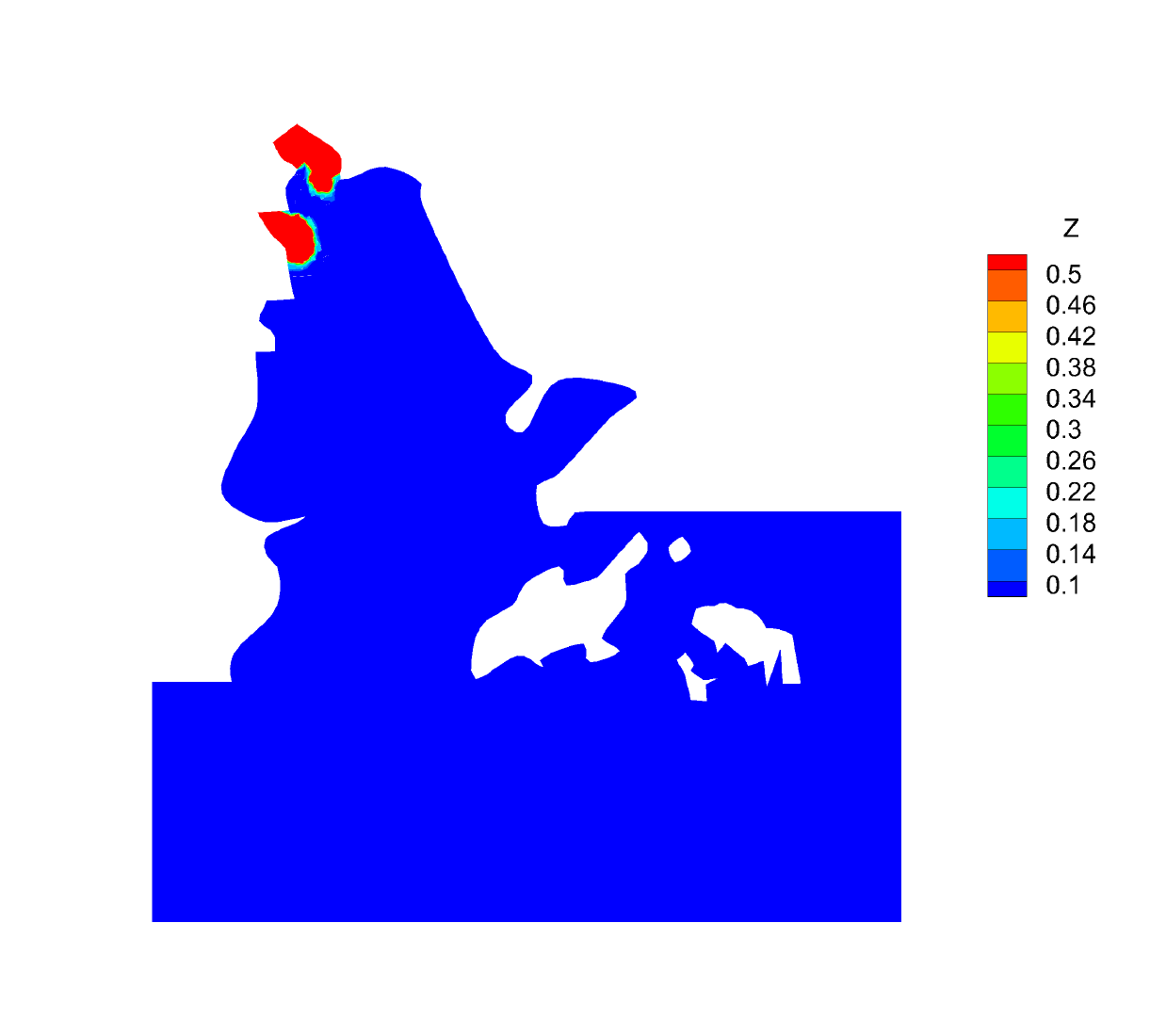}
\includegraphics[width=0.45\textwidth,trim=30 20 30 20,clip]{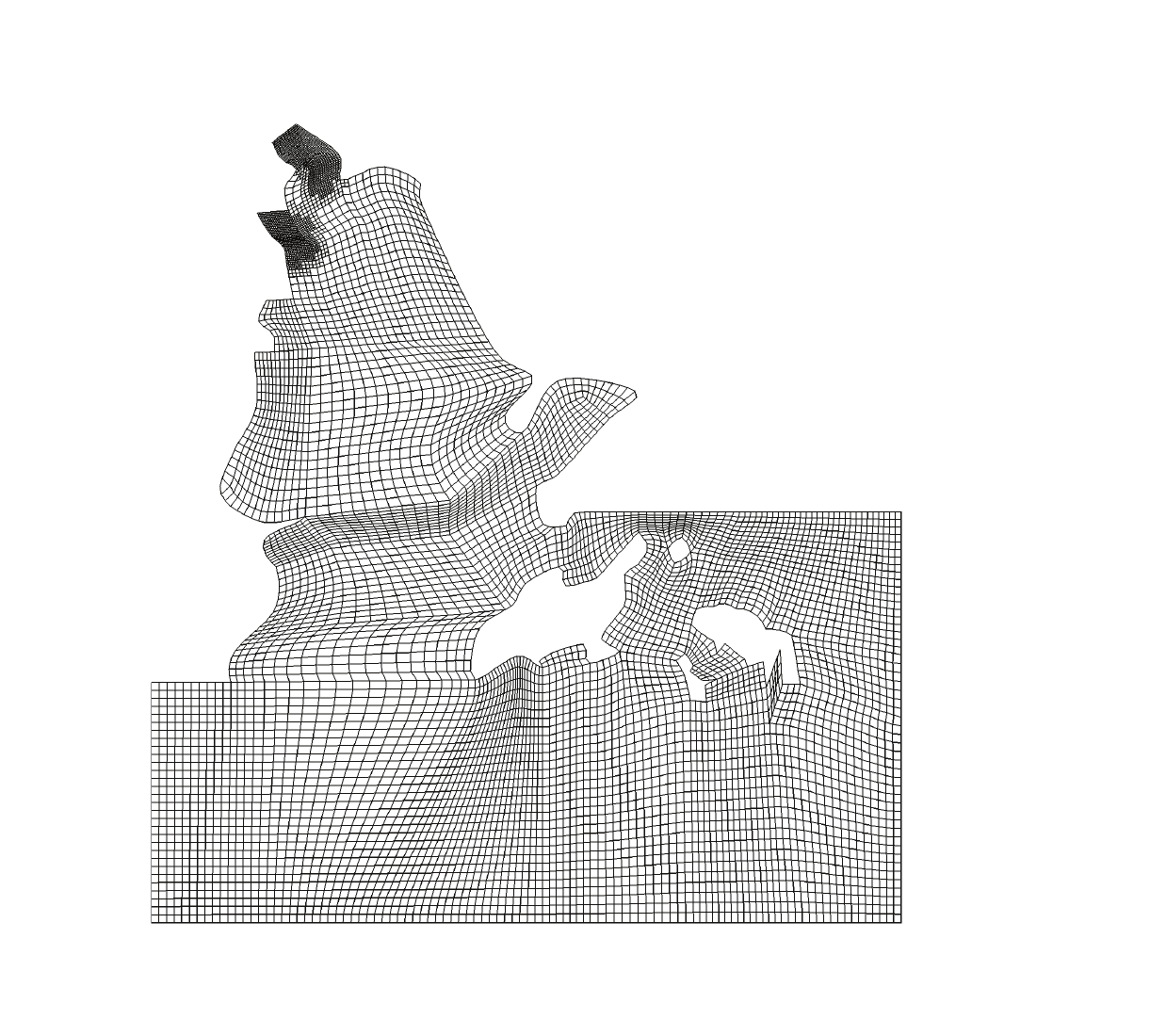}
\includegraphics[width=0.45\textwidth,trim=30 20 30 20,clip]{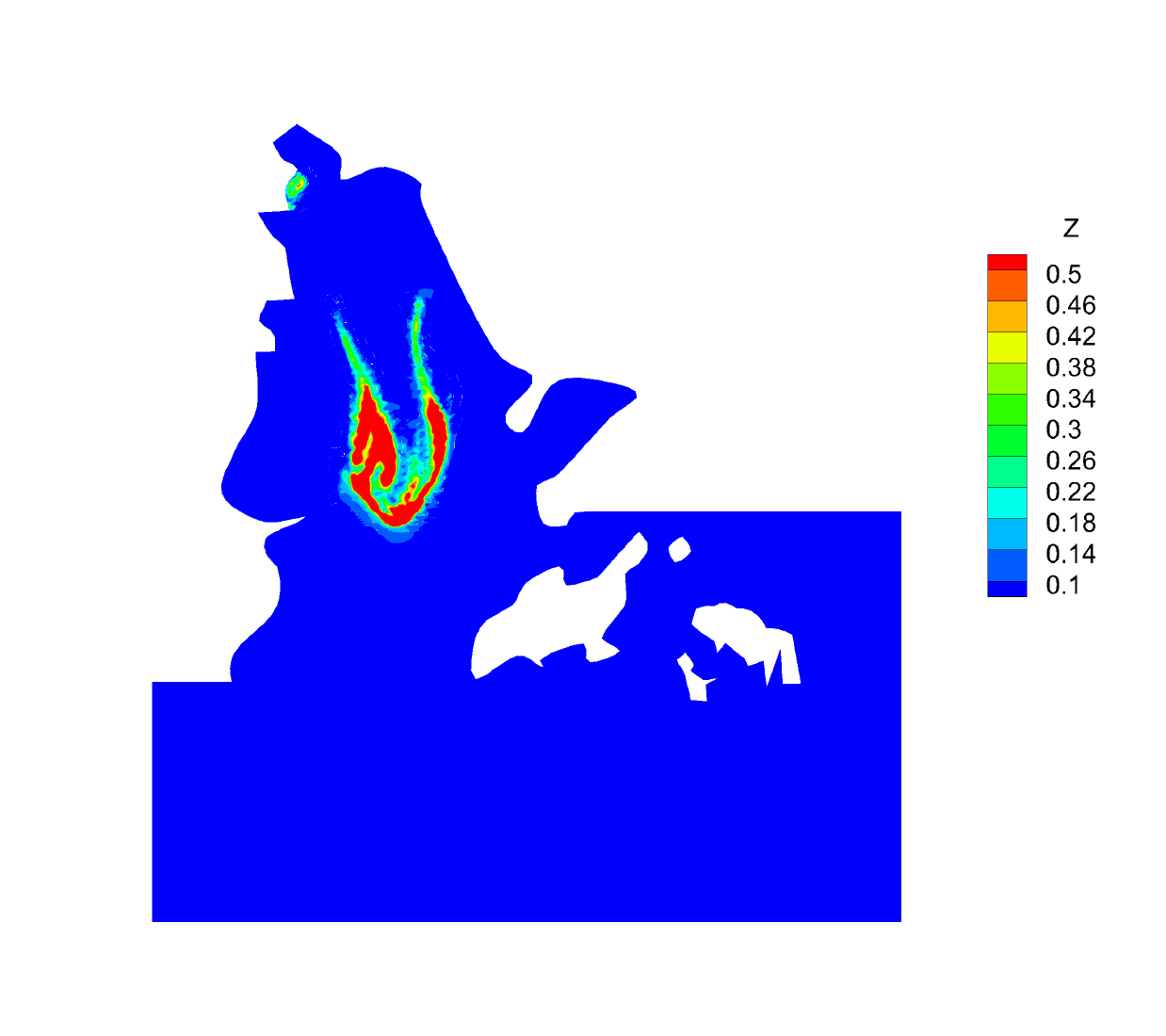}
\includegraphics[width=0.45\textwidth,trim=30 20 30 20,clip]{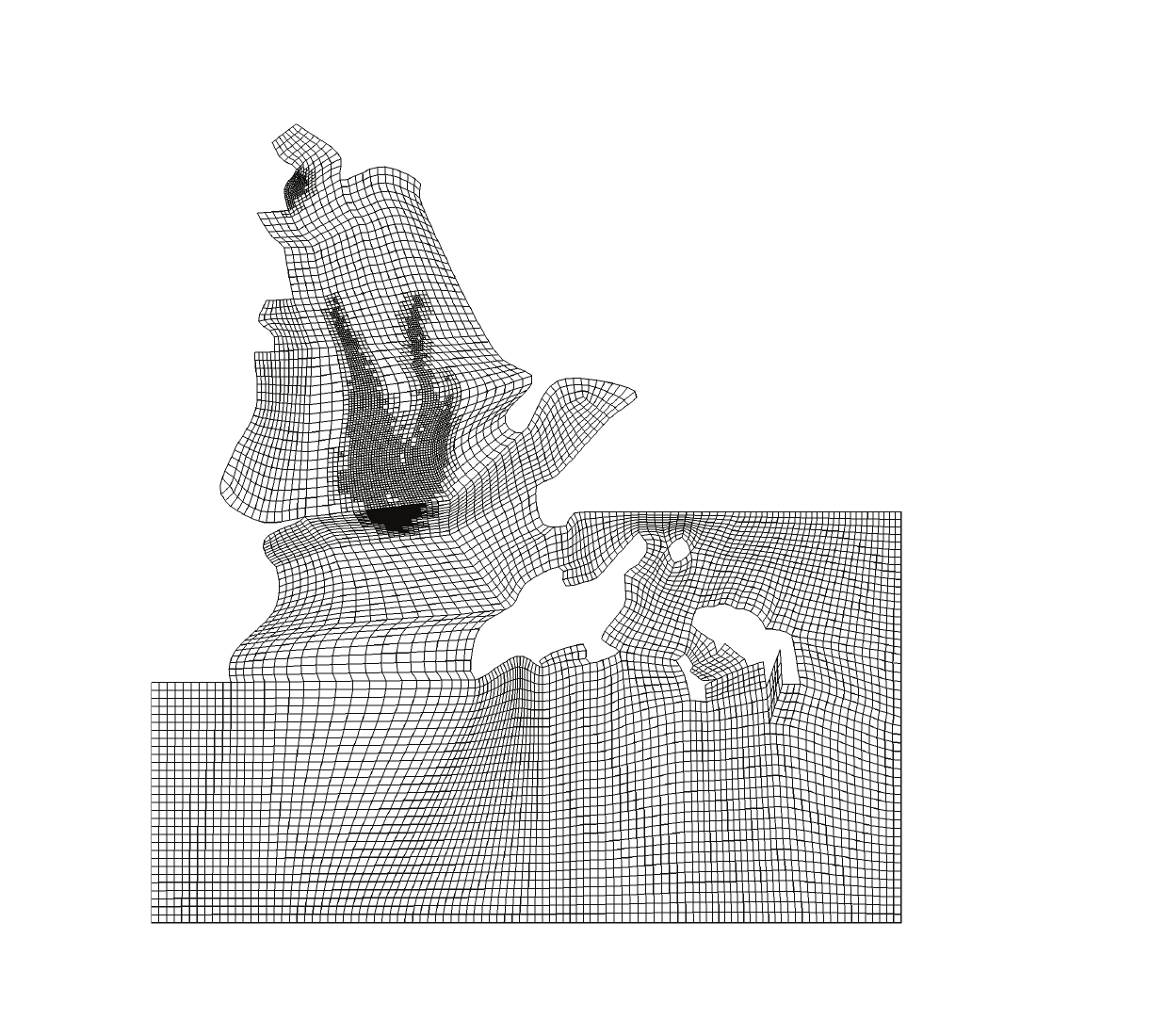}
\includegraphics[width=0.45\textwidth,trim=30 20 30 20,clip]{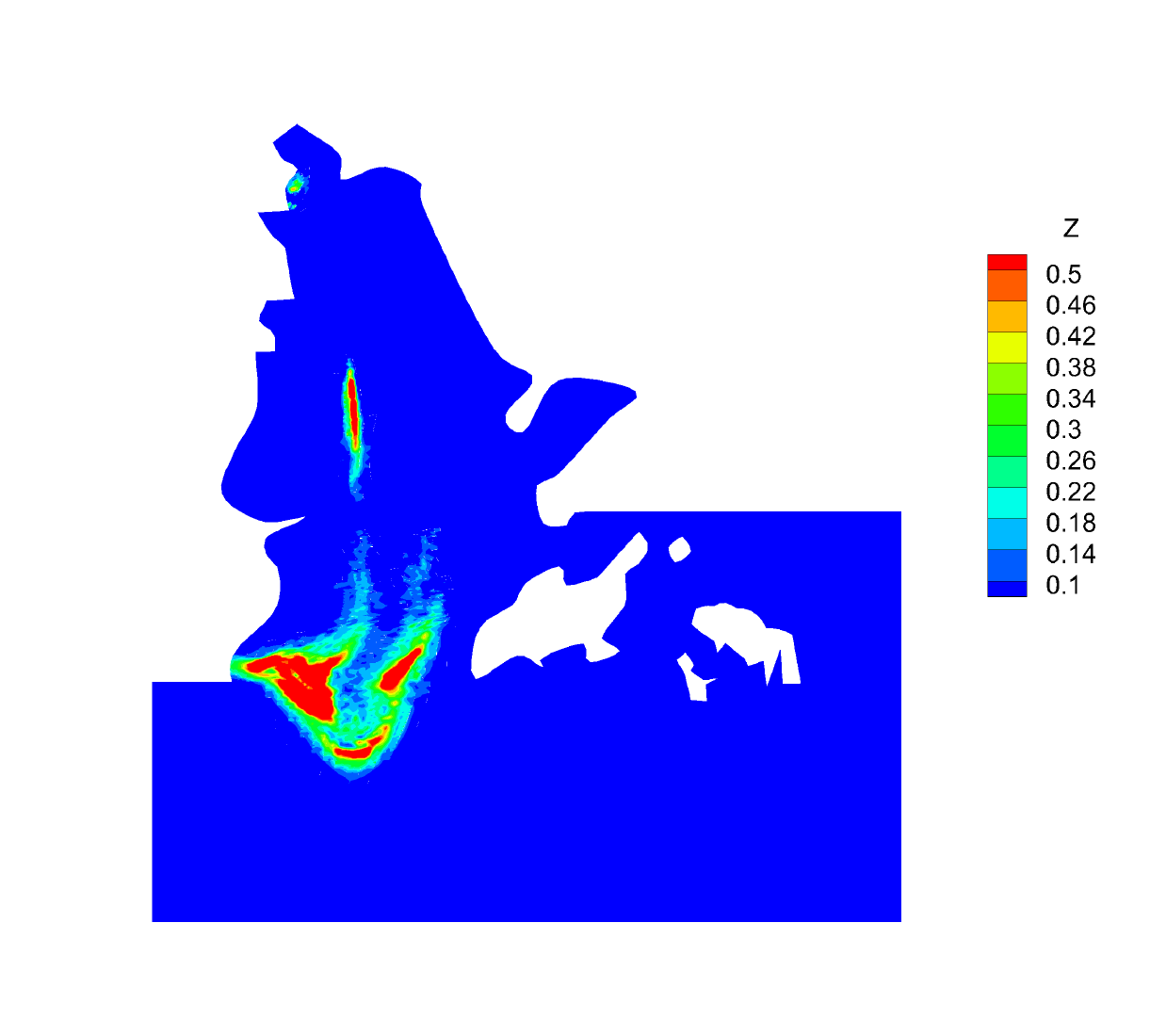}
\includegraphics[width=0.45\textwidth,trim=30 20 30 20,clip]{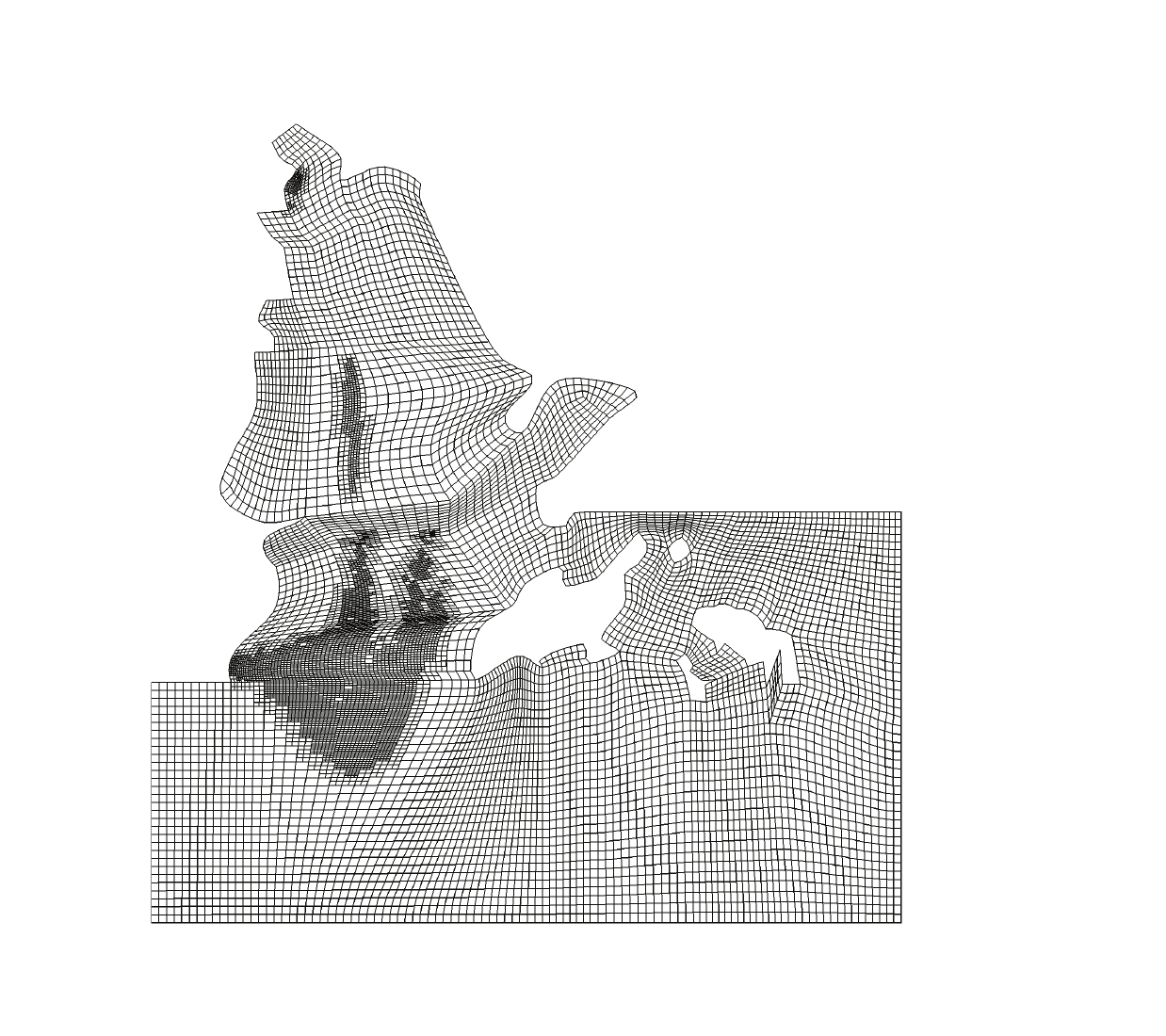}
\caption{\label{pearl-river-flat} Flow in an estuary: pollutant contours and mesh over flat bottom topography.}
\end{figure}

\begin{figure}[!htb]
\centering
\includegraphics[width=0.45\textwidth,trim=30 20 30 20,clip]{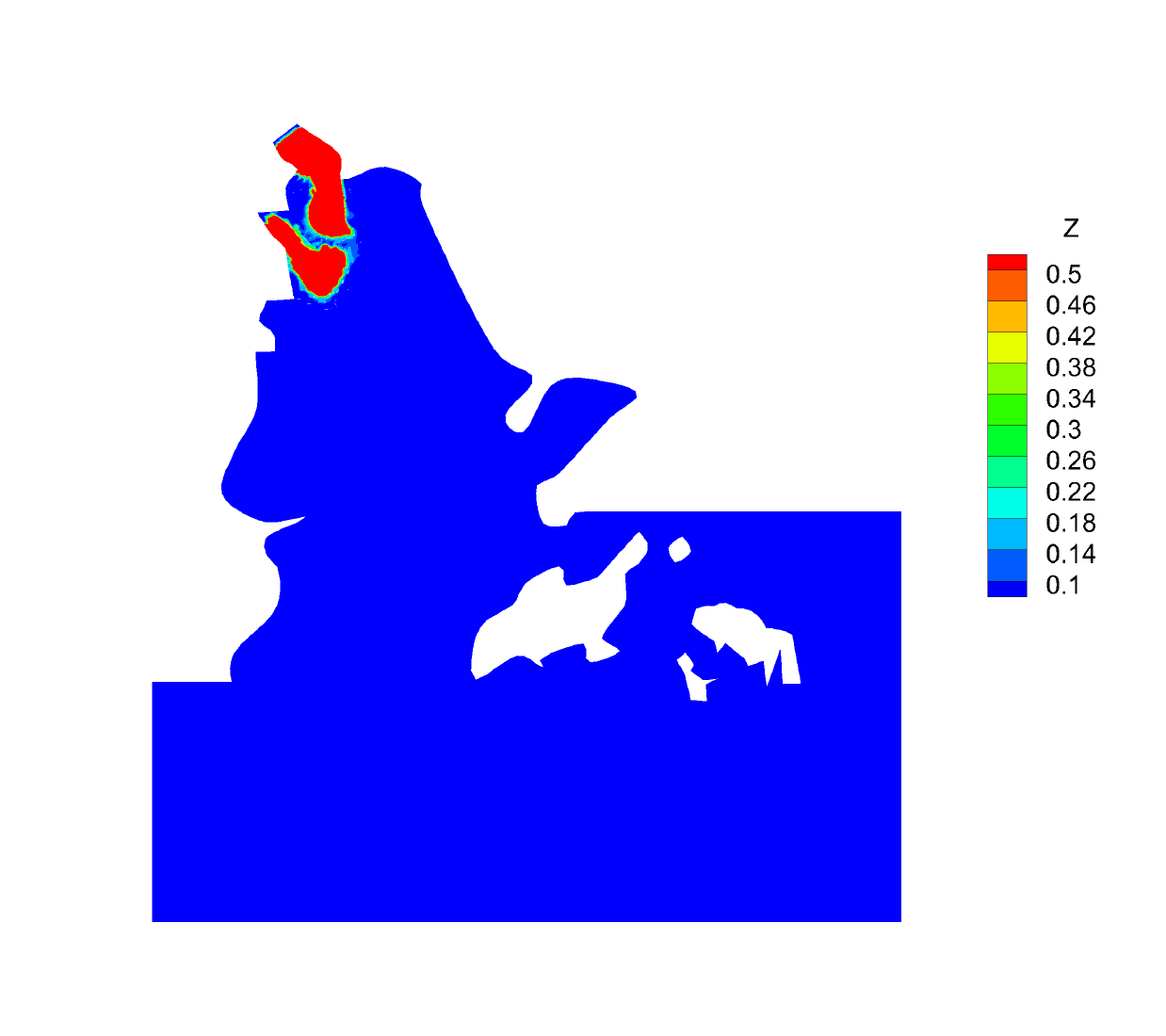}
\includegraphics[width=0.45\textwidth,trim=30 20 30 20,clip]{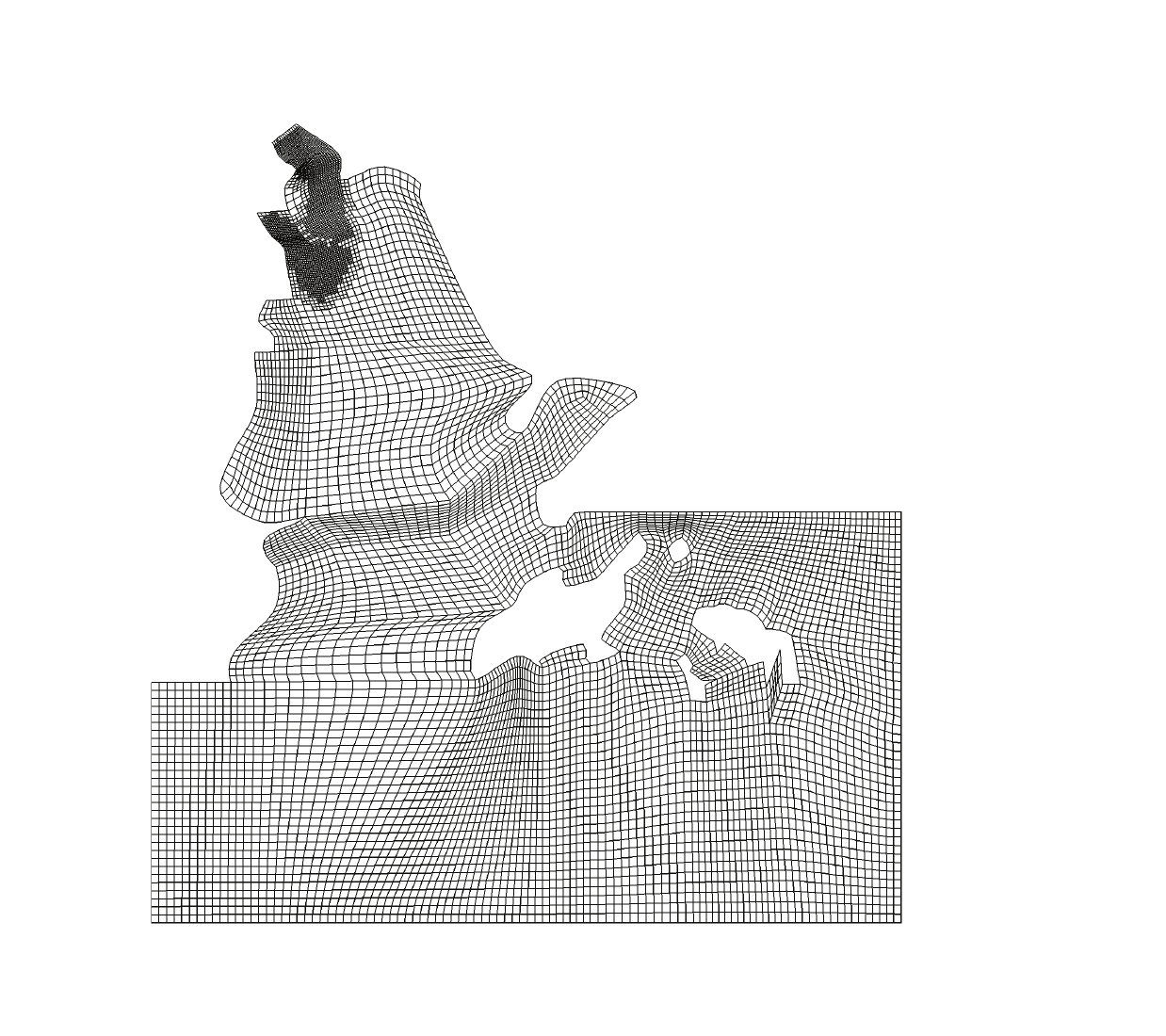}
\includegraphics[width=0.45\textwidth,trim=30 20 30 20,clip]{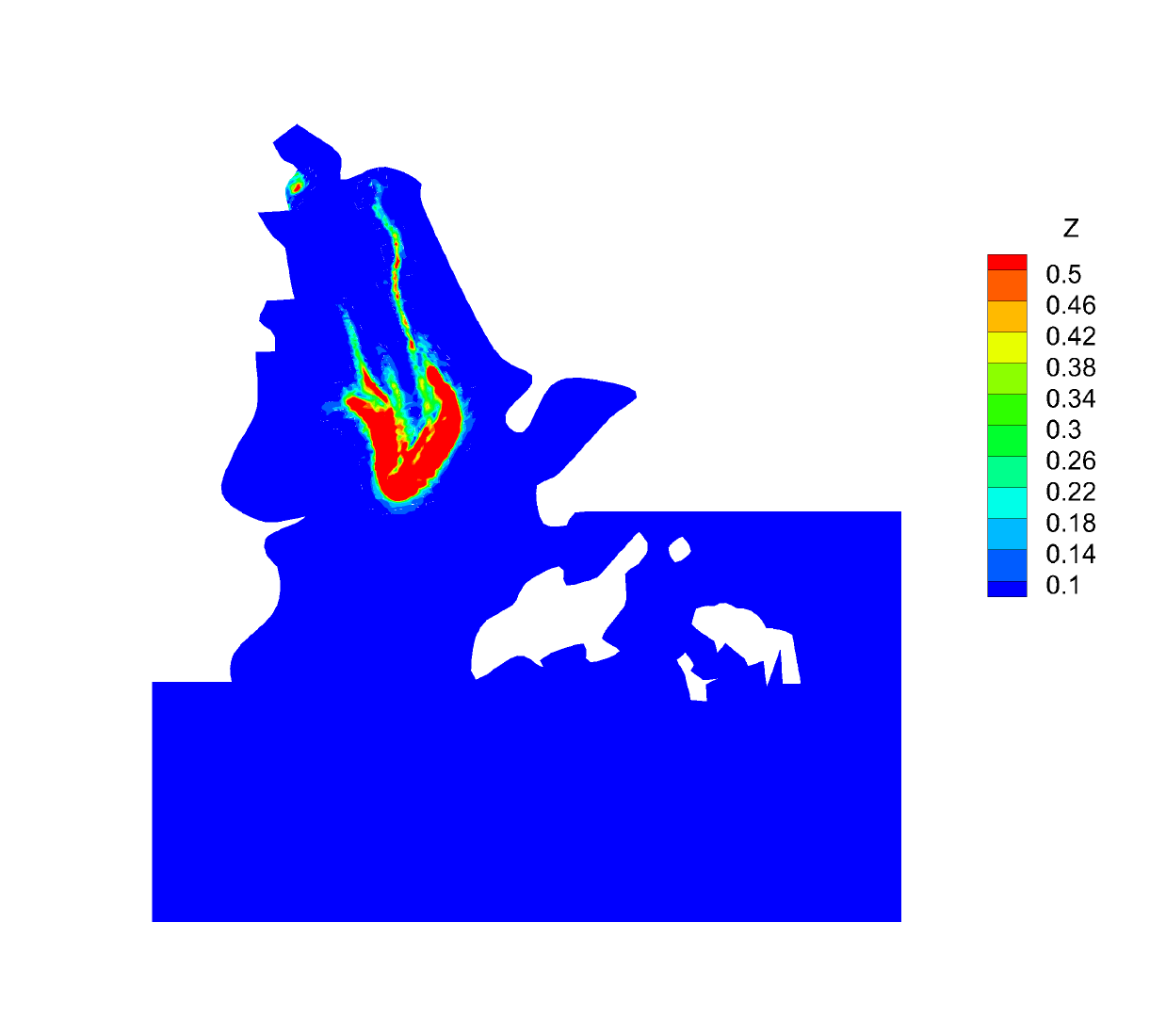}
\includegraphics[width=0.45\textwidth,trim=30 20 30 20,clip]{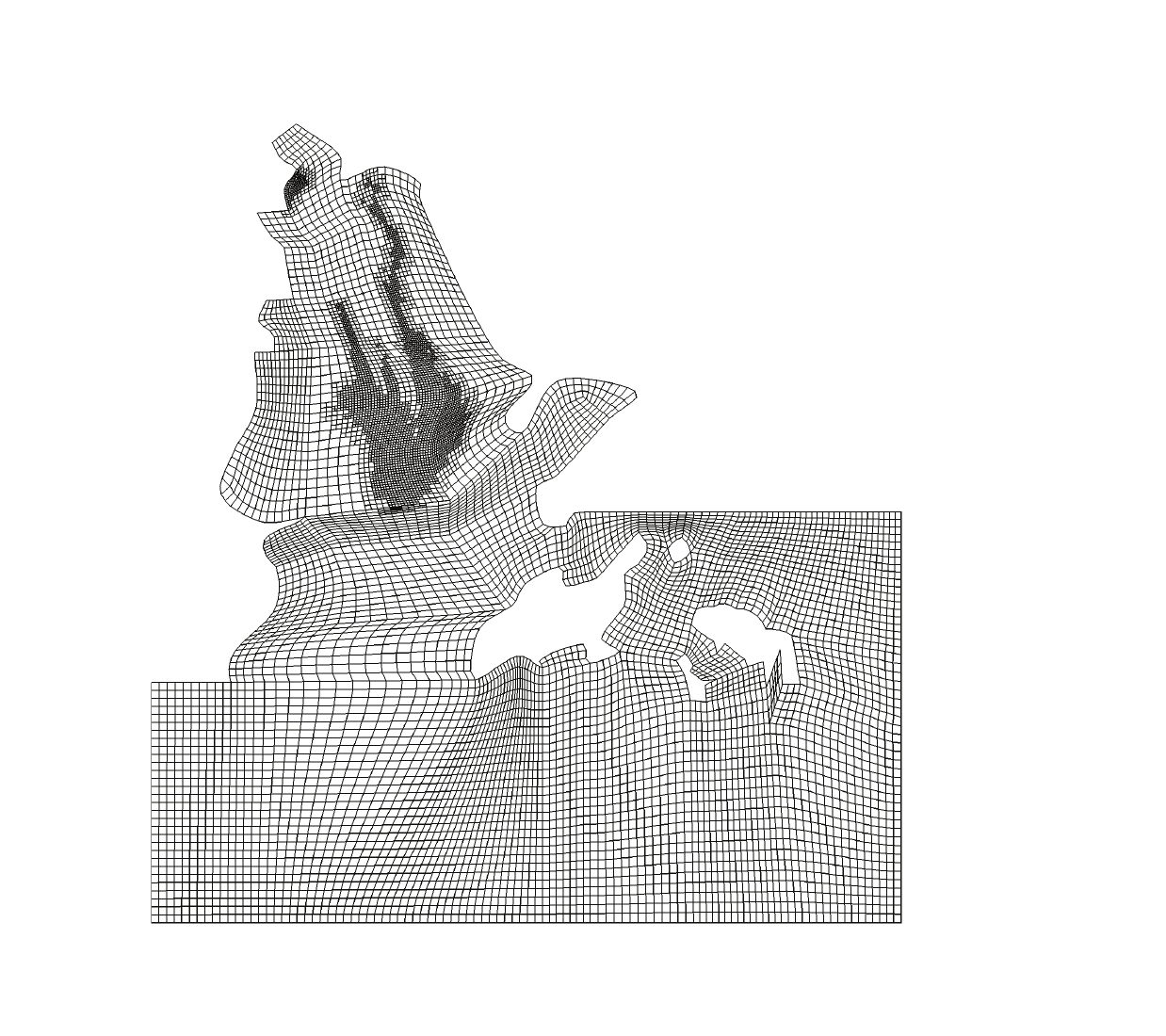}
\includegraphics[width=0.45\textwidth,trim=30 20 30 20,clip]{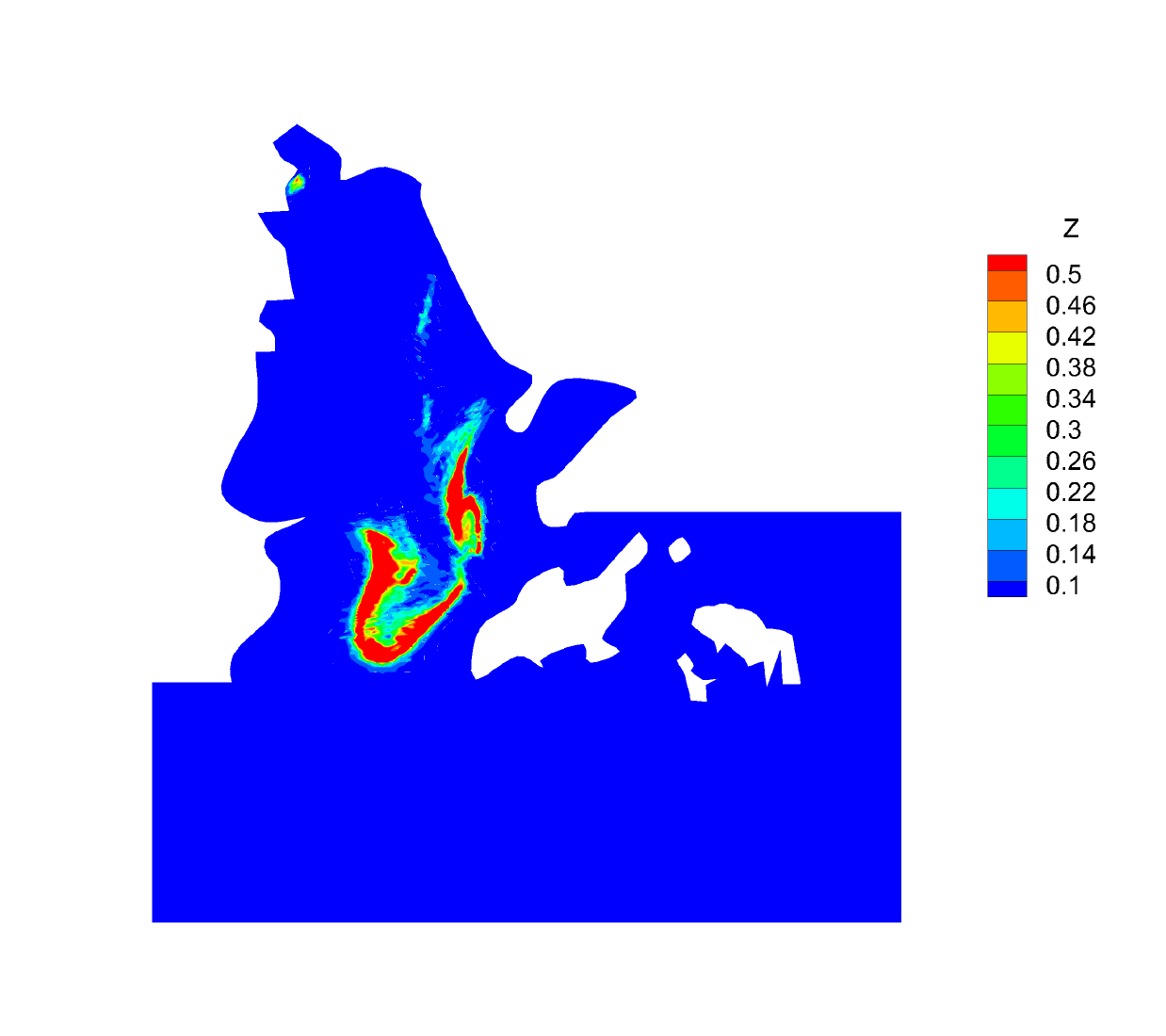}
\includegraphics[width=0.45\textwidth,trim=30 20 30 20,clip]{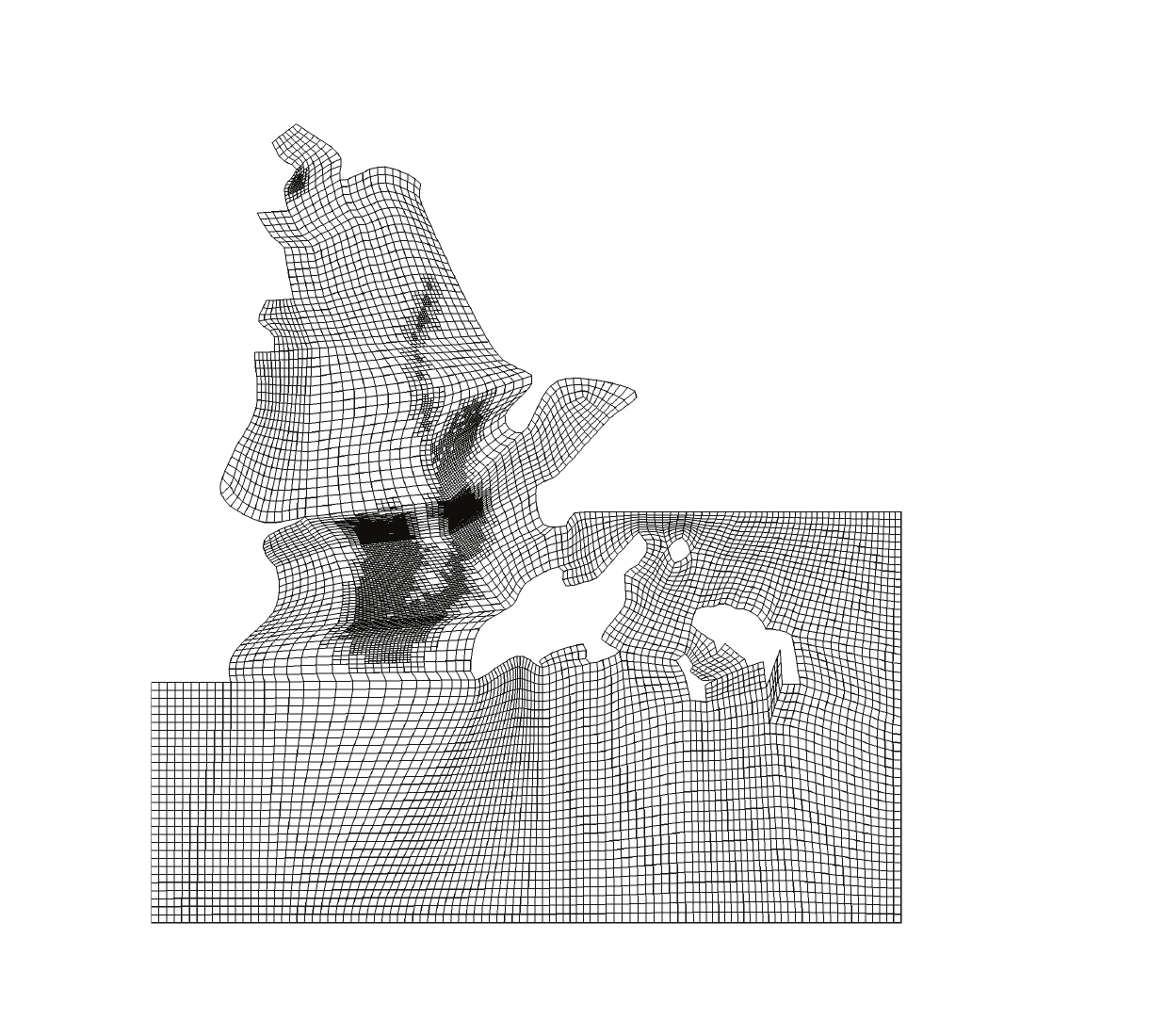}
\caption{\label{pearl-river-nonflat} Flow in an estuary: pollutant contours and mesh over non-flat bottom topography.}
\end{figure}

\subsection{Flow in an estuary}
In order to further verify the application potential of the current GKS-STAMR for SWE in complex real flows, the shallow water flow with pollutant transport in an estuary is simulated.
The estuary and the computational mesh are shown in Fig.~\ref{pearl-river-mesh}, where 6632 cells are used as the base mesh. This estuary refers to the Pearl River estuary, which has two river inlets (inlet C1 and inlet C2) in the upper reaches and one channel outlet (outlet C3) on the right side of the middle of the estuary. The boundary below the line $L: y = 0.3894x - 31.43$ within the domain is also treated as an outlet boundary. There are some islands in the estuary area, represented as the white regions in the computational domain.

We also test the influence of bottom topography, considering two cases: a simple artificial bottom topography, with the contour shown on the left of Fig. \ref{pearl-river-mesh}. The discharges at the two inlets are set to 0.003. The water surface heights are set to $0.033$ for the flat bottom and $0.0423 - B(x, y)$ for the non-flat bottom. Wall boundary conditions are applied to the other boundaries. Pollutants are emitted from the inlet between $t = 10^4$ and $t = 1.1 \times 10^4$. In order to trace pollutants, the criterion is selected as $\beta=|\nabla hZ|/|\nabla hZ|_{max}$, and the threshold is set at 0.1.

In Fig. \ref{pearl-river-compare}, the streamlines and water heights are presented. For the flat bottom, the contours of pollutants and mesh distribution are shown in Fig. \ref{pearl-river-flat}; for the non-flat bottom, the contours of pollutants and mesh distribution are given in Fig. \ref{pearl-river-nonflat}. In both cases, as the pollutants flow, the mesh adapts according to their movement. For this unsteady long-term flow, the current AMR significantly demonstrates its advantages, improving the accuracy of pollutant capture with fewer mesh cells.
Moreover, when considering the influence of the bottom topography, a comparison between Fig. \ref{pearl-river-nonflat} and Fig. \ref{pearl-river-flat} reveals that the non-flat bottom alters the acceleration of nearby fluids, especially near the shore, causing pollutants to accumulate near the islands.

\section{Conclusion}

In this study, a well-balanced GKS for the SWE is developed and integrated with local time stepping within a STAMR framework. To preserve the well-balanced property on adaptive meshes, quadrilateral meshes with hanging nodes are managed through subdivision into triangular subcells and a tailored source terms treatment.
Numerical experiments demonstrate that the proposed method maintains the well-balanced property under STAMR and improves computational efficiency by focusing refinement in critical regions. To illustrate its practical utility, we simulate pollutant transport in the Pearl River under both flat and non-flat bathymetries. The results show that the method accurately captures pollutant migration using substantially fewer cells and highlight the significant influence of bottom topography on pollutant dispersion.

\section*{Acknowledgments}
This research was supported by the Research Grants Council Areas of Excellence (AoE) Scheme (AoE/P-601/23N-D - MATH),
and by CORE as a joint research center for ocean research between Laoshan Laboratory and HKUST.

%\section*{References}
\bibliographystyle{ieeetr}
\bibliography{STAMR-GKS-SWE}

\begin{thebibliography}{10}

\bibitem{leveque-1998}
R.~J. LeVeque, ``Balancing source terms and flux gradients in high-resolution
  godunov methods: the quasi-steady wave-propagation algorithm,'' {\em Journal
  of computational physics}, vol.~146, no.~1, pp.~346--365, 1998.

\bibitem{zhou2001-surface}
J.~G. Zhou, D.~M. Causon, C.~G. Mingham, and D.~M. Ingram, ``The surface
  gradient method for the treatment of source terms in the shallow-water
  equations,'' {\em Journal of Computational physics}, vol.~168, no.~1,
  pp.~1--25, 2001.

\bibitem{xu2002-swe}
K.~Xu, ``A well-balanced gas-kinetic scheme for the shallow-water equations
  with source terms,'' {\em Journal of Computational Physics}, vol.~178, no.~2,
  pp.~533--562, 2002.

\bibitem{CHO2007}
Y.-S. Cho, D.-H. Sohn, and S.~O. Lee, ``Practical modified scheme of linear
  shallow-water equations for distant propagation of tsunamis,'' {\em Ocean
  Engineering}, vol.~34, no.~11, pp.~1769--1777, 2007.

\bibitem{zhao2021-swe}
F.~Zhao, J.~Gan, and K.~Xu, ``The study of shallow water flow with bottom
  topography by high-order compact gas-kinetic scheme on unstructured mesh,''
  {\em Physics of Fluids}, vol.~33, no.~8, p.~083613, 2021.

\bibitem{zhao2024}
F.~Zhao, J.~Gan, and K.~Xu, ``High-order compact gas-kinetic scheme for
  two-layer shallow water equations on unstructured mesh,'' {\em Journal of
  Computational Physics}, vol.~498, p.~112651, 2024.

\bibitem{delg2024}
A.~{Del Grosso}, M.~J. Castro, A.~Chan, G.~Gallice, R.~Loubère, and P.-H.
  Maire, ``A well-balanced, positive, entropy-stable, and
  multi-dimensional-aware finite volume scheme for 2d shallow-water equations
  with unstructured grids,'' {\em Journal of Computational Physics}, vol.~503,
  p.~112829, 2024.

\bibitem{xing2006-DG}
Y.~Xing and C.-W. Shu, ``High order well-balanced finite volume weno schemes
  and discontinuous galerkin methods for a class of hyperbolic systems with
  source terms,'' {\em Journal of Computational Physics}, vol.~214, no.~2,
  pp.~567--598, 2006.

\bibitem{dambreak-2007}
M.~Ricchiuto, R.~Abgrall, and H.~Deconinck, ``Application of conservative
  residual distribution schemes to the solution of the shallow water equations
  on unstructured meshes,'' {\em Journal of Computational Physics}, vol.~222,
  no.~1, pp.~287--331, 2007.

\bibitem{highorder-efficiency}
D.~Wirasaet, E.~Kubatko, C.~Michoski, S.~Tanaka, J.~Westerink, and C.~Dawson,
  ``Discontinuous galerkin methods with nodal and hybrid modal/nodal
  triangular, quadrilateral, and polygonal elements for nonlinear shallow water
  flow,'' {\em Computer methods in applied mechanics and engineering},
  vol.~270, pp.~113--149, 2014.

\bibitem{Berger1984}
M.~J. Berger and J.~Oliger, ``{Adaptive mesh refinement for hyperbolic partial
  differential equations},'' {\em J. Comput. Phys.}, vol.~53, pp.~484--512,
  1984.

\bibitem{BENK2007}
F.~Benkhaldoun, I.~Elmahi, and M.~Seaı¨d, ``Well-balanced finite volume
  schemes for pollutant transport by shallow water equations on unstructured
  meshes,'' {\em Journal of Computational Physics}, vol.~226, no.~1,
  pp.~180--203, 2007.

\bibitem{XU2024}
C.~Xu, X.~Zhou, H.~Ren, S.~Sutulo, and C.~{Guedes Soares}, ``Real-time
  calculation of ship to ship hydrodynamic interaction in shallow waters with
  adaptive mesh refinement,'' {\em Ocean Engineering}, vol.~295, p.~116943,
  2024.

\bibitem{ZHANG2024}
Z.~Zhang, H.~Tang, and J.~Duan, ``High-order accurate well-balanced energy
  stable finite difference schemes for multi-layer shallow water equations on
  fixed and adaptive moving meshes,'' {\em Journal of Computational Physics},
  vol.~517, p.~113301, 2024.

\bibitem{ZHANG2025}
Z.~Zhang, H.~Tang, and K.~Wu, ``High-order accurate structure-preserving finite
  volume schemes on adaptive moving meshes for shallow water equations:
  Well-balancedness and positivity,'' {\em Journal of Computational Physics},
  vol.~527, p.~113801, 2025.

\bibitem{Osher1983}
S.~Osher and R.~Sanders, ``Numerical approximations to nonlinear conservation
  laws with locally varying time and space grids,'' {\em Math. Comput.},
  vol.~41, no.~164, pp.~321--336, 1983.

\bibitem{Dawson2001}
C.~Dawson and R.~Kirby, ``High resolution schemes for conservation laws with
  locally varying time steps,'' {\em SIAM J. Sci. Comput.}, vol.~22, no.~6,
  pp.~2256--2281, 2001.

\bibitem{xu2}
K.~Xu, ``{A gas-kinetic BGK scheme for the Navier--Stokes equations and its
  connection with artificial dissipation and Godunov method},'' {\em Journal of
  Computational Physics}, vol.~171, no.~1, pp.~289--335, 2001.

\bibitem{CGKSAIA}
F.~Zhao, X.~Ji, W.~Shyy, and K.~Xu, ``Compact higher-order gas-kinetic schemes
  with spectral-like resolution for compressible flow simulations,'' {\em
  Advances in Aerodynamics}, vol.~1, no.~1, pp.~1--34, 2019.

\bibitem{zhao2023direct}
F.~Zhao, X.~Ji, W.~Shyy, and K.~Xu, ``Direct modeling for computational fluid
  dynamics and the construction of high-order compact scheme for compressible
  flow simulations,'' {\em Journal of Computational Physics}, p.~111921, 2023.

\bibitem{XU1997}
K.~Xu, ``Bgk-based scheme for multicomponent flow calculations,'' {\em Journal
  of Computational Physics}, vol.~134, no.~1, pp.~122--133, 1997.

\bibitem{LI2005}
Q.~Li, S.~Fu, and K.~Xu, ``A compressible navier–stokes flow solver with
  scalar transport,'' {\em Journal of Computational Physics}, vol.~204, no.~2,
  pp.~692--714, 2005.

\bibitem{XU_LI_2004}
K.~XU and Z.~LI, ``Microchannel flow in the slip regime: gas-kinetic
  bgk–burnett solutions,'' {\em Journal of Fluid Mechanics}, vol.~513,
  p.~87–110, 2004.

\bibitem{TAN2018}
S.~Tan, Q.~Li, Z.~Xiao, and S.~Fu, ``Gas kinetic scheme for turbulence
  simulation,'' {\em Aerospace Science and Technology}, vol.~78, pp.~214--227,
  2018.

\bibitem{P4est2011}
C.~Bursteddeand, L.~Wilcox, and O.~Ghattas, ``{{\texttt{p4est}}: scalable
  algorithms for parallel adaptive mesh refinement on forests of octrees},''
  {\em SIAM J. Sci. Comput.}, vol.~33(3), pp.~1103--1133, 2011.

\bibitem{venkatakrishnan1995}
V.~Venkatakrishnan, ``Convergence to steady state solutions of the euler
  equations on unstructured grids with limiters,'' {\em Journal of
  computational physics}, vol.~118, no.~1, pp.~120--130, 1995.

\bibitem{Briggs1995}
S.~C. E. H. G. S. G. D.~R. Briggs, Michael~J., ``Laboratory experiments of
  tsunami runup on a circular island,'' {\em pure and applied geophysics},
  vol.~144, no.~3, pp.~569--593, 1995.

\bibitem{Garc2019}
M.~J. F.-P.~J. Garc{\'i}a-Navarro, P., ``The shallow water equations and their
  application to realistic cases,'' {\em Environmental Fluid Mechanics},
  vol.~19, no.~5, pp.~1235--1252, 2019.

\end{thebibliography}

\end{document}